\documentclass[11pt,a4paper,reqno]{amsart}
\usepackage[a4paper,top=3cm, bottom=3cm, left=2.8cm, right=2.8cm]{geometry}

\usepackage[american]{babel} 
\usepackage[utf8]{inputenc} 

\usepackage{url}
\usepackage[backend=biber, style=numeric, giveninits=true, sorting=nyt, doi=true, date=year, maxbibnames=99,maxcitenames=99]{biblatex}
\bibliography{ref.bib}

\usepackage{csquotes}

\usepackage{nth}

\usepackage{xcolor}
\usepackage[shortlabels]{enumitem}
\usepackage{hyperref}\hypersetup{colorlinks,linkcolor=.,citecolor=.,urlcolor=.}

\usepackage{booktabs}

\usepackage{amsmath,amssymb,amsthm,mathtools}

\newtheorem{thm}{Theorem}
\numberwithin{thm}{section}
\newtheorem{lem}[thm]{Lemma}
\newtheorem{prop}[thm]{Proposition}
\newtheorem{cor}[thm]{Corollary}

\theoremstyle{remark}
\newtheorem{rem}[thm]{Remark}

\newcommand{\thistheoremname}{}
\newtheorem*{genericthm*}{\thistheoremname}

\newcommand{\RR}{\mathbb{R}}
\newcommand{\FF}{\mathbb{F}}
\newcommand{\cP}{\mathcal{P}}

\newcommand{\hamming}{H(n,q)}
\newcommand{\johnson}{J(n,m)}
\newcommand{\qjohnson}{J_q(n,m)}
\newcommand{\bil}{\operatorname{Bil}_q(n,m)}
\newcommand{\alt}{\operatorname{Alt}_q(m)}
\newcommand{\her}{\operatorname{Her}_q(n)}

\newcommand{\hermitianpol}[1]{\prescript{2}{}{A}_{{#1}}}
\newcommand{\hyperbolichalf}{\tfrac{1}{2}D_m}

\newcommand{\elliptic}[1]{\prescript{2}{}{D}_{{#1}}}

\newcommand{\qbin}[2]{\genfrac{[}{]}{0pt}{}{{#1}}{{#2}}_q}
\newcommand{\bbin}[2]{\genfrac{[}{]}{0pt}{}{{#1}}{{#2}}_b}
\newcommand{\pbin}[2]{\genfrac{[}{]}{0pt}{}{{#1}}{{#2}}_p}
\newcommand{\qqbin}[2]{\genfrac{[}{]}{0pt}{}{{#1}}{{#2}}_{q^2}}

\newcommand{\poch}[3]{({#1};{#2})_{{#3}}}
\newcommand{\pochq}[2]{\poch{{#1}}{q}{{#2}}}
\newcommand{\pochb}[2]{({#1})_{{#2}}}
\newcommand{\pochp}[2]{\poch{{#1}}{p}{{#2}}}

\newcommand{\hyp}[3]{\,\mbox{}_{3}\phi_{2}
	\left(\!\!\left.\begin{array}{c}#1\\#2\end{array}\right|#3\right)}
\newcommand{\hyper}[5]{\,\mbox{}_{{#1}}\phi_{{#2}}
	\left(\!\!\left.\begin{array}{c}#3\\#4\end{array}\right|#5\right)}

\DeclareMathOperator{\rk}{rank}
\DeclareMathOperator{\LP}{LP}
\DeclareMathOperator{\GL}{GL}
\DeclareMathOperator{\sign}{sign}

\renewcommand{\epsilon}{\varepsilon}

\newcommand{\qandq}{\quad\text{and}\quad}

\makeatletter
\@namedef{subjclassname@2020}{%
	\textup{2020} Mathematics Subject Classification}
\makeatother

\title[The LP optimum for packings in classical association schemes]{The linear programming optimum for packings in classical association schemes}

\author{Kai-Uwe Schmidt}
\address[K.-U. Schmidt]{Department of Mathematics, Paderborn University, Warburger Str.\ 100, 33098 Paderborn, Germany.}

\author{Charlene Weiss}
\address[C. Wei\ss]{Department of Mathematics, Paderborn University, Germany, and, 
Korteweg-de Vries Institute for Mathematics, University of Amsterdam, The Netherlands.}
\email{chweiss@math.upb.de}

\thanks{The second author was partially supported during this research project by the Deutsche Forschungsgemeinschaft (DFG, German Research Foundation)--Project-ID 491392403--TRR 358, as well as by a fellowship of the German Academic Exchange Service (DAAD)}

\date{18 August 2025, updated 27 March 2026}

\subjclass[2020]{Primary 05E30; Secondary 94B65, 90C05, 51E23, 33C80, 05D99} 

\keywords{}

\begin{document}
	\begin{abstract}
		Association schemes are central objects in algebraic combinatorics, with the classical schemes lying at their core. These classical association schemes essentially consist of the Hamming and Johnson schemes, and their $q$-analogs: bilinear forms scheme, alternating bilinear forms scheme, Hermitian forms scheme, $q$-Johnson scheme, and polar space schemes. Each of them gives rise to a distance-regular graph on a vertex set $X$, naturally endowed with the path metric. We study $d$-codes in these schemes, that is, subsets $Y$ of $X$ in which every pair of distinct elements has path distance at least $d$. A powerful tool for deriving upper bounds on the size of $d$-codes is the linear programming method. In the case of the Hamming and Johnson schemes, the linear program has been studied since the 1970s, but its optimum is still unknown. We determine the optimum of the linear program for nearly all classical association schemes distinct from the Hamming and Johnson schemes. As a corollary, we obtain upper bounds on $t$-intersecting sets in classical association schemes, providing new proofs of several known results and, in particular, improving earlier bounds on $t$-intersecting sets of generators in polar spaces. These results can be viewed as analogs of the classical Erd\H{o}s--Ko--Rado Theorem in extremal set theory. Our proofs draw on techniques from algebraic combinatorics and the duality theory of linear programming.
	\end{abstract}
	
	\maketitle
	
	\section{Introduction}
	The linear programming method has proven to be a very mighty tool for solving extremal problems across various areas of mathematics, for example see \cite{DelsartePhD}, \cite{MRRW}, \cite{DelGoeSei}, \cite{CohnElkies}, \cite{Viazovska}. Originally, Delsarte \cite{DelsartePhD} established this method to study packing problems in association schemes and used it to obtain upper bounds on packings in several classical association schemes such as on the error-correcting codes in the Hamming scheme and on the binary constant-weight codes in the Johnson scheme. The classical association schemes are distance-regular graphs with classical parameters and give naturally rise to a finite metric space $(X,\rho)$, where $X$ is the vertex set and $\rho$ is the path metric of the graph. A subset $Y$ of $X$ is then called a \emph{$d$-code} if the distance $\rho(x,y)$ between all distinct elements $x,y\in Y$ is at least $d$. The central problem is to determine, for fixed $d$, the maximum size of a $d$-code. One is thus interested in deriving upper bounds on the size of codes. To this day, almost all known bounds on the size of codes in classical association schemes can be derived from Delsarte's linear program (see, e.g., \cite[Ch.~17, \S~4]{MacWilliamsSloane}, \cite{DelsarteBilinear}, \cite{DelGoe1975}, \cite{SchmidtHermitian}, \cite{ZhangLPqJohn}, \cite{SchmidtWeissSteiner}). It is therefore a natural problem to determine the optimum of this linear program. For the Hamming and Johnson schemes, this is a classical problem, which has been studied since the 1970s, but is still unsolved. In this paper, we compute the optimum of the linear program in most classical association schemes, distinct from the Hamming and Johnson schemes.

	Namely, we study several $q$-analogs of the Hamming and Johnson schemes, which, together with these two schemes, essentially comprise the list of all classical association schemes. These $q$-analogs can be split into two groups---the \emph{ordinary $q$-analogs} and the \emph{affine $q$-analogs}. For the former ones, the vertex set $X$ of the distance-regular graph consists of all $n$-dimensional subspaces in a finite vector space $V$ of type $A_{m+n-1}$, $\hermitianpol{2n-1}$, $\hermitianpol{2n}$, $B_n$, $C_n$, $D_n$, or $\elliptic{n+1}$. The notation comes from the type of the group acting on the vector space $V$, namely from the groups $\GL(m+n,q)$, $U(2n,q^2)$, $U(2n+1,q^2)$, $O(2n+ 1,q)$, $Sp(2n,q)$, $O^+(2n,q)$, or $O^-(2n+2,q)$, respectively, see \cite{Carter}. In the case of $A_{m+n-1}$, one obtains the \emph{$q$-Johnson scheme} (also known as \emph{Grassmann scheme}). In the other cases, the space $V$ is equipped with a nondegenerate form and the notation is chosen such that the maximal totally isotropic subspaces of $V$ have dimension $n$ (see Table~\ref{table:polarspaces}). A collection of all totally isotropic subspaces with respect to a given form is a finite classical polar space (or just polar space) of rank $n$ and the arising association schemes are thus called \emph{polar space schemes}. We denote the polar spaces and their schemes by the same symbol as the type of the underlying vector space.
	The $q$-Johnson scheme and the polar space schemes have the same metric, namely $(x,y)\mapsto n-\dim(x\cap y)$ (which agrees with the subspace metric used by coding theorists).
	
	Then there are the affine $q$-analogs, where the set $X$ consists of all unrestricted $m\times n$ matrices over a finite field $\FF_q$, all Hermitian $n\times n$ matrices over $\FF_{q^2}$, or all alternating matrices $m\times m$ matrices over $\FF_q$. These schemes are known as the \emph{bilinear forms scheme}, \emph{alternating bilinear forms scheme}, or \emph{Hermitian forms scheme}, respectively. In all three of them, the metric is given by $(x,y)\mapsto \rk(x-y)$ (which agrees with the rank metric used in coding theory). 

	Our first main result, Theorem~\ref{thm:LPopt_matrix_qJohn_herm_hyp}, contains the linear programming optimum for codes in all affine $q$-analogs and their $q$-ordinary counterparts, namely the $q$-Johnson scheme, the polar space scheme of type $\hermitianpol{2n-1}$, and the bipartite half of the polar space scheme of type $D_n$ (see Section~\ref{sec:codesAS} for background on this association scheme). For both groups, the linear programming optimum is expressed in a unified way. Our second main result, Corollary~\ref{cor:LPopt_Bn_Cn_Dn}, gives the linear programming optimum in the polar space schemes of type $B_n$, $C_n$, and $D_n$ for most parameter sets. It turns out that all the derived optima coincide with the bounds in \cite{DelsarteBilinear}, \cite{DelGoe1975}, \cite{SchmidtHermitian}, \cite{WangXingSafavi-Naini}, \cite{SchmidtWeissSteiner}.

	As a byproduct, we obtain Erd\H{o}s--Ko--Rado-type bounds in the classical association schemes in Corollary \ref{cor:EKR}, which immediately follow from the computed linear programming optima. These bounds are analogs of the famous result on $t$-intersecting sets given in \cite{EKR}, which in its generalized form states that for all fixed $n$ and $t$ and all sufficiently large $v$, every $t$-intersecting set of $n$-subsets of $\{1,2,\dots,v\}$ has size at most $\binom{v-t}{n-t}$. Several analogs of \cite{EKR} are known (see \cite{GodsilMeagher} for an overview), in particular in the classical association schemes (see, e.g., \cite{FranklWilson}, \cite{Moon}, \cite{Huang}, \cite{StantonEKR}, \cite{IhringerMetsch}). However, we give a new proof of these bounds and moreover, we improve the so-far best known upper bounds on $t$-intersecting sets from \cite{IhringerMetsch} in the polar spaces $\hermitianpol{2n-1}$, $B_n$, $C_n$, and $D_n$ for infinitely many cases.

	We organize the paper as follows. In Section~\ref{sec:codesAS}, we recall relevant background on classical association schemes and the linear programming method. These are required to state our main results in Section~\ref{sec:mainres}, where we first give the optimum of the linear programs and then use them to establish upper bounds on $t$-intersecting sets in the classical association schemes. In Section~\ref{sec:PQnum}, we look at the orthogonal polynomials associated to the classical association schemes, which will play a central role in computing the linear programming optima. In Section~\ref{sec:linearprogramming}, we give some facts from optimization theory. The proof of the main result, Theorem~\ref{thm:LPopt_matrix_qJohn_herm_hyp}, is then given in the remaining sections. First, we compute a feasible solution of the dual linear program in Section~\ref{sec:SolDualLP} and then of the primal linear program in Section~\ref{sec:SolPrimalLP}.

	\section{Codes in classical association schemes}\label{sec:codesAS}
	In this section, we give some facts about association schemes, codes, and the linear programming method that are needed to state our main results in Section~\ref{sec:mainres}.

	A (symmetric) \emph{association scheme $(X,(R_i))$ with $n$ classes} is a finite set $X$ together with $n+1$ nonempty relations $R_0,R_1,\dots,R_n$ such that
	\begin{itemize}
		\item $R_0$ is the identity relation and all $n+1$ relations partition $X\times X$,
		\item each relation is symmetric, that is, if $(x,y)\in R_i$, then $(y,x)\in R_i$,
		\item for every pair $(x,y)\in R_k$, the number of $z\in X$ with $(x,z)\in R_i$ and $(z,y)\in R_j$ is a constant $p_{ij}^k$ depending only on $i$, $j$, and $k$, but not on the particular choice of $(x,y)$.
	\end{itemize}

	We will focus on \emph{classical} association schemes, which are distance-regular graphs with classical parameters (cf.\ \cite[Table~6.1]{BrouwerCohenNeumaier}). The list of classical association schemes essentially consists of the Hamming scheme, the Johnson scheme, the \emph{ordinary $q$-analogs} of the Hamming or Johnson scheme, and the \emph{affine $q$-analogs} of the Hamming scheme, which are defined below. For more information on the classical association schemes, we refer to \cite[\S~3.6]{BannaiIto1984}, \cite[\S~6.4.1]{Bannaietal2021}, and \cite[\S~9]{BrouwerCohenNeumaier}.
	
	\subsubsection*{Hamming scheme}
	Let $n$ and $q$ be positive integers with $q\geq 2$ and let $X=Q^n$, where $Q$ is a set with $q$ elements. The \emph{Hamming distance} $d_H(x,y)$ between two points $x,y\in X$ is the number of positions in which $x$ and $y$ differ. 
	Then $X$ together with the relations
	\begin{align*}
	R_i=\{(x,y)\in X\times X : d_H(x,y)=i\}
	\end{align*}
	for $i=0,1,\dots,n$ forms an association scheme with $n$ classes, known as the \emph{Hamming scheme} and denoted by $\hamming$. We have $|X|=q^n$.
	
	\subsubsection*{Johnson scheme}
	Let $m$ and $n$ be positive integers with $m\geq n$. Then the set $X$ of all $n$-subsets of a given $(m+n)$-set together with the relations
	\begin{align*}
		R_i=\{(x,y)\in X\times X : |x\cap y|=n-i\}
	\end{align*}
	for $i=0,1,\dots,n$ forms an association scheme with $n$ classes, known as the \emph{Johnson scheme} and denoted by $\johnson$. We have $|X|=\binom{m+n}{n}$.
	
	\vspace{\baselineskip}The affine $q$-analogs are the following three association schemes.
	
	\subsubsection*{Bilinear forms scheme}
	Let $m$ and $n$ be integers with $m\geq n\geq 1$ and let $X$ be the set of $m\times n$ matrices over the finite field~$\FF_q$. Then $X$ together with the relations
	\begin{align}\label{eq:Ri_bil}
	R_i=\{(A,B)\in X\times X : \rk(A-B)=i\}
	\end{align}
	for $i=0,1,\dots,n$ forms an association scheme with $n$ classes, known as the \emph{bilinear forms scheme} and denoted by $\bil$. We have $|X|=q^{mn}$. Here, we also refer to \cite{DelsarteBilinear}.
	
	\subsubsection*{Alternating bilinear forms scheme}
	An $m\times m$ matrix $A=(a_{ij})$ over $\FF_q$ is called \emph{alternating} if it is skew-symmetric with zero main diagonal, that is, $a_{i,i}=0$ and $a_{i,j}=-a_{j,i}$ for all $i,j$. Let $X$ be the set of alternating $m\times m$ matrices over $\FF_q$. Then $X$ together with the relations
	\begin{align}\label{eq:Ri_alt}
	R_i=\{(A,B)\in X\times X : \rk(A-B)=2i\}
	\end{align}
	for $i=0,1,\dots,n$ with $n=\lfloor m/2 \rfloor$ is an association scheme with $n$ classes, known as the \emph{alternating bilinear forms scheme} and denoted by $\alt$. We have $|X|=q^{\binom{m}{2}}$. We also refer to \cite{DelGoe1975}.
	
	\subsubsection*{Hermitian forms scheme}
	An $n\times n$ matrix~$A$ over $\FF_{q^2}$ is called \emph{Hermitian} if $A^\ast=A$, where~$A^\ast$ is obtained from~$A$ by conjugation $x\mapsto x^q$ of each entry~$x$ of $A$ and transposition. Let~$X$ be the set of $n\times n$ Hermitian matrices over $\FF_{q^2}$. Then $X$ together with the relations
	\begin{align}\label{eq:Ri_her}
	R_i=\{(A,B)\in X\times X : \rk(A-B)=i\}
	\end{align}
	for $i=0,1,\dots,n$ gives rise to an association scheme with $n$ classes, known as the \emph{Hermitian forms scheme} and denoted by $\her$. We have $|X|=q^{n^2}$. See also~\cite{SchmidtHermitian} for more information.

	\vspace{\baselineskip}We now come to the ordinary $q$-analogs.
	
	\subsubsection*{$q$-Johnson scheme}
	Let $m$ and $n$ be positive integers with $m\geq n$. Let~$X$ be the set of all $n$-dimensional subspaces ($n$-spaces for short) of an $(m+n)$-dimensional vector space over $\FF_q$. Then $X$ together with the relations
	\begin{align}\label{eq:Ri_qJohnson}
	R_i=\{(U,W)\in X\times X : \dim(U\cap W)=n-i\}
	\end{align}
	for $i=0,1,\dots,n$ forms an association scheme with $n$ classes, known as the \emph{$q$-Johnson scheme} (also called \emph{Grassmann scheme}) and denoted by $\qjohnson$. We have $|X|=\qbin{m+n}{n}$, where $\qbin{n}{k}$ is the \emph{$q$-binomial coefficient} defined by
	\[
	\qbin{n}{k}=\prod_{j=1}^k \frac{q^{n-j+1}-1}{q^j-1}
	\]
	for nonnegative integers $n$ and $k$.
	
	\subsubsection*{Polar space schemes}
	Let~$V$ be a vector space over a finite field~$\FF_p$ with~$p$ elements  equipped with a nondegenerate form~$f$. A subspace $U$ of $V$ is called \emph{totally isotropic} if \mbox{$f(u,w)=0$} for all $u,w\in U$, or in the case of a quadratic form, if $f(u)=0$ for all $u\in U$. A \emph{finite classical polar space} with respect to a form~$f$ consists of all totally isotropic subspaces of~$V$. It is well known that all maximal (with respect to the dimension) totally isotropic spaces in a polar space have the same dimension, which is called the \emph{rank} of the polar space. The maximal totally isotropic spaces are called \emph{generators}. A finite classical polar space~$\cP$ has the \emph{parameter}~$e$ if every $(n-1)$-space in~$\cP$ lies in exactly $p^{e+1}+1$ generators. Up to isomorphism, there are exactly six finite classical polar spaces of rank~$n$, which are listed together with their parameter~$e$ in Table~\ref{table:polarspaces}. We refer to~\cite{Cam1992},~\cite{Taylor},~\cite[\S~9.4]{BrouwerCohenNeumaier},~\cite[\S~4.2]{Bal2015}, and~\cite[\S~5.1]{HirschfeldThas} for further background on polar spaces. 
	\begin{table}[ht]
		\caption{List of all six finite classical polar spaces.}
		\centering
		\renewcommand*{\arraystretch}{1.25}
		\begin{tabular}{cccccccc}
			\toprule[0.4mm]
			name & form & type & group & $\dim(V)$ & $p$ & $e$\\
			\midrule[0.4mm]
			Hermitian & Hermitian & $\hermitianpol{2n-1}$ & $U(2n,q^2)$ & $2n$ & $q^2$ & $-1/2$\\ 
			Hermitian & Hermitian & $\prescript{2}{}{A}_{2n}$ & $U(2n+1,q^2)$ & $2n+1$ & $q^2$ & $1/2$\\
			symplectic & alternating & $C_n$ & $Sp(2n,q)$ & $2n$ & $q$ &  $0$\\
			hyperbolic & quadratic & $D_n$& $O^+(2n,q)$ & $2n$ & $q$ & $-1$\\
			parabolic & quadratic & $B_n$ & $O(2n+1,q)$ & $2n+1$ & $q$ & $0$\\
			elliptic & quadratic & $\prescript{2}{}{D}_{n+1}$ & $O^-(2n+2,q)$ & $2n+2$ & $q$ & $1$\\
			\bottomrule[0.4mm]
		\end{tabular}
		\label{table:polarspaces}
	\end{table}

	Let $X$ consist of all generators of a polar space of rank~$n$ and define the relations
	\begin{equation}\label{eq:Ri_polarspaces}
		R_i=\{(U,W)\in X\times X:\dim(U\cap W)=n-i\}
	\end{equation}
	for $i=0,1,\dots,n$. Then $(X,(R_i))$ is an association scheme with $n$ classes. It is well known that
	\begin{align*}
		|X|=\prod_{i=1}^n (1+p^{i+e}).
	\end{align*}
	
	We note that~$D_n$ gives rise to another association scheme, called the \emph{bipartite half of~$D_n$}, in the following way. Let $X$ be the set of generators of $D_n$ and define two generators in $X$ to be equivalent if the dimension of their intersection has the same parity as $n$. This induces two equivalence classes, $X_1$ and $X_2$, and each pair $(X_i,(R_{2j})_{0\le j\le \lfloor\frac n2\rfloor})$ is an association scheme with $\lfloor n/2\rfloor$ classes~\cite[\S~9.4.C]{BrouwerCohenNeumaier}, denoted by~$\frac 12 D_n$.
	
	The affine $q$-analogs $\bil$, $\her$, and $\alt$ are related to the ordinary $q$-analogs $\qjohnson$, $\hermitianpol{2n-1}$, and $\hyperbolichalf$ in the following way. For a vector space~$V$ over $\FF_q$, let $P_n(V)$ be the set of $n$-spaces of~$V$. Define the mapping
	\begin{equation*}
			v\colon\FF_q^{m\times n}\to P_n(\FF_q^{m+n}),\qquad			A\mapsto\left\{\binom{x}{Ax}: x\in\FF_q^n\right\}.
	\end{equation*}
	It is well known~\cite[\S~9.5.E]{BrouwerCohenNeumaier} that, for $m=n$, after an appropriate choice of the form, $v(A)$ is in $\hermitianpol{2n-1}$ if and only if $A$ is Hermitian (here, $q$ has to be a square), and $v(A)$ is in $D_n$ if and only if~$A$ is alternating. The mapping~$v$ satisfies
	\[
	n-\dim(v(A)\cap v(B))=\rk(A-B)
	\]
	for all $A,B\in\FF_q^{m\times n}$, so in particular $v$ is injective. We therefore obtain the following embeddings:
	\begin{equation}\label{eq:embeddings_AS}
		\begin{aligned}
			\bil&\hookrightarrow \qjohnson\\
			\her&\hookrightarrow \hermitianpol{2n-1}\\
			\alt &\hookrightarrow \hyperbolichalf
		\end{aligned}
	\end{equation}
	
	We now give some well-known facts about association schemes and then introduce codes in association schemes.

	Let $(X,(R_i))$ be an association scheme with $n$ classes. For each relation $R_i$, the adjacency matrix of the graph $(X,R_i)$ is denoted by $D_i$ (with $D_0$ being the identity matrix). The complex vector space spanned by the matrices $D_0,D_1,\dots,D_n$ is a commutative matrix algebra of dimension $n+1$, called the \emph{Bose-Mesner algebra} of the association scheme. There exists a unique basis of minimal idempotent matrices $E_0(=1/|X|J),E_1,\dots,E_n$ for this algebra, where $J$ denotes the all-ones matrix. A change of basis gives the existence of unique real numbers $P_i(k)$ and $Q_k(i)$ such that
	\[
	D_i=\sum_{k=0}^n P_i(k) E_k \qandq E_k=\frac{1}{|X|}\sum_{i=0}^n Q_k(i)D_i.
	\]
	The numbers $P_i(k)$ and $Q_k(i)$ are called \emph{$P$-numbers} and \emph{$Q$-numbers} of the association scheme $(X,(R_i))$, respectively. Write $v_i=P_i(0)$ and $\mu_k=Q_k(0)$, which are called \emph{valencies} and \emph{multiplicities} of the association scheme, respectively. Indeed $P_i(k)$ is an eigenvalue of the graph $(X,R_i)$, each column of $E_k$ is a corresponding eigenvector, and the rank of $E_k$ equals $\mu_k$. In particular $v_i$ is the valency of the graph $(X,R_i)$. 
		
	An association scheme is called \emph{$P$-polynomial} with respect to the ordering $R_0,R_1, \dots,R_n$ if there exist polynomials $f_i\in\RR[x]$ of degree $i$ and distinct real numbers $y_0, y_1,\dots,y_n$ such that $P_i(k)=f_i(y_k)$ for all $i,k=0,1,\dots,n$. Accordingly, an association scheme is called \emph{$Q$-polynomial} with respect to the ordering $E_0,E_1,\dots,E_n$ if there exist polynomials $g_k\in\RR[x]$ of degree $k$ and distinct real numbers $z_0,z_1,\dots,z_n$ such that $Q_k(i)=g_k(z_i)$ for all $i,k=0,1,\dots,n$. We will later see that all the previous introduced ordinary and affine $q$-analogs are ($P$ and $Q$)-polynomial.
	
	A \emph{$D$-code} is a subset $Y$ of~$X$ such that $A_i=0$ for all $i\in\{1,2,\dots,n\}\setminus D$, where~$D$ is a subset of $\{1,\dots,n\}$. If the association scheme $(X,(R_i))$ is $P$-polynomial with respect to the ordering $R_0,R_1,\dots,R_n$, then a subset $Y$ of $X$ is a \emph{$d$-code} if no pair $(x,y)\in Y\times Y$ lies in one of the relations $R_1,R_2,\dots,R_{d-1}$, that means, $Y$ is a $D$-code with $D=\{d,d+1,\dots,n\}$. A central coding-theoretic problem is to obtain upper bounds on the size of $D$-codes. One way to achieve this is by combining the theory of association schemes with linear programming. To do so, we associate with each subset~$Y$ of $X$ two tuples of rational numbers. The \emph{inner distribution} of~$Y$ is the tuple $(A_0,A_1,\dots,A_n)$, where
	\[
	A_i=\frac{|(Y\times Y)\cap R_i|}{|Y|}.
	\]
	We then have $A_0=1$ and $\sum_{i=0}^nA_i=|Y|$. Note that $A_i=0$ for all $i\in\{1,\dots,n\}\setminus D$ if and only if~$Y$ is a $D$-code. The \emph{dual distribution} of $Y$ is the tuple $(A'_0,A'_1,\dots,A'_n)$, where
	\begin{align}\label{eq:dualdistr}
		A'_k=\sum_{i=0}^n Q_k(i)A_i.
	\end{align}
	Since $Q_0(i)=1$, we obtain $A'_0=|Y|$. Moreover, we have $A'_k\ge 0$ for all $k=0,1,\dots,n$, see~\cite[Thm.~3.3]{DelsartePhD}, for example. We call $Y$ a \emph{$T$-design} if $A'_k=0$ for all $k\in T$, where $T$ is a subset of $\{1,2,\dots,n\}$. If the association scheme $(X,(R_i))$ is $Q$-polynomial with respect to an ordering $E_0,E_1,\dots,E_n$ of the idempotents, then $Y$ is a \emph{$t$-design} if $A'_k=0$ for all $k=1,2,\dots,t$, that means, $Y$ is a $T$-design with $T=\{1,2,\dots,t\}$.
	
	Delsarte combined all properties of the distributions and established the following \emph{linear program~(LP)}:
	\begin{align}\label{eq:LPprimalCodes}
		\begin{array}{lrcl}
			\operatorname*{maximize}\limits_{x_i\in\RR} & \sum\limits_{i=0}^n x_i & &\\
			\text{subject to} & x_0&=& 1\\
			& x_i&\geq & 0\quad\text{for all $i\in D$}\\
			&x_i&=&0\quad\text{for all $i\in\{1,\dots,n\}\setminus D$}\\
			& \sum\limits_{i=0}^n Q_k(i)x_i &\geq & 0\quad\text{for all $k=1,2,\dots,n$.}
		\end{array}
	\end{align}
	If a vector $x\in\RR^{n+1}$ satisfies all the constraints of the linear program~\eqref{eq:LPprimalCodes}, then $x$ is called a \emph{feasible solution} of~\eqref{eq:LPprimalCodes}. A feasible solution $x^\ast$ of~\eqref{eq:LPprimalCodes} is called \emph{optimal} if $\sum_{i=0}^n x_i\leq \sum_{i=0}^n x^\ast_i$ for all feasible solutions $x$. The objective function value $\sum_{i=0}^n x^\ast_i$ of an optimal solution is then called the \emph{LP optimum} of~\eqref{eq:LPprimalCodes} and denoted by $\LP(D)$. In the case of a $P$-polynomial association scheme and $D=\{d,d+1,\dots,n\}$, the LP optimum is also denoted by $\LP(d)$. The linear program~\eqref{eq:LPprimalCodes} gives upper bounds on codes. Namely, if $Y$ is a $D$-code in $(X,(R_i))$, then the inner distribution of $Y$ satisfies the constraints of~\eqref{eq:LPprimalCodes} and in particular, we have $|Y|\leq \LP(D)$ (see \cite[\S~3.2, Thm.~3.8]{DelsartePhD}).

	\section{Main results}\label{sec:mainres}
	In this section, we will state our main results---first the LP optima for the classical association schemes, and second, the implied upper bounds on $t$-intersecting sets in classical association schemes.

	In accordance with the embeddings~\eqref{eq:embeddings_AS} of the classical association schemes, we set
	\begin{align}\label{eq:bc}
		(b,c)=\begin{cases}
			(q,q^{m-n})&\text{for $\bil$ and $\qjohnson$}\\
			(-q,-1) &\text{for $\her$ and $\hermitianpol{2n-1}$}\\
			(q^2,1/q) &\text{for $\alt$ and $\hyperbolichalf$ if $m$ is even}\\
			(q^2,q) &\text{for $\alt$ and $\hyperbolichalf$ if $m$ is odd}.
		\end{cases}
	\end{align}
	We can now give our first main result.
	\begin{thm}\label{thm:LPopt_matrix_qJohn_herm_hyp}
		~\begin{enumerate}[(a)]
			\item\label{LPoptJohnson} The LP optimum for $d$-codes with $1\leq d\leq n$ in $\johnson$ is given by
			\[
			\LP(d)=\frac{\binom{m+n}{n-d+1}}{\binom{n}{n-d+1}}
			\]
			if $m$ is sufficiently large and $\binom{n-i}{d-1}$ divides $\binom{m+n-i}{m+d-1}$ for all $i=0,1,\dots,n-d$.
			\item\label{LPoptHamming} The LP optimum for $d$-codes with $1\leq d\leq n$ in $\hamming$ is given by
			\[
			\LP(d)=q^{n-d+1}
			\]
			if $q\geq \max\{d,n-d+2\}$.
			\item\label{LPoptordinaryq} Let $X$ be the set of $n$-spaces in $\qjohnson$, or generators in $\hermitianpol{2n-1}$ or $\hyperbolichalf$. Then the LP optimum for $d$-codes with $1\leq d\leq n$ in $\qjohnson$, $\hermitianpol{2n-1}$, and $\hyperbolichalf$ is given by
			\begin{align}\label{eq:LPbound_qJohn_herm_hyp}
				\LP(d)=|X|\prod\limits_{i=0}^{d-2}\frac{qb^i-1}{qcb^{n+i}-1},
			\end{align}
			where $d$ is required to be odd in the case of $\hermitianpol{2n-1}$ and $n=\lfloor m/2\rfloor$ in the case of $\hyperbolichalf$. For even $d$ with $2\leq d\leq n$, the LP optimum for $d$-codes in $\hermitianpol{2n-1}$ is given by
			\begin{align}\label{eq:LPbound_hermpol_deven}
				\LP(d)= |X|\left(\prod\limits_{i=0}^{d-2}\frac{qb^i-1}{qcb^{n+i}-1}\right)\epsilon(n,d),
			\end{align}
			where
			\begin{align}\label{eq:epsnd_Hermpol_deven}
				\epsilon(n,d)=\frac{((-q)^{n-d+2}-1)+q\frac{(-q)^{n+d-2}-1}{q(-q)^{d-2}-1}((-q)^{n-d+1}-1)}{((-q)^{n-d+2}-1)+q\frac{(-q)^{n+d-2}-1}{(-q)^{n+d-1}-1}((-q)^{n-d+1}-1)}.
			\end{align}
			Moreover, both LP optima~\eqref{eq:LPbound_qJohn_herm_hyp} and~\eqref{eq:LPbound_hermpol_deven} also hold for association schemes with the same $P$- and $Q$-numbers as $\qjohnson$, $\hermitianpol{2n-1}$, and $\hyperbolichalf$.
			\item\label{LPoptaffineq} The LP optimum for $d$-codes with $1\leq d\leq n$ in $\bil$, $\her$, and $\alt$ is given by
			\begin{align}\label{eq:LPbound_matrixschemes}
				\LP(d)=(cb^n)^{n-d+1},
			\end{align}
			where $d$ is required to be odd in the case of $\her$ and $n=\lfloor m/2\rfloor$ in the case of $\alt$. For even $d$ with $2\leq d\leq n$, the LP optimum for $d$-codes in $\her$ is given by
			\begin{align}\label{eq:LPbound_Hermat_deven}
				\LP(d)=(cb^n)^{n-d+1}\,\frac{(b^{n-d+2}-1)+b^n(b^{n-d+1}-1)}{b^{n-d+2}-b^{n-d+1}}.
			\end{align}
			Moreover, both LP optima~\eqref{eq:LPbound_matrixschemes} and~\eqref{eq:LPbound_Hermat_deven} also hold for association schemes with the same $P$- and $Q$-numbers as $\bil$, $\her$, and $\alt$.
		\end{enumerate}
	\end{thm}

	Part \ref{LPoptJohnson} and \ref{LPoptHamming} of Theorem~\ref{thm:LPopt_matrix_qJohn_herm_hyp} are essentially known, but added for completeness. Namely, in \cite[\S~4.3.2]{DelsartePhD}, Delsarte derived a feasible solution of the dual linear program (see \eqref{eq:LPdualCodes}) for $d$-codes in $\johnson$ and $\hamming$ with objective function value $\binom{m+n}{n-d+1}/\binom{n}{n-d+1}$ and $q^{n-d+1}$, respectively. In the case of $\johnson$, it was shown in \cite{Kee2014} and \cite{GloKuhLoOst2016} that for all $n$ and $d$, $(n-d+1)$-$(m+n,n,1)$ designs exist if $m$ is sufficiently large enough and if the basic divisibility conditions are satisfied. These designs are $d$-codes of size $\binom{m+n}{n-d+1}/\binom{n}{n-d+1}$, and therefore, their inner distribution is a feasible solution of the primal LP~\eqref{eq:LPprimalCodes} for $d$-codes in $\johnson$. Using the well-known strong duality theorem from optimization theory (see Theorem \ref{thm:strongduality}), Theorem~\ref{thm:LPopt_matrix_qJohn_herm_hyp}~\ref{LPoptJohnson} follows. In the case of $\hamming$, according to \cite[\S~4.3.2]{DelsartePhD}, Piret constructed a feasible solution with objective function value $q^{n-d+1}$ of the primal LP~\eqref{eq:LPprimalCodes} for $d$-codes in $\hamming$ if $q\geq \max\{d,n-d+2\}$. Hence, the strong duality theorem implies Theorem~\ref{thm:LPopt_matrix_qJohn_herm_hyp}~\ref{LPoptHamming}. We note that in the case of the Hamming scheme with $q=2$, it is still unknown what the optimum of the linear program is. However, numerical evidence \cite{BargJaffe} indicates that the bounds from \cite{MRRW} are optimal on the exponential scale. 

	Theorem~\ref{thm:LPopt_matrix_qJohn_herm_hyp} also gives the LP optimum in $B_n$, $C_n$, and $D_n$ in the following way. First, $B_n$ and $C_n$ induce new association schemes with the classes $R_0,R_1\cup R_2, R_3\cup R_4,\dots$, which have the same $P$- and $Q$-numbers as $\frac12D_{n+1}$ \cite{Ivanovetal}. Therefore, without loss of generality for odd $d$, we can study the linear program for $\frac{d+1}{2}$-codes in $\frac12D_{n+1}$ instead of the linear program for $d$-codes in $B_n$ and $C_n$. Second for even~$d$, 
	it was shown in \cite[Prop.~3.3]{SchmidtWeissSteiner} that every $d$-code in $D_n$ induces a $\frac d2$-code in $\frac12D_n$ of the same size.
	We can therefore require without loss of generality that the inner distribution $(A_0,A_1,\dots,A_n)$ of a $d$-code in $D_n$ satisfies $A_i=0$ for all odd $i$, and thus, by adding this constraint to the linear program~\eqref{eq:LPprimalCodes} for $d$-codes in $D_n$, we can study the linear program for $\frac d2$-codes in $\frac12 D_n$ instead. These observations imply the following corollary from Theorem~\ref{thm:LPopt_matrix_qJohn_herm_hyp}.
	\begin{cor}\label{cor:LPopt_Bn_Cn_Dn}
		~\begin{enumerate}[(a)]
			\item Let $X$ be the set of generators in $B_n$ or $C_n$. Assume that $d$ is odd with $1\leq d\leq n$. Then the LP optimum for $d$-codes in $C_n$ and $B_n$ is given by
			\begin{align*}
				\LP(d)=\begin{cases}
					|X|\prod\limits_{i=1}^{\frac{d-1}{2}}\dfrac{q^{2i-1}-1}{q^{n+2i-1}-1}&\text{for odd $n$}\\[10pt]
					|X|\prod\limits_{i=1}^{\frac{d-1}{2}}\dfrac{q^{2i-1}-1}{q^{n+2i}-1}&\text{for even $n$}.
				\end{cases}
			\end{align*}
			\item Let $X$ be the set of generators in $D_n$. Assume that $d$ is even with $2\leq d\leq n$. Then the LP optimum for $d$-codes in $D_n$ is given by
			\[
			\LP(d)=\begin{cases}
				\dfrac{|X|}{2}\prod\limits_{i=1}^{\frac d2 -1}\dfrac{q^{2i-1}-1}{q^{n+2i-1}-1}&\text{for odd $n$}\\[10pt]
				\dfrac{|X|}{2}\prod\limits_{i=1}^{\frac d2 -1}\dfrac{q^{2i-1}-1}{q^{n+2i-2}-1}&\text{for even $n$}.
			\end{cases}
			\]
		\end{enumerate}
	\end{cor}
	
	We note that the LP optima in Theorem~\ref{thm:LPopt_matrix_qJohn_herm_hyp} \ref{LPoptordinaryq} and Corollary~\ref{cor:LPopt_Bn_Cn_Dn} for the polar spaces are exactly the bounds obtained in~\cite{SchmidtWeissSteiner}. We conjecture that the bound from Corollary~3.4.~(d) in~\cite{SchmidtWeissSteiner} for~$D_n$ is precisely the LP optimum if $n$ and $d$ are odd. This was checked with a computer for many small values of $q$, $n$, and~$d$. We also note that for $\qjohnson$, the LP optimum is the well-known Singleton bound (also called \emph{Wang-Xing-Safavi-Naini}, \emph{anticode bound}, or \emph{packing bound}):
	\[
	\LP(d)=\frac{\qbin{m+n}{n-d+1}}{\qbin{n}{n-d+1}},
	\]
	see~\cite[Thm.~5.2]{WangXingSafavi-Naini}, \cite{ZhangLPqJohn}, and \cite[Thm.~1]{EtzionVardy}. Moreover, it was known that \eqref{eq:LPbound_matrixschemes} is an upper bound on the size of $d$-codes in $\bil$ and $\alt$, and for odd $d$ in $\her$, which was proved by using Delsarte's linear program in \cite{DelsarteBilinear}, \cite{DelGoe1975}, and \cite{SchmidtHermitian}. Because of this and since sharp constructions were given in the aforementioned papers except for $\alt$ if $m$ is even and $q$ is odd, the LP optimum \eqref{eq:LPbound_matrixschemes} was basically known, except in the latter case (see Section~\ref{sec:linearprogramming} for a detailed explanation of this argument). However, we derive the LP optimum in $\bil$, $\alt$, and $\her$ for all parameters without using the known constructions.
	
	We now show how upper bounds on intersecting sets in classical association schemes follow immediately from the LP optimum for codes. These bounds are commonly known as Erd\H{o}s--Ko--Rado-type bounds. To establish them, we require the following lemma.
	
	\begin{lem}[{\cite[Thm.~3]{Tarnanen1999}}]\label{lem:LP_D_Dbar}
		Let $(X,(R_i))$ be an association scheme with $n$ classes. Let $D$ be a subset of $\{1,\dots,n\}$ and let $\overline{D}=\{1,\dots,n\}\setminus D$. Then the LP optima $\LP(D)$ and $\LP(\overline{D})$ of the LP~\eqref{eq:LPprimalCodes} for $D$-codes and $\overline{D}$-codes, respectively, satisfy $\LP(D)\LP(\overline{D})\leq |X|$.
	\end{lem}

	We can now directly derive the next corollary from Theorem~\ref{thm:LPopt_matrix_qJohn_herm_hyp} and Corollary~\ref{cor:LPopt_Bn_Cn_Dn}. 
	\begin{cor}\label{cor:EKR}
		Let $n$ and $t$ be integers with $1\leq t\leq n$.
		\begin{enumerate}[(a)]
			\item\label{EKR_John} Let $Y$ be a set of $n$-subsets of a $v$-set such that $|x\cap y|\geq t$ for all $x,y\in Y$. If $v$ is sufficiently large and if $\binom{n-i}{n-t}$ divides $\binom{v-i}{v-t}$ for all $i=0,1,\dots,t-1$, then we have $|Y|\leq \binom{v-t}{n-t}$.
			\item\label{EKR_Hamming} Let $Y$ be a set of $n$-tuples over a set with $q$ elements such that all pairs in $Y\times Y$ coincide in at least $t$ positions. If $q\geq\max\{t+1,n-t+1\}$, then we have $|Y|\leq q^{n-t}$.
			\item\label{EKR_qJohn} Let $Y$ be a set of $n$-dimensional subspaces of $\FF_q^v$ with $v\geq 2n$ such that $\dim(x\cap y)\geq t$ for all $x,y\in Y$. Then we have $|Y|\leq \qbin{v-t}{n-t}$.
			\item\label{EKR_bil} Let $Y$ be a set of $m\times n$ matrices over $\FF_q$ with $m\geq n$ such that $\rk(x-y)\leq n-t$ for all $x,y\in Y$. Then we have $|Y|\leq q^{m(n-t)}$.
			\item\label{EKR_alt} Let $Y$ be a set of alternating $n\times n$ matrices over $\FF_q$ such that $\rk(x-y)\leq n-t$ for all $x,y\in Y$. Then we have
			\[
			|Y|\leq \begin{cases}
				q^{(n-t)(n-1)/2}&\text{if $n$ and $t$ are even},\\
				q^{n(n-t)/2}&\text{if $n$ and $t$ are odd},\\
				q^{(n-t-1)(n-1)/2}&\text{if $n$ is even and $t$ is odd},\\
				q^{n(n-t-1)/2}&\text{if $n$ is odd and $t$ is even}.\\
			\end{cases}
			\]
			\item Let $Y$ be a set of Hermitian $n\times n$ matrices over $\FF_{q^2}$ such that $\rk(x-y)\leq n-t$ for all $x,y\in Y$. Then we have
			\[
			|Y|\leq \begin{cases}
				q^{n(n-t)}&\text{if $n-t$ is even,}\\
				q^{n(n-t)}\,\dfrac{q^{t+1}+q^t}{q^{n+t}+q^n-q^{t+1}+1}&\text{if $n$ is even and $t$ is odd,}\\[10pt]
				q^{n(n-t)}\,\dfrac{q^{t+1}+q^t}{q^{n+t}-q^n+q^{t+1}+1}&\text{if $n$ is odd and $t$ is even.}\\
			\end{cases}
			\]
			\item\label{EKR_polar} Let $Y$ be a set of generators in a polar space $\cP$ of rank $n$ such that $\dim(x\cap y)\geq t$ for all $x,y\in Y$. If $\cP=\hermitianpol{2n-1}$, then we have
			\begin{align*}
				|Y|\leq \begin{cases}
					\prod\limits_{i=0}^{n-t-1}\dfrac{q^{n+1+i}+(-1)^{n+i}}{q^{i+1}+(-1)^{i+1}}&\text{if $n-t$ is even,} \\[10pt]
					\dfrac{(-1)^{n+1}}{\epsilon(n,n-t+1)}\prod\limits_{i=0}^{n-t-1}\dfrac{q^{n+1+i}+(-1)^{n+i}}{q^{i+1}+(-1)^{i+1}}&\text{if $n-t$ is odd.}
				\end{cases}
			\end{align*}
			If $\cP=B_n$ or $C_n$, then we have
			\begin{align*}
				|Y|\leq \begin{cases}
					\prod\limits_{i=1}^{\frac{n-t}{2}} \dfrac{q^{n+2i-1}-1}{q^{2i-1}-1}&\text{if $n$ and $t$ are odd,}\\[10pt]
					\prod\limits_{i=1}^{\frac{n-t}{2}} \dfrac{q^{n+2i}-1}{q^{2i-1}-1}&\text{if $n$ and $t$ are even.}
				\end{cases}
			\end{align*}
			If $\cP=D_n$, then we have
			\begin{align*}
				|Y|\leq \begin{cases}
					2\prod\limits_{i=1}^{\frac{n-t-1}{2}} \dfrac{q^{n+2i-1}-1}{q^{2i-1}-1}&\text{if $n$ is odd and $t$ is even,} \\[10pt]
					2\prod\limits_{i=1}^{\frac{n-t-1}{2}} \dfrac{q^{n+2i-2}-1}{q^{2i-1}-1}&\text{if $n$ is even and $t$ is odd.}
				\end{cases}
			\end{align*}
			\item\label{EKR_half} Let $Y$ be a set of generators in one of the bipartite halves $\frac12 D_n$ of the hyperbolic polar space of rank $n$ such that $\dim(x\cap y)\geq t$ for all $x,y\in Y$. Then we have
			\[
			|Y|\leq\begin{cases}
				\prod\limits_{i=0}^{\lfloor \frac{n-t-2}{2}\rfloor} \dfrac{q^{2n+2i}-1}{q^{2i+1}-1}&\text{if $n$ is even,}\\[10pt]
				\prod\limits_{i=0}^{\lfloor \frac{n-t-2}{2}\rfloor} \dfrac{q^{2n+2i+2}-1}{q^{2i+1}-1}&\text{if $n$ is odd.}
			\end{cases}
			\]
		\end{enumerate}
	\end{cor}
	\begin{proof}
		In all the cases except \ref{EKR_alt} and \ref{EKR_half}, the set $Y$ is a $D$-code with $D=\{1,\dots,n-t\}$. Applying Lemma~\ref{lem:LP_D_Dbar} gives
		\[
		|Y|\leq \frac{|X|}{\LP(n-t+1)},
		\]
		where $\LP(n-t+1)$ is the corresponding LP optimum for $(n-t+1)$-codes from Theorem~\ref{thm:LPopt_matrix_qJohn_herm_hyp} or Corollary~\ref{cor:LPopt_Bn_Cn_Dn}. This gives the required bounds. In cases \ref{EKR_alt} and \ref{EKR_half}, the bounds follow similarly by using that the set~$Y$ is a $D$-code in $\operatorname{Alt}_q(n)$ and in $\frac12D_n$ with $D=\{1,\dots,\lfloor\frac{n-t}{2}\rfloor\}$, respectively.
	\end{proof}

	We note that the results in Corollary~\ref{cor:EKR} \ref{EKR_John}, \ref{EKR_Hamming}, \ref{EKR_qJohn}, and \ref{EKR_bil} are known, see \cite{EKR}, \cite{Moon}, \cite{FranklWilson}, and \cite{Huang}, respectively. Moreover, for $t=1$ in Corollary~\ref{cor:EKR} \ref{EKR_polar}, our bounds coincide with the ones given in \cite{StantonEKR}, which can be shown by a short proof of induction. However, our proof shows a new approach to obtain Erd\H{o}s--Ko--Rado-type bounds for subsets in association schemes.
	
	The following proposition provides a more useful version of the bounds on $t$-intersecting sets in polar spaces derived from Corollary~\ref{cor:EKR}. This version allows for a direct comparison with the so-far best known bounds for $t>1$ obtained in \cite{IhringerMetsch}. In that work, the maximal $t$-intersecting sets in all polar spaces were classified for $t\geq n-\sqrt{8n/5}+2$ if $q\geq 3$ and for $t\geq n-\sqrt{8n/9}+2$ if $q=2$. For smaller $t$, the exact size of the maximal $t$-intersecting sets remains unknown, but \cite[Thm.~49]{IhringerMetsch} also gives general upper bounds for all cases. Comparing the bounds from Proposition~\ref{prop:bounds_on_EKRbounds} with those in \cite[Thm.~49]{IhringerMetsch} shows that the former improve upon the latter in all cases except for $t=2$ in $B_n$ and $C_n$ with even $n\geq 4$, and for $t=2$ with $q=3,4,5,7$ in $D_n$ with odd $n\geq 5$. However, the bounds in Corollary~\ref{cor:EKR} are still far away from the size of the largest examples of $t$-intersecting sets currently known.
	
	\begin{prop}\label{prop:bounds_on_EKRbounds}
		Let $Y$ be a set of generators in a polar space $\cP$ of rank $n$ such that $\dim(x\cap y)\geq t$ for all $x,y\in Y$. Suppose that $1<t<n$ and $q\geq 2$.
		\begin{enumerate}[(a)]
			\item If $\cP=\hermitianpol{2n-1}$, then we have
			\begin{align*}
				|Y|\leq\begin{cases}
					8q^{n(n-t)}&\text{if $n-t$ is even},\\
					43q^{n(n-t-1)+1}&\text{if $n$ is odd and $t$ is even},\\
					26q^{n(n-t-1)+1}&\text{if $n$ is even and $t$ is odd}.
				\end{cases}
			\end{align*}
			\item If $\cP=B_n$ or $C_n$, then we have
			\begin{align*}
				|Y|\leq\begin{cases}
					4 q^{n(n-t)/2}&\text{if $n$ and $t$ are odd},\\
					4 q^{(n+1)(n-t)/2}&\text{if $n$ and $t$ are even}.
				\end{cases}
			\end{align*} 
			\item If $\cP=D_n$, then we have
			\begin{align*}
				|Y|\leq\begin{cases}
					8q^{n(n-t-1)/2}&\text{if $n$ is odd and $t$ is even},\\
					8q^{(n-1)(n-t-1)/2}&\text{if $n$ is even and $t$ is odd}.
				\end{cases}
			\end{align*}
		\end{enumerate}
	\end{prop}
	\begin{proof}
		We require the bound 
		\begin{align}\label{eq:ineq_product_minus}
			\prod_{i=1}^n \left(1-\frac{1}{q^i}\right)\geq \frac 14
		\end{align}
		for all integers $n\geq 1$ and $q\geq 2$, which follows by using $1-x\geq 4^{-x}$ for $0\leq x\leq\frac 12$ and thus
		\[
		\prod_{i=1}^n \left(1-\frac{1}{q^i}\right)
		\geq \prod_{i=1}^n \left(1-\frac{1}{2^i}\right)
		\geq \prod_{i=1}^n 4^{-1/2^i}
		\geq 4^{-\sum_{i=1}^\infty 2^{-i}}
		=\frac 14.
		\]

		Assume that $\cP=B_n$ or $C_n$, and $n,t$ are odd. From Corollary~\ref{cor:EKR} \ref{EKR_polar} and \eqref{eq:ineq_product_minus}, it follows that
		\begin{align*}
			|Y|
			&\leq\prod\limits_{i=1}^{\frac{n-t}{2}} \frac{q^{n+2i-1}-1}{q^{2i-1}-1}
			\leq\prod\limits_{i=1}^{\frac{n-t}{2}} \frac{q^n}{1-\frac{1}{q^{2i-1}}}
			\leq 4q^{n(n-t)/2}.
		\end{align*}
		The other bounds in (b) and (c) can be proved similarly.

		Assume now that $\cP=\hermitianpol{2n-1}$ and $n-t$ is even. From Corollary~\ref{cor:EKR} \ref{EKR_polar}, we find that
		\begin{align*}
			|Y|&\leq \prod\limits_{i=0}^{n-t-1}\dfrac{q^{n+1+i}+(-1)^{n+i}}{q^{i+1}+(-1)^{i+1}}
			\leq \prod\limits_{i=1}^{\frac{n-t}{2}}\frac{(q^{n+2i}-1)(q^{n+2i-1}+1)}{(q^{2i}+1)(q^{2i-1}-1)}
			\leq \prod\limits_{i=1}^{\frac{n-t}{2}} q^{2n} \left(\frac{1+\frac{1}{q^{n+2i-1}}}{1-\frac{1}{q^{2i-1}}}\right).
		\end{align*}
		By using \eqref{eq:ineq_product_minus} and 
		\[
		\prod_{i=1}^{\frac{n-t}{2}} \left(1+\frac{1}{q^{n+2i-1}}\right)<2,
		\]
		see \cite[Lem.~3.6]{SchmidtWeissSteiner}, we obtain the required bound $|Y|\leq 8q^{n(n-t)}$. 

		Lastly, assume that $\cP=\hermitianpol{2n-1}$ and $n-t$ is odd. Using similar steps as in the case of even $n-t$, we obtain
		\begin{align}\label{eq:upperbound_tintersect_hermitian}
			|Y|&\leq 8q^{n(n-t-1)} \frac{(-1)^{(n+1)}}{\epsilon(n,n-t+1)} \frac{q^{2n-t}+(-1)^n}{q^{n-t}-1}.
		\end{align}
		We will later give bounds on $\epsilon(n,d)$ in Lemma~\ref{lem:eps_hermpol}. Namely, we will see that
		\begin{align}\label{eq:upperbound_epsilon}
			\frac{(-1)^{n+1}}{\epsilon(n,n-t+1)}<\begin{cases}
				2q^{-n+t+1}\dfrac{1}{q^t-1}&\text{for odd $n$}\\[10pt]
				\dfrac{q+1}{q^n}&\text{for even $n$}.
			\end{cases}
		\end{align}
		Combining \eqref{eq:upperbound_tintersect_hermitian} and \eqref{eq:upperbound_epsilon} gives
		\begin{align}\label{eq:upperbound_tintersect_hermitian2}
			|Y|\leq \begin{cases}
				16 q^{n(n-t-1)-n+t+1}\dfrac{q^{2n-t}-1}{(q^t-1)(q^{n-t}-1)}&\text{for odd $n$}\\[15pt]
				8 q^{n(n-t-1)} \dfrac{(q+1)(q^{2n-t}+1)}{q^n(q^{n-t}-1)}&\text{for even $n$}.
			\end{cases}
		\end{align}
		Since $q\geq 2$, $t\geq 2$, and $n-t\geq 1$, we obtain
		\begin{align}\label{eq:frac1_bound}
			\frac{q^{2n-t}-1}{(q^t-1)(q^{n-t}-1)}&\leq q^{n-t}\frac{1}{\left(1-\frac{1}{q^t}\right)\left(1-\frac{1}{q^{n-t}}\right)}\leq \frac83 q^{n-t}.
		\end{align}
		Moreover, since also $n\geq 3$, we have
		\begin{align}\label{eq:frac2_bound}
			\frac{(q+1)(q^{2n-t}+1)}{q^n(q^{n-t}-1)}
			=q\;\frac{\left(1+\frac1q\right)\left(1+\frac{1}{q^{2n-t}}\right)}{1-\frac{1}{q^{n-t}}}
			\leq \frac{51}{16}q.
		\end{align}
		Combining \eqref{eq:upperbound_tintersect_hermitian2}, \eqref{eq:frac1_bound}, and \eqref{eq:frac2_bound} gives the required bounds.
	\end{proof}

	\section{\texorpdfstring{$P$- and $Q$-numbers of the ordinary and affine $q$-analogs}{P- and Q-numbers of the ordinary and affine q-analogs}}\label{sec:PQnum}

	Let $(X,(R_i))$ be an association scheme with $n$ classes with $P$-numbers $P_i(k)$ and $Q$-numbers $Q_k(i)$. Then for all $i,j,k=0,1,\dots,n$, these numbers satisfy
	\begin{align}
		\frac{1}{|X|}\sum_{k=0}^n P_i(k)Q_k(j)=\delta_{ij}\label{eq:PQnumorth}\\
		\mu_k P_i(k)=v_i Q_k(i)\label{eq:AskeyWilsonDuality}
	\end{align}
	and the two orthogonality relations
	\begin{align*}
		\frac{1}{|X|}\sum_{k=0}^n \mu_k P_i(k) P_j(k)=\delta_{ij}v_i
		\quad\text{and}\quad
		\frac{1}{|X|}\sum_{i=0}^n v_i Q_k(i) Q_j(i)=\delta_{jk}\mu_k.
	\end{align*}
	In the case of a $P$- or $Q$-polynomial association scheme with $P_i(k)=f_i(x_k)$ or $Q_k(i)=g_k(z_i)$, the associated polynomials $f_i$ or $g_k$ are orthogonal with respect to the inner product
	\begin{align}\label{eq:innerprod_Qpoly}
		(f,g)=\sum_{k=0}^n \mu_k f(y_k) g(y_k)
		\quad\text{or}\quad
		(f,g)=\sum_{i=0}^n v_i f(z_i) g(z_i) 
	\end{align}
	for all $f,g\in\RR[x]$,	respectively.
	
	To give the $P$- and $Q$-numbers of the ordinary and affine $q$-analogs, we need the \emph{$q$-Pochhammer symbol}~$\pochq{a}{n}$ defined by
	\[
	\pochq{a}{0}=1,\quad \pochq{a}{n}=\prod_{i=0}^{n-1}(1-aq^i)
	\]
	for a positive integer $n$ and a real number $a$. We note that the $q$-Pochhammer symbol is connected to the $q$-binomial coefficient as follows
	\begin{align}\label{eq:qbin_poch}
		\qbin{n}{k}=&\,\frac{\pochq{q^{-n}}{k}}{\pochq{q}{k}} (-1)^k q^{kn-\binom{k}{2}}.
	\end{align}
	Moreover, we need the \emph{$q$-hypergeometric function}~$\mbox{}_{r}\phi_{s}$ defined by
	\[
	\hyper{r}{s}{a_1,\dots,a_r}{b_1,\dots,b_s}{q,z}=\sum_{\ell=0}^{\infty} \frac{(a_1;q)_\ell\cdots(a_r;q)_\ell}{(b_1;q)_\ell\cdots(b_s;q)_\ell}\,(-1)^{(1+s-r)\ell} q^{(1+s-r)\binom{\ell}{2}}\,\frac{z^\ell}{(q;q)_\ell}.
	\]
	In the remainder of the paper, we write
	\[
		(a)_k=(a;b)_k,
	\]
	where $b$ is given as in \eqref{eq:bc}.

	We can now look at the $P$- and $Q$-numbers of the ordinary and affine $q$-analogs and we will see that in every case, the $P$- and $Q$-numbers are given by polynomials, establishing the $P$- and $Q$-polynomial property of the association schemes. The valencies and multiplicities can be found in Table~\ref{table:ClassicalASvalmul}.

	\begin{table}[ht]
		\caption{Valencies $v_i$ and multiplicities $\mu_k$ of the ordinary and affine $q$-analogs.}
		\centering
		\renewcommand*{\arraystretch}{1.15}
		\begin{tabular}{ccc}
			\toprule[0.4mm]
			\begin{tabular}{@{}c@{}}association \\ scheme\end{tabular} & valency $v_i$ & multiplicity $\mu_k$\\
			\midrule[0.4mm]\\[-2ex]
			\begin{tabular}{@{}c@{}} $\hermitianpol{2n-1},\, \hermitianpol{2n}$\\ $C_n$,\, $D_n$\\ $B_n$,\, $\elliptic{n+1}$ \end{tabular} & $\displaystyle p^{\binom{i+1}{2}+ie}\pbin{n}{i}$ & $\displaystyle p^{k(k-n)}\pbin{n}{k}\frac{\pochp{-p^{e+1}}{n}}{\pochp{-p^{e-k+1}}{n-k}\pochp{-p^{k-n-e-1}}{k}}$\\[5ex]
			$J_q(n,m)$ & $\displaystyle q^{i^2}\qbin{n}{i}\qbin{m}{i}$ & $\displaystyle \qbin{m+n}{k}-\qbin{m+n}{k-1}$ \\[3ex]
			\begin{tabular}{@{}c@{}}$\bil$\\ $\alt$\\ $\her$ \end{tabular} & $\displaystyle b^{\binom{i}{2}}\bbin{n}{i}\prod_{j=0}^{i-1}(cb^{n-j}-1)$ & $\displaystyle b^{\binom{k}{2}}\bbin{n}{k}\prod_{j=0}^{k-1}(cb^{n-j}-1)$\\[4ex]
			\bottomrule[0.4mm]
		\end{tabular}
		\label{table:ClassicalASvalmul}
	\end{table}
			
	\subsubsection*{Affine schemes}
	Consider $\bil$, $\alt$, and $\her$ with the parameters $b$ and $c$ given as in~\eqref{eq:bc}. Then there exists a unique ordering of the primitive idempotents $E_0,E_1,\dots,E_n$ such that the $P$- and $Q$-numbers are given by
	\begin{align}\label{eq:Pnum_affineschemes}
		P_i(k)=Q_i(k)=\sum_{j=0}^i (-1)^{i-j} b^{\binom{i-j}{2}} \bbin{n-j}{n-i}\bbin{n-k}{j} (cb^n)^j,
	\end{align}
	see \cite{DelsarteBilinear}, \cite{DelGoe1975}, and \cite{SchmidtHermitian}. Moreover, $P_i(k)$ is a polynomial of degree~$i$ in~$b^{-k}$, called \emph{affine $q$-Krawtchouk polynomial}, namely we have
	\begin{align*}
		P_i(k)=Q_i(k)=v_i\hyp{b^{-k},b^{-i},0}{c^{-1}b^{-n}, b^{-n}}{b;b},
	\end{align*}
	cf.\ \cite[\S~14.16]{Koekoek}. The affine schemes are thus $P$-polynomial with respect to the ordering $R_0,R_1,\dots,R_n$ and $Q$-polynomial with respect to the ordering $E_0,E_1,\dots,E_n$ of the primitive idempotents that is imposed by~\eqref{eq:Ri_bil}, \eqref{eq:Ri_alt}, or \eqref{eq:Ri_her}, and~\eqref{eq:Pnum_affineschemes}.
	
	\subsubsection*{$q$-Johnson scheme}
	Consider $\qjohnson$. Then there exists a unique ordering of the primitive idempotents $E_0,E_1,\dots,E_n$ such that the $P$-numbers are given by
	\begin{align}\label{eq:Pnum_qJohnson}
		P_i(k)=\sum_{j=0}^i (-1)^{i-j} \qbin{n-j}{i-j}\qbin{n-k}{j}\qbin{m+j-k}{j} q^{jk+\binom{i-j}{2}}
	\end{align}
	and the $Q$-numbers are given by
	\begin{align*}
		Q_k(i)=\mu_k \sum_{j=0}^k (-1)^j q^{\binom{j}{2}}\qbin{k}{j}\qbin{m+n+1-k}{j} \qbin{n}{j}^{-1}\qbin{m}{j}^{-1}\qbin{i}{j} q^{-ij},
	\end{align*}
	see \cite[Thm.~10]{DelsarteSemilattice}, \cite[\S~2]{DelsarteHahnpoly}. Moreover, $P_i(k)$ is a polynomial of degree $i$ in $[k]_q[m+n+1-k]_q$, called \emph{dual $q$-Hahn polynomial}, where $[n]_q$ is the \emph{$q$-number} defined by $[n]_q=(q^n-1)/(q-1)$,
	and $Q_k(i)$ is a polynomial of degree $k$ in $q^{-i}$, called \emph{$q$-Hahn polynomial}. Namely, we have
	\begin{align}\label{eq:Pnum_qJohnson_hypgeo}
		P_i(k)=v_i\hyp{q^{-i},q^{-k},q^{k-m-n-1}}{q^{-m},q^{-n}}{q;q}
	\end{align}
	and
	\begin{align}\label{eq:Qnum_qJohnson_hypgeo}
		Q_k(i)=\mu_k\hyp{q^{-i},q^{-k},q^{k-m-n-1}}{q^{-m},q^{-n}}{q;q},
	\end{align}
	cf.\ \cite[\S~14.5, 14.6]{Koekoek}. The $q$-Johnson scheme is thus $P$-polynomial with respect to the ordering $R_0,R_1,\dots,R_n$ and $Q$-polynomial with respect to the ordering $E_0,E_1,\dots,E_n$ of the primitive idempotents that is imposed by~\eqref{eq:Ri_qJohnson} and~\eqref{eq:Pnum_qJohnson}.
	
	\subsubsection*{Polar spaces scheme}
	Consider the polar spaces schemes. Then there exists a unique ordering of the primitive idempotents $E_0,E_1,\dots,E_n$ such that the $P$- and $Q$-numbers are given by
	\begin{align}\label{eq:Pnum_polarspaces}
		P_i(k)=v_i\pbin{n}{k}^{-1}\sum_{\ell=0}^i (-1)^\ell \pbin{n-i}{k-\ell}\pbin{i}{\ell} p^{\ell(\ell-i-e-1)}
	\end{align}
	and
	\begin{align*}
		Q_k(i)=\mu_k \pbin{n}{k}^{-1} \sum_{\ell=0}^i (-1)^\ell \pbin{n-i}{k-\ell} \pbin{i}{\ell} p^{\ell(\ell-i-e-1)},
	\end{align*}
	see \cite[Eq.~(8.1)]{StantonSome}, \cite[Prop.~2.4]{StantonThree}.\footnote{It should be noted that $p$ is assumed to be odd in~\cite{StantonSome} and~\cite{StantonThree}. However, all parameters of the association scheme as well as~$P_i(k)$ and $Q_k(i)$ are polynomials in $p$. Hence, the expressions for $P_i(k)$ and $Q_k(i)$ hold for all~$p$.} Moreover, $P_i(k)$ is a polynomial of degree $i$ in $p^{-k}$ and $Q_k(i)$ is a polynomial of degree $k$ in $p^{-i}$, both called \emph{$q$-Krawtchouk polynomial}, namely we have
	\begin{align*}
		P_i(k)=v_i\hyp{p^{-k},p^{-i},-p^{-n-e-1+k}}{0,p^{-n}}{p;p}
	\end{align*}
	and
	\begin{align*}
		Q_k(i)=\mu_k\hyp{p^{-k},p^{-i},-p^{-n-e-1+k}}{0,p^{-n}}{p;p},
	\end{align*}
	cf.\ \cite[\S~14.15]{Koekoek}. The polar space schemes are thus $P$-polynomial with respect to the ordering $R_0,R_1,\dots,R_n$ and $Q$-polynomial with respect to the ordering $E_0,E_1,\dots,E_n$ of the primitive idempotents that is imposed by~\eqref{eq:Ri_polarspaces} and~\eqref{eq:Pnum_polarspaces}. We call both orderings the \emph{standard orderings}.
	
	In \cite{SchmidtWeissSteiner}, it was shown that the $P$- and $Q$-numbers of $\hermitianpol{2n-1}$ and $\hyperbolichalf$ can be written in a unified way as $q$-Hahn polynomials by using the parameters $b$ and $c$ from \eqref{eq:bc}. First, the association scheme~$\hermitianpol{2n-1}$ is $Q$-polynomial with respect to two different orderings: the standard ordering $E_0,E_1,\dots,E_n$ and the second ordering $E_0,E_n,E_1,E_{n-1},E_2,E_{n-2},\dots$ \cite{ChiharaStanton}. We continue to use $P_i(k)$ and $Q_k(i)$ to denote the $P$- and $Q$-numbers with respect to the standard ordering and we use $P'_i(k)$ and $Q'_k(i)$ to denote the $P$- and $Q$-numbers with respect to the second ordering. Similarly, let $P'_i(k)$ and $Q'_k(i)$ denote the $P$- and $Q$-numbers, respectively, of $\hyperbolichalf$ and of $\qjohnson$. Then, due to~\cite[Eq.~(3.5), (3.6)]{SchmidtWeissSteiner}, \eqref{eq:Pnum_qJohnson_hypgeo}, and \eqref{eq:Qnum_qJohnson_hypgeo}, we have
	\begin{align}\label{eq:Pnum_qJohn_Herm_hyp}
		P'_i(k)=v'_i  \hyp{b^{-i},b^{-k},q^{-1}c^{-1}b^{-2n+k}}{b^{-n},c^{-1}b^{-n}}{b;b},
	\end{align}
	which is a polynomial of degree $i$ in $y_k=b^{-k}$,
	and
	\begin{align}\label{eq:Qnum_qJohn_Herm_hyp}
		Q'_k(i)=\mu'_k \hyp{b^{-i},b^{-k},q^{-1}c^{-1}b^{-2n+k}}{b^{-n},c^{-1}b^{-n}}{b;b},
	\end{align}
	which is a polynomial of degree $k$ in $z_i=b^{-i}$. The corresponding valencies $v'_i$ and multiplicities $\mu'_k$ are stated in Table~\ref{table:Pnum_parameters} and in the case of $\qjohnson$, they are equal to $v_i$ and $\mu_k$ as in Table~\ref{table:ClassicalASvalmul}.
	\begin{table}[ht]
		\caption{Valencies and multiplicities for $\hermitianpol{2n-1}$ and $\hyperbolichalf$ occurring in~\eqref{eq:Pnum_qJohn_Herm_hyp} and~\eqref{eq:Qnum_qJohn_Herm_hyp}.}
		\centering
		\renewcommand*{\arraystretch}{1.25}
		\begin{tabular}{lccc} 
			\toprule[0.4mm]
			& $v'_i$ & $\mu'_k$ \\
			\midrule[0.4mm]
			$\hermitianpol{2n-1}$ & $q^{i^2}\qqbin{n}{i}$ & \begin{tabular}{@{}c@{}}$\mu_{k/2}$ for even $k$ \\ $\mu_{n-(k-1)/2}$ for odd $k$\end{tabular} \\
			$\hyperbolichalf$ & $q^{\binom{2i}{2}}\qbin{m}{2i}$ & $\mu_k$ for $D_n$ \\
			\bottomrule[0.4mm]
		\end{tabular}
		\label{table:Pnum_parameters}
	\end{table}
	
	To apply Delsarte's linear programming method, we need some identities for the $P$-numbers, which are summarized in the next proposition.
	\begin{prop}~
		\begin{enumerate}[(a)]
			\item The $P$-numbers~\eqref{eq:Pnum_qJohn_Herm_hyp} of $\qjohnson$, $\hermitianpol{2n-1}$ (with respect to the second ordering), and $\hyperbolichalf$ satisfy
				\begin{align}\label{eq:Pnumidentity_qJohn_Herm_hyp}
				\sum_{i=0}^n \bbin{n-i}{j}P'_i(k)
				=b^{k(n-j)}\bbin{n-k}{n-j}\frac{\pochb{qcb^{n-k}}{n-j}}{\pochb{q}{n-j}}
			\end{align}
			for all $j,k=0,1,\dots,n$, where $n=\lfloor m/2\rfloor$ in the case of $\hyperbolichalf$.
			\item The $P$-numbers~\eqref{eq:Pnum_affineschemes} of $\bil$, $\her$, and $\alt$ satisfy
			\begin{align}\label{eq:Pnumidentity_affineschemes}
				\sum_{i=0}^n\bbin{n-i}{j} P_i(k)=\bbin{n-k}{n-j} (cb^n)^{n-j}
			\end{align}
			for all $j,k=0,1,\dots,n$, where $n=\lfloor m/2\rfloor$ in the case of $\alt$. 
		\end{enumerate}
	\end{prop}
	\begin{proof}
		The identity in (a) was proved for $\hermitianpol{2n-1}$ and $\hyperbolichalf$ in \cite{SchmidtWeissSteiner}. The same proof holds for $\qjohnson$. Alternatively, the identity for $\qjohnson$ can also be found in~\cite[Thm.~9]{DelsarteSemilattice}. The identity in (b) follows by using the $P$-numbers~\eqref{eq:Pnum_affineschemes} and the $q$-binomial inversion formula of the form
		\begin{align*}
			\sum_{j=i}^k  (-1)^{k-j}b^{\binom{k-j}{2}}\bbin{j}{i}\bbin{k}{j}=\delta_{ik}.
		\end{align*}
	\end{proof}
	We note that~\eqref{eq:Pnumidentity_qJohn_Herm_hyp} immediately implies an identity for the $Q$-numbers of $\qjohnson$, $\hermitianpol{2n-1}$ (with respect to the second ordering), and $\hyperbolichalf$ by using~\eqref{eq:PQnumorth}. Namely, we have
	\begin{align}\label{eq:Qnumidentity_qJohn_Herm_hyp}
		\sum_{k=0}^n b^{k(n-j)}\bbin{n-k}{n-j}\frac{\pochb{qcb^{n-k}}{n-j}}{\pochb{q}{n-j}} Q'_k(i)=|X| \bbin{n-i}{j}
	\end{align}
	for all $i,j=0,1,\dots,n$, see \cite[Lem.~3.2]{SchmidtWeissSteiner}.
	
	\section{Linear programming}\label{sec:linearprogramming}
	
	Here, we briefly summarize some facts from linear programming and explain the ideas behind the proof of our main result, Theorem~\ref{thm:LPopt_matrix_qJohn_herm_hyp}. For further information on linear programming, we refer to \cite{SchrijverLPbook} and \cite{Vanderbei}.

	Let $(X,(R_i))$ be an association scheme with $n$ classes and denote its $P$- and $Q$-numbers by $P_i(k)$ and $Q_k(i)$, respectively. It is well known from the theory of linear programming that the (primal) linear program~\eqref{eq:LPprimalCodes} for $D$-codes in $X$ has a \emph{dual linear program (dual LP)} of the form
	\begin{align}\label{eq:LPdualCodes}
		\begin{array}{lrcl}
			\operatorname*{minimize}\limits_{y_i\in\RR} & \sum\limits_{k=0}^n \mu_k y_k & &\\
			\text{subject to} & y_0&=& 1\\
			& y_k&\geq & 0\quad\text{for all $k=1,2,\dots,n$}\\
			& \sum\limits_{k=0}^n Q_k(i)y_k &\leq & 0\quad\text{for all $i\in D$},
		\end{array}
	\end{align}
	where $\mu_k$ denote the multiplicities of the association scheme. The terms feasible solution and optimal are analogously defined for the dual LP as for the primal LP. By using the well-known weak duality theorem \cite[Cor.~7.1g]{SchrijverLPbook} and complementary slackness conditions \cite[\S~7.9]{SchrijverLPbook} from optimization theory, we have the following powerful connection between the primal and dual LPs.
	\begin{cor}[{\cite[Lem.~3.5, \S~3.3.2]{DelsartePhD}}]\label{cor:dualLPcodes}
		Let $(X,(R_i))$ be an association scheme with $n$~classes and $Q$-numbers $Q_k(i)$, and let $D$ be a subset of $\{1,2,\dots,n\}$. If $\LP(D)$ denotes the LP optimum of~\eqref{eq:LPprimalCodes}, then the objective function value of every feasible solution $y_k$ of~\eqref{eq:LPdualCodes} gives an upper bound on $\LP(D)$, that is,
		\begin{align}\label{eq:dualLPbound}
			\LP(D)\leq \sum\limits_{k=0}^n \mu_k y_k.
		\end{align}
		Moreover, if $Y$ is a $D$-code with dual distribution $(A'_k)$ such that $|Y|$ equals the right-hand side of~\eqref{eq:dualLPbound}, then $y_kA'_k=0$ for all $k=1,2,\dots,n$.
	\end{cor}
	
	The proof of Theorem~\ref{thm:LPopt_matrix_qJohn_herm_hyp} relies on the following \emph{strong duality theorem}.
	
	\begin{thm}[{\cite[\S~7.9]{SchrijverLPbook}}]\label{thm:strongduality}
		Let $x$ and $y$ be feasible solutions of the primal LP~\eqref{eq:LPprimalCodes} and dual LP~\eqref{eq:LPdualCodes}, respectively. Then both solutions $x$ and $y$ are optimal if and only if their objective function values coincide, that is, $\sum\limits_{i=0}^n x_i=\sum\limits_{k=0}^n \mu_k y_k$.
	\end{thm}
	
	To prove Theorem~\ref{thm:LPopt_matrix_qJohn_herm_hyp}, we will derive a feasible solution of the primal LP~\eqref{eq:LPprimalCodes} and of the dual LP~\eqref{eq:LPdualCodes} for each considered association scheme such that the objective function values of both feasible solutions coincide. By the strong duality theorem, this will imply Theorem~\ref{thm:LPopt_matrix_qJohn_herm_hyp}.
	
	To construct a feasible solution of the dual LP~\eqref{eq:LPdualCodes}, we will apply the following theorem. It restates Corollary~\ref{cor:dualLPcodes} for $Q$-polynomial association schemes by using the orthogonal polynomials associated with the $Q$-polynomial structure.
	
	\begin{thm}[{\cite[\S~4.3]{DelsartePhD}}]\label{thm:polydualLP}
		Let $(X,(R_i))$ be a $Q$-polynomial association scheme with $n$~classes, where $Q_k(i)=g_k(z_i)$ for some $g_k\in\RR[x]$ of degree $k$ and some real numbers~$z_i$. Let $D$ be a subset of $\{1,2,\dots,n\}$. Suppose that $F\in\RR[x]$ is a polynomial of degree at most $n$ whose coefficients $F_k$ from the expansion $F=F_0g_0+F_1g_1+\cdots+F_ng_n$ satisfy $F_0=1$, $F_k\geq 0$ for all $k=1,2,\dots,n$, and $F(z_i)\leq 0$ for all $i\in D$. Then $(F_0,F_1,\dots,F_n)$ is a feasible solution of the dual LP~\eqref{eq:LPdualCodes} with objective function value $F(z_0)$. In particular, we have $\LP(D)\leq F(z_0)$. 
	\end{thm}
	
	The next remark shows how to compute the coefficients $F_k$ for a chosen polynomial~$F$.
	\begin{rem}
		Assume that $(X,(R_i))$ is $Q$-polynomial with $Q$-numbers given by $Q_k(i)=g_k(z_i)$ as in Theorem~\ref{thm:polydualLP}. Let $P_i(k)$, $v_i$, and $\mu_k$ be the $P$-numbers, the valencies, and the multiplicities of $(X,(R_i))$, respectively. Then, for a polynomial~$F$ with $F=F_0g_0+F_1g_1+\cdots+F_ng_n$, the coefficients $F_k$ can be computed by using the inner product~\eqref{eq:innerprod_Qpoly} and~\eqref{eq:AskeyWilsonDuality}, namely
		\begin{align}\label{eq:Fk}
			F_k=\frac{(F,g_k)}{(g_k,g_k)}=\frac{1}{|X|\mu_k}\sum_{i=0}^nv_i F(z_i)Q_k(i)=\frac{1}{|X|}\sum_{i=0}^n F(z_i)P_i(k).
		\end{align}
	\end{rem}
	
	\section{Feasible solution of the dual LP}\label{sec:SolDualLP}
	In this section, we will apply Theorem~\ref{thm:polydualLP} for $D=\{d,d+1,\dots,n\}$ together with a Singleton polynomial of the form
	\begin{align*}
	F(z)=c\prod_{i=d}^n (z-z_i)\quad\text{for some real constant $c$},
	\end{align*}
	or variations of it, to construct a feasible solution of the dual~LP~\eqref{eq:LPdualCodes} for the ordinary $q$-analogs $\qjohnson$, $\hermitianpol{2n-1}$, and $\hyperbolichalf$ as well as for the affine $q$-analogs $\bil$, $\her$, and $\alt$. We start with the latter three. By using the definition of $b$ and $c$, the number of elements in~$X$ for the affine $q$-analogs can be written as
	\begin{align}\label{eq:size_X_aff}
			|X|=(cb^n)^n.
	\end{align}
	\begin{prop}\label{prop:fsol_dual_matrixschemes}
		There exists a feasible solution of the dual LP~\eqref{eq:LPdualCodes} for $d$-codes in $\bil$, $\her$, and $\alt$ with objective function value~\eqref{eq:LPbound_matrixschemes} for $1\leq d\leq n$, where $d$ is required to be odd in the case of $\her$ and $n=\lfloor m/2\rfloor$ in the case of~$\alt$.
	\end{prop}
	\begin{proof}
		Recall from~\eqref{eq:Pnum_affineschemes} that $\bil$, $\her$, and $\alt$ are $Q$-polynomial with $z_i=b^{-i}$. Let $P_i(k)$ be given in \eqref{eq:Pnum_affineschemes}. Take the Singleton polynomial
		\[
		F(z)=(cb^n)^{n-d+1}\bbin{n}{d-1}\prod_{j=d}^n b^j\,\frac{z-z_j}{b^j-1}
		\]
		with $z_j=b^{-j}$. Then we have $F(z_i)=0$ for all $i=d,d+1,\dots,n$. Moreover, we obtain
		\[
		F(z_i)
		=(cb^n)^{n-d+1}\bbin{n}{d-1}\,\prod_{j=d}^n b^j\,\frac{b^{-i}-b^{-j}}{b^j-1}
		=(cb^n)^{n-d+1}\bbin{n-i}{n-d+1}
		\]
		for all $i=0,1,\dots,n$. Together with~\eqref{eq:Fk}, \eqref{eq:size_X_aff}, and \eqref{eq:Pnumidentity_affineschemes}, we have 
		\begin{align*}
			F_k
			=\frac{(cb^n)^{n-d+1}}{(cb^n)^n}\sum_{i=0}^n \bbin{n-i}{n-d+1} P_i(k)
			=\bbin{n-k}{d-1}.
		\end{align*}
		Observe that in the case of $\her$, the sign of $\bbin{n-k}{d-1}$ is $(-1)^{(d-1)(n-k+d-1)}=1$ since $d$ is odd. We thus obtain $F_k\geq 0$ for all $k=0,1,\dots,n$. Therefore, the polynomial $F/F_0$ satisfies all conditions of Theorem~\ref{thm:polydualLP} with $F(z_0)/F_0=(cb^n)^{n-d+1}$, which implies the stated result.
	\end{proof}
	
	We now look at $\her$ with even $d$ and proceed similarly as in Proposition~\ref{prop:fsol_dual_matrixschemes}, but we take a linear combination of two Singleton polynomials instead of just one.
	
	\begin{prop}\label{prop:fsol_dual_hermat_deven} 
		Let $n$ and $d$ be integers with $2\leq d\leq n$ and even $d$. Then there exists a feasible solution of the dual LP~\eqref{eq:LPdualCodes} for $d$-codes in $\her$ with objective function value~\eqref{eq:LPbound_Hermat_deven}.
	\end{prop}
	\begin{proof}
		Take
		\[
		F(z)=\bbin{n-1}{d-2}\beta_1(z)-\bbin{n-1}{d-1}\beta_2(z)
		\]
		with
		\begin{align*}
			\beta_1(z)&=q^{n(n-d+1)}\bbin{n}{d-1}\,\prod_{j=d}^n b^j\,\frac{z-z_j}{b^j-1}\\
			\beta_2(z)&=(-1)^{n+1}q^{n(n-d+2)}\bbin{n}{d-2}\,\prod_{j=d-1}^n b^j\,\frac{z-z_j}{b^j-1},
		\end{align*}
		where $z_j=b^{-j}$.	We then have $F(z_i)=0$ for all $i=d,d+1,\dots,n$. We also obtain
		\begin{align*}
			\beta_1(z_i)=q^{n(n-d+1)}\bbin{n-i}{n-d+1},
			\qquad
			\beta_2(z_i)=(-1)^{n+1}q^{n(n-d+2)}\bbin{n-i}{n-d+2}
		\end{align*}
		for all $i=0,1,\dots,n$. This gives
		\begin{align*}
			F(z_i)=q^{n(n-d+1)}\bbin{n-1}{d-2}&\bbin{n-i}{n-d+1}
			+(-1)^n q^{n(n-d+2)}\bbin{n-1}{d-1}\bbin{n-i}{n-d+2}
		\end{align*}
		for all $i=0,1,\dots,n$. From~\eqref{eq:Fk}, \eqref{eq:size_X_aff}, and \eqref{eq:Pnumidentity_affineschemes}, we find
		\begin{align*}
			F_k=(-1)^{n+1}\left(\bbin{n-1}{d-2}\bbin{n-k}{d-1}-\bbin{n-1}{d-1}\bbin{n-k}{d-2}\right).
		\end{align*}
		The sign of $\bbin{n-1}{d-2}\bbin{n-k}{d-1}$ and $\bbin{n-1}{d-1}\bbin{n-k}{d-2}$ is $(-1)^{n-k+1}$ and $(-1)^n$, respectively, which gives
		\[
		F_k=(-1)^k\Bigg\lvert \bbin{n-1}{d-2}\bbin{n-k}{d-1}\Bigg\rvert +\Bigg\lvert \bbin{n-1}{d-1}\bbin{n-k}{d-2} \Bigg\rvert.
		\]
		For even $k$, we immediately have $F_k\geq 0$. For $k=1$, we obtain $F_1=0$ and for all odd $k\geq 3$, we have
		\[
		\Bigg\lvert\frac{ \bbin{n-1}{d-1}\bbin{n-k}{d-2}}{ \bbin{n-1}{d-2}\bbin{n-k}{d-1}}\Bigg\rvert
		=\Bigg\lvert\frac{(-q)^{n-d+1}-1}{(-q)^{n-k-d+2}-1}\Bigg\rvert
		\geq \frac{q^{n-d+1}-1}{q^{n-k-d+2}+1}
		\geq 1.
		\]
		Hence, we obtain $F_k\geq 0$ for all $k=0,1,\dots,n$. Therefore, the polynomial $F/F_0$ satisfies all conditions of Theorem~\ref{thm:polydualLP} with
		\[
		\frac{F(z_0)}{F_0}=\frac{q^{n(n-d+1)}\left(\bbin{n-1}{d-2}\bbin{n}{d-1}+b^n\bbin{n-1}{d-1}\bbin{n}{d-2}\right)}{(-1)^{n+1}\left(\bbin{n-1}{d-2}\bbin{n}{d-1}-\bbin{n-1}{d-1}\bbin{n}{d-2}\right)}.
		\]
		After some elementary manipulations, this gives the stated objective function value and thus proves the proposition.
	\end{proof}
	
	We will now look at the ordinary $q$-analogs $\qjohnson$, $\hermitianpol{2n-1}$, and $\hyperbolichalf$ and again use Theorem~\ref{thm:polydualLP} with a Singleton polynomial. Similarly to $\her$, we will distinguish between even and odd $d$ in the case of $\hermitianpol{2n-1}$. Moreover, for the sake of convenience, we will rewrite the LP optima~\eqref{eq:LPbound_qJohn_herm_hyp} and~\eqref{eq:LPbound_hermpol_deven} as
	\[
	\LP(d)=\frac{|X|\pochb{q}{d-1}}{\pochb{qcb^n}{d-1}}
	\qquad\text{and}\qquad
	\LP(d)=\frac{|X|\pochb{q}{d-1}}{\pochb{qcb^n}{d-1}}\;\epsilon(n,d),
	\]
	respectively, in the remainder of this work.
	\begin{prop}\label{prop:fsol_dual_qJohn_herm_hyp}
		Let $X$ be the set of $n$-spaces in $\qjohnson$ or generators in $\hermitianpol{2n-1}$ or $\hyperbolichalf$, where $n=\lfloor m/2\rfloor$ in the case of $\hyperbolichalf$.  Then there exists a feasible solution of the dual LP~\eqref{eq:LPdualCodes} for $d$-codes in $\qjohnson$, $\hermitianpol{2n-1}$, and $\hyperbolichalf$ with objective function value~\eqref{eq:LPbound_qJohn_herm_hyp} for $1\leq d\leq n$, where $d$ is required to be odd in the case of $\hermitianpol{2n-1}$.
	\end{prop}
	\begin{proof}
		Recall from~\eqref{eq:Qnum_qJohn_Herm_hyp} that the association schemes $\qjohnson$, $\hermitianpol{2n-1}$, and $\hyperbolichalf$ are $Q$-polynomial with $z_i=b^{-i}$, where we take the second ordering for $\hermitianpol{2n-1}$. Let $P'_i(k)$ be given in~\eqref{eq:Pnum_qJohn_Herm_hyp}. Take the Singleton polynomial
		\[
		F(z)=\bbin{n}{d-1}\prod_{j=d}^n b^j\,\frac{z-z_j}{b^j-1}
		\]
		with $z_j=b^{-j}$. This gives $F(z_i)=0$ for all $i=d,d+1,\dots,n$ and
		\[
		F(z_i)=\bbin{n}{d-1}\prod_{j=d}^n b^j\,\frac{b^{-i}-b^{-j}}{b^j-1}=\bbin{n-i}{n-d+1}
		\]
		for all $i=0,1,\dots,n$. From~\eqref{eq:Fk} and~\eqref{eq:Pnumidentity_qJohn_Herm_hyp}, we find that
		\[
		F_k=\frac{1}{|X|}\sum_{i=0}^n \bbin{n-i}{n-d+1}P'_i(k)
		=\frac{1}{|X|} b^{k(d-1)}\bbin{n-k}{d-1}\frac{\pochb{qcb^{n-k}}{d-1}}{\pochb{q}{d-1}}.
		\]
		For $\qjohnson$ (since $m\geq k$) and $\hyperbolichalf$, we see that
		\[
		\frac{\pochb{qcb^{n-k}}{d-1}}{\pochb{q}{d-1}}\geq 0
		\]
		and thus, we have $F_k\geq 0$ for all $k=0,1,\dots,n$. For $\hermitianpol{2n-1}$, the sign of
		\[
		\frac{\pochb{qcb^{n-k}}{d-1}}{\pochb{q}{d-1}}=\frac{\poch{(-q)^{n-k+1}}{-q}{d-1}}{\poch{q}{-q}{d-1}}
		\]
		is $(-1)^{(n-k+1)(d-1)}$ and $\bbin{n-k}{d-1}$ has the sign $(-1)^{(d-1)(n-k-d+1)}$. Since $d$ is odd for $\hermitianpol{2n-1}$, we also obtain $F_k\geq 0$ for all $k=0,1,\dots,n$. Therefore, the polynomial~$F/F_0$ satisfies all conditions of Theorem~\ref{thm:polydualLP} with
		\[
		\frac{F(z_0)}{F_0}=\frac{|X|\pochb{q}{d-1}}{\pochb{qcb^n}{d-1}}.
		\]
		This concludes the proof.
	\end{proof}
	
	It remains to look at $\hermitianpol{2n-1}$ with even $d$, where we take a linear combination of two Singleton polynomials as in the proof of Proposition~\ref{prop:fsol_dual_hermat_deven} for $\her$ with even $d$.

	\begin{prop}\label{prop:fsol_dual_hermpol_deven}
		Let $X$ be the set of generators in $\hermitianpol{2n-1}$ and let $d$ be an even integer with $2\leq d\leq n$. Then there exists a feasible solution of the dual LP~\eqref{eq:LPdualCodes} for $d$-codes in $\hermitianpol{2n-1}$ with objective function value~\eqref{eq:LPbound_hermpol_deven}.
	\end{prop}
	\begin{proof}
		Take
		\[
		F(z)=\bbin{n-1}{d-2}\beta_1(z)-\bbin{n-1}{d-1}\beta_2(z)
		\]
		with
		\begin{align*}
			\beta_1(z)&=\bbin{n}{d-1}\prod_{j=d}^n b^j \,\frac{z-z_j}{b^j-1}\\
			\beta_2(z)&=b\,\frac{(b^{n+d-2}-1)}{(qb^{d-2}-1)}\bbin{n}{d-2}\prod_{j=d-1}^n b^j \,\frac{z-z_j}{b^j-1},
		\end{align*}
		where $z_j=b^{-j}$. We then obtain $F(z_i)=0$ for all $i=d,d+1,\dots,n$. Moreover, we have
		\begin{align*}
			\beta_1(z_i)=\bbin{n-i}{n-d+1},
			\qquad
			\beta_2(z_i)=b\,\frac{(b^{n+d-2}-1)}{(qb^{d-2}-1)}\bbin{n-i}{n-d+2}
		\end{align*}
		for all $i=0,1,\dots,n$. Therefore, we find
		\[
		F(z_i)=\bbin{n-1}{d-2}\bbin{n-i}{n-d+1}-b\,\frac{(b^{n+d-2}-1)}{(qb^{d-2}-1)}\bbin{n-1}{d-1}\bbin{n-i}{n-d+2}
		\]
		for all $i=0,1,\dots,n$. Using~\eqref{eq:Fk} and~\eqref{eq:Pnumidentity_qJohn_Herm_hyp} gives
		\begin{multline*}
			F_k=\frac{1}{|X|}\Bigg(b^{k(d-1)}\frac{\pochb{b^{n-k+1}}{d-1}}{\pochb{q}{d-1}}\bbin{n-1}{d-2}\bbin{n-k}{d-1}\\
			-b^{k(d-2)+1}\frac{(b^{n+d-2}-1)}{(qb^{d-2}-1)}\;\frac{\pochb{b^{n-k+1}}{d-2}}{\pochb{q}{d-2}}\bbin{n-1}{d-1}\bbin{n-k}{d-2}\Bigg).
		\end{multline*}
		Since the sign of $\bbin{m}{\ell}$ and $\pochb{b^m}{\ell}/\pochb{q}{\ell}$ is $(-1)^{\ell(m-\ell)}$ and $(-1)^{m\ell}$, respectively, we obtain
		\begin{multline*}
			F_k=\frac{1}{|X|}\Bigg((-1)^k\left\lvert b^{k(d-1)}\frac{\pochb{b^{n-k+1}}{d-1}}{\pochb{q}{d-1}}\bbin{n-1}{d-2}\bbin{n-k}{d-1}\right\rvert\\
			+\left\lvert b^{k(d-2)+1}\frac{(b^{n+d-2}-1)}{(qb^{d-2}-1)}\;\frac{\pochb{b^{n-k+1}}{d-2}}{\pochb{q}{d-2}}\bbin{n-1}{d-1}\bbin{n-k}{d-2}\right\rvert\Bigg).
		\end{multline*}
		We have $F_k\geq 0$ for all even $k$. Assume that $k$ is odd. This gives
		\begin{align*}
			\left\lvert\frac{ b^{k(d-2)+1}\frac{(b^{n+d-2}-1)}{(qb^{d-2}-1)}\,\frac{\pochb{b^{n-k+1}}{d-2}}{\pochb{q}{d-2}}\bbin{n-1}{d-1}\bbin{n-k}{d-2}}{b^{k(d-1)}\frac{\pochb{b^{n-k+1}}{d-1}}{\pochb{q}{d-1}}\bbin{n-1}{d-2}\bbin{n-k}{d-1}}\right\rvert
			=\left\lvert\frac{b^{-k+1}(b^{n+d-2}-1)(b^{n-d+1}-1)}{(b^{n-k+d-1}-1)(b^{n-k-d+2}-1)}\right\rvert.
		\end{align*}
		For $k=1$, this becomes $1$ and for $k\geq 3$, it can be bounded from below by
		\[
		q^{-k+1}\frac{(q^{n+d-2}+1)(q^{n-d+1}-1)}{(q^{n-k+d-1}-1)(q^{n-k-d+2}+1)}\geq 1.
		\]
		Hence, we also have $F_k\geq 0$ for all odd $k$. Thus, the polynomial~$F/F_0$ satisfies all conditions of Theorem~\ref{thm:polydualLP} with
		\[
		\frac{F(z_0)}{F_0}=\frac{\bbin{n-1}{d-2}\bbin{n}{d-1}-b\,\frac{(b^{n+d-2}-1)}{(qb^{d-2}-1)}\bbin{n-1}{d-1}\bbin{n}{d-2}}
		{\frac{1}{|X|}\Bigg(\frac{\pochb{b^{n+1}}{d-1}}{\pochb{q}{d-1}}\bbin{n-1}{d-2}\bbin{n}{d-1}\\
			-b\,\frac{(b^{n+d-2}-1)}{(qb^{d-2}-1)}\,\frac{\pochb{b^{n+1}}{d-2}}{\pochb{q}{d-2}}\bbin{n-1}{d-1}\bbin{n}{d-2}\Bigg)}.
		\]
		After some elementary manipulations, this gives the stated objective function value and thus proves the proposition.
	\end{proof}
	
	\section{Feasible solution of the primal LP}\label{sec:SolPrimalLP}
	The goal of this section is to derive a feasible solution of the primal LP~\eqref{eq:LPprimalCodes} whose objective function value equals the stated LP optimum in Theorem~\ref{thm:LPopt_matrix_qJohn_herm_hyp} for the affine $q$-analogs $\bil$, $\her$, $\alt$, and for the ordinary $q$-analogs $\qjohnson$, $\hermitianpol{2n-1}$, $\hyperbolichalf$. The strategy is to first compute the inner and dual distribution of a code whose size equals the respective stated LP optimum in Theorem~\ref{thm:LPopt_matrix_qJohn_herm_hyp}. Afterwards, we will show that these distributions are nonnegative and therefore, the inner distribution $(A_i)$ is a feasible solution of the primal LP~\eqref{eq:LPprimalCodes} such that its objective function value---the sum of the entries~$A_i$---is precisely the stated LP optimum in Theorem~\ref{thm:LPopt_matrix_qJohn_herm_hyp}. This is done in Section~\ref{subsec:feassol_primalLP_affine} and~\ref{subsec:feassol_primalLP_ordinaryq} for the affine and ordinary $q$-analogs, respectively, where $d$ is required to be odd for $\her$ and $\hermitianpol{2n-1}$. The case where $d$ is even in the latter two association schemes has to be treated separately in Section~\ref{subsec:feassol_primalLP_Hermat} and~\ref{subsec:feassol_primalLP_Hermpol}.
	
	In the following, we will frequently use the standard inequality
	\begin{align}\label{eq:ineqfraction1}
		\frac{y-1}{x-1}\geq \frac{y}{x}\quad\text{for }y\geq x> 1,
	\end{align}
	which we also rewrite for the sake of convenience as
	\begin{align}\label{eq:ineqfraction}
		\frac{x-1}{y-1}\le \frac{x}{y}\quad\text{for $y\ge x\geq1$ with $y\neq 1$}.
	\end{align}
	We will also apply the well-known identity
	\begin{align}\label{eq:qbinproduct}
		\bbin{k}{j}\bbin{j}{i}=\bbin{k}{i}\bbin{k-i}{j-i}
	\end{align}
	and the $q$-binomial inversion formula
	\begin{align}\label{eq:qbininversion}
		\sum_{j=i}^k  (-1)^{j-i}b^{\binom{j-i}{2}}\bbin{j}{i}\bbin{k}{j}=\sum_{j=i}^k  (-1)^{k-j}b^{\binom{k-j}{2}}\bbin{j}{i}\bbin{k}{j}=\delta_{ik}.
	\end{align}
	Moreover, we need two properties of the $q$-Pochhammer symbol, which can be found in~\cite[\S~1.8]{Koekoek}, for example. For a real number~$a$ and nonnegative integers $n,k$, we have
	\begin{align}
		\pochq{a}{n+k}=&\,\pochq{a}{n}\pochq{aq^n}{k}\label{eq:pochidentity_indexsum}\\
		\pochq{a}{n-k}=\frac{\pochq{a}{n}}{\pochq{a^{-1}q^{1-n}}{k}}& (-a)^{-k} q^{\binom{k}{2}-nk+k}\quad\text{for $a\neq 0$} \label{eq:pochidentity_indexdiff}.
	\end{align}
	
	\subsection{Affine \texorpdfstring{$q$}{q}-analogs}\label{subsec:feassol_primalLP_affine}
	In this subsection, we look at $\bil$, $\her$, and $\alt$, where $d$ is required to be odd in the case of $\her$. In all these three cases, the inner distribution of a code whose size equals~\eqref{eq:LPbound_matrixschemes} was computed in~\cite{DelsarteBilinear}, \cite{DelGoe1975}, and~\cite{SchmidtHermitian}. Moreover, it was shown therein that a $d$-code of size~\eqref{eq:LPbound_matrixschemes} is an $(n-d+1)$-design. Here, we use these results to determine the corresponding dual distributions and show that both distributions are nonnegative, implying that the inner distribution is a feasible solution of the primal LP~\eqref{eq:LPprimalCodes}. It actually suffices to look at $\alt$ with even $m$ and odd $q$ since for all the other cases, there are known constructions of codes in $\bil$, $\her$, and $\alt$ whose sizes equal the respective stated LP optimum~\eqref{eq:LPbound_matrixschemes}, see \cite{DelsarteBilinear}, \cite{DelGoe1975}, and \cite{SchmidtHermitian}. Nevertheless, we will give a proof for the nonnegativity of the inner and dual distribution without using the known constructions. The case of $\her$ with even $d$ is handled in Section~\ref{subsec:feassol_primalLP_Hermat}.
	
	The main result of this subsection is the following proposition.
	\begin{prop}\label{prop:fsol_primal_matrixschemes}
		There exists a feasible solution of the primal LP~\eqref{eq:LPprimalCodes} for \mbox{$d$-codes} with $1\leq d\leq n$ in $\bil$, $\her$, and $\alt$ with objective function value~\eqref{eq:LPbound_matrixschemes}, where $d$ is required to be odd in the case of $\her$ and $n=\lfloor m/2\rfloor$ in the case of $\alt$.
	\end{prop}
	Observe that Proposition~\ref{prop:fsol_primal_matrixschemes} and~\ref{prop:fsol_dual_matrixschemes} together with Theorem~\ref{thm:strongduality} imply the first part of Theorem~\ref{thm:LPopt_matrix_qJohn_herm_hyp} \ref{LPoptaffineq}. To prove Proposition~\ref{prop:fsol_primal_matrixschemes}, we first derive the dual distribution of a code of size~\eqref{eq:LPbound_matrixschemes}.
	\begin{prop}\label{prop:distr_matrixschemes}
		Let $Y$ be a $d$-code with $1\leq d\leq n$ in $\bil$, $\her$, or $\alt$ of size~\eqref{eq:LPbound_matrixschemes}, where $d$ is required to be odd in the case of $\her$ and $n=\lfloor m/2\rfloor$ in the case of $\alt$. Then the inner distribution $(A_i)$ of $Y$ satisfies
		\[
		A_{n-i}=\sum_{j=i}^{n-d} (-1)^{j-i} b^{\binom{j-i}{2}}\bbin{j}{i}\bbin{n}{j} ((cb^n)^{n-d+1-j}-1)
		\]
		for all $i=0,1,\dots,n-1$, and the dual distribution $(A'_k)$ of $Y$ satisfies
		\[
		A'_{n-i}=(cb^n)^{n-d+1}\sum_{j=i}^{d-2}(-1)^{j-i} b^{\binom{j-i}{2}} \bbin{j}{i}\bbin{n}{j}((cb^n)^{d-1-j}-1)
		\]
		for all $i=0,1,\dots,n-1$. In particular, $Y$ is an $(n-d+1)$-design.
	\end{prop}
	\begin{proof}
		Let $(A_i)$ and $(A'_k)$ be the inner and dual distribution of $Y$, respectively. Then $(A_i)$ was determined in~\cite[Thm.~5.6]{DelsarteBilinear}, \cite[Thm.~4]{DelGoe1975}, and~\cite[Thm.~3]{SchmidtHermitian} and moreover, it was shown that $Y$ is an $(n-d+1)$-design. It remains to compute $(A'_k)$. By using~\eqref{eq:dualdistr}, $Q_k(i)=P_k(i)$, and \eqref{eq:Pnumidentity_affineschemes}, we obtain
		\begin{align*}
			\sum_{k=0}^{j} \bbin{n-k}{n-j}A'_k
			=\sum_{i=0}^n A_i \sum_{k=0}^j \bbin{n-k}{n-j} P_k(i)
			=(cb^n)^j\sum_{i=0}^{n}\bbin{n-i}{j}A_i
		\end{align*}
		for all $j=0,1,\dots,n$. We have $A_1=A_2=\cdots=A_{d-1}=0$ and $A'_1=A'_2=\cdots=A'_{n-d+1}=0$.
		Because of $A_0=1$, $A'_0=|Y|$, and $\bbin{n-i}{j}=0$ if $i\geq d$ and $j\geq n-d+2$, we obtain
		\begin{align*}
			\sum_{k=n-d+2}^{j}\bbin{n-k}{n-j}A'_k=\bbin{n}{j}((cb^n)^j-(cb^n)^{n-d+1})
		\end{align*}
		for all $j=n-d+2,\dots,n$. Interchanging the order of summation gives
		\begin{align*}
			\sum_{k=j}^{d-2}\bbin{k}{j}A'_{n-k}=(cb^n)^{n-d+1}\bbin{n}{j}((cb^n)^{d-1-j}-1)
		\end{align*}
		for all $j=0,1,\dots,d-2$. Applying the $q$-binomial inversion formula~\eqref{eq:qbininversion} implies the desired expression of $A'_{n-k}$.
	\end{proof}
	
	We need the following lemma to show that both distributions $(A_i)$ and $(A'_k)$ are nonnegative.
	\begin{samepage}
	\begin{lem}
		Let $q\geq 2$ be an integer and $b=-q$.
		\begin{enumerate}[(a)]
			\item For all nonnegative integers $n,i,j$ with $n-i\geq j+2$, we have
			\begin{align}\label{eq:frac_qbin}
				\frac{\big\lvert\bbin{n-i}{j}\big\rvert}
				{\big\lvert\bbin{n-i}{j+2}\big\rvert}
				\geq q^{-2n+4j+2i+2}.
			\end{align}
			\item Let $n$ and $i$ be nonnegative integers. If $n-i\geq 1$, then we have
			\begin{align}\label{eq:ineq_bbin1}
				\left\lvert \bbin{n-i}{1} \right\rvert\leq q^{n-i-1}.
			\end{align}
			If $n-i\geq 2$, then we have
			\begin{align}\label{eq:ineq_bbin2}
				\left\lvert \bbin{n-i}{2} \right\rvert\leq\frac 13 q^{2n-2i-2}.
			\end{align}
		\end{enumerate}
	\end{lem}
	\end{samepage}
	\begin{proof}
		\begin{enumerate}[(a)]
			\item The required bound follows from the inequality
			\begin{align*}
				\frac{\big\lvert\bbin{n-i}{j}\big\rvert}
				{\big\lvert\bbin{n-i}{j+2}\big\rvert}
				=\Bigg\lvert \frac{((-q)^{j+2}-1)((-q)^{j+1}-1)}{((-q)^{n-i-j}-1)((-q)^{n-i-j-1}-1)}\Bigg\rvert 
				\geq\frac{(q^{j+2}+1)\,(q^{j+1}-1)}{(q^{n-i-j}-1)\,(q^{n-i-j-1}+1)}.
			\end{align*}
			\item For $n-i\geq 1$, by using~\eqref{eq:ineqfraction1}, we have
			\begin{align*}
				\left\lvert \bbin{n-i}{1} \right\rvert=\frac{|(-q)^{n-i}-1|}{q+1}\leq\frac{q^{n-i}+1}{q+1}\leq q^{n-i-1}.
			\end{align*}
			For $n-i\geq 2$, by again using~\eqref{eq:ineqfraction1}, we obtain
			\begin{align*}
				\left\lvert \bbin{n-i}{2} \right\rvert\leq \frac{(q^{n-i}-1)(q^{n-i-1}+1)}{(q+1)(q^2-1)}\leq \frac 13 q^{2n-2i-2},
			\end{align*}
			as wanted.
		\end{enumerate}
	\end{proof}
	
	We can now prove Proposition~\ref{prop:fsol_primal_matrixschemes}.
	
	\begin{proof}[Proof of Proposition~\ref{prop:fsol_primal_matrixschemes}]
		Let $1\leq d\leq n$. For $d=1$, the set of all matrices in the respective affine scheme is a $1$-code and thus, there exists a feasible solution of the primal LP with the required objective function value. Assume now that $d\geq 2$. Let $(A_i)$ and $(A'_k)$ be given in Proposition~\ref{prop:distr_matrixschemes}. We will show that all entries of $(A_i)$ and $(A'_k)$ are nonnegative, which implies that $(A_i)$ is a feasible solution of the primal LP~\eqref{eq:LPprimalCodes}. First, we rewrite $A_{n-i}$ by applying~\eqref{eq:qbinproduct} and interchanging the order of summation to obtain for all $i=0,1,\dots,n-d$ that
		\[
		A_{n-i}=\bbin{n}{i}\sum_{j=0}^{n-d-i} (-1)^j b^{\binom{j}{2}}\bbin{n-i}{j} ((cb^n)^{n-d+1-j-i}-1).
		\]
		
		We start with $\bil$ and $\alt$. Observe that it suffices to show that the inner distribution is nonnegative because by taking $n-d+2$ instead of~$d$, the dual distribution becomes a positive multiple of the inner distribution. For $n=d$, we immediately have $A_n\geq 0$. Assume now that $2\leq d<n$. Write
		\[
		a_{i,j}=b^{\binom{j}{2}}\bbin{n-i}{j} ((cb^n)^{n-d+1-j-i}-1)
		\]
		for all $i=0,1,\dots,n-d$ and $j=0,1,\dots,n-d-i$. Take $i\in\{0,1,\dots,n-d\}$. For all $j=0,1,\dots,n-d-i-1$, we have
		\begin{align*}
			\frac{a_{i,j}}{a_{i,j+1}}
			=\frac{b^{\binom{j}{2}}\bbin{n-i}{j}((cb^n)^{n-d+1-j-i}-1)}{b^{\binom{j+1}{2}}\bbin{n-i}{j+1}((cb^n)^{n-d-j-i}-1)}
			=b^{-j}\frac{(b^{j+1}-1)}{(b^{n-i-j}-1)}\,\frac{((cb^n)^{n-d+1-j-i}-1)}{((cb^n)^{n-d-j-i}-1)}.
		\end{align*}
		From~\eqref{eq:ineqfraction1} we find
		\[
		\frac{a_{i,j}}{a_{i,j+1}}>cb^{i+j}\geq 1
		\]
		for all $i,j\geq 0$ except for $\alt$ with even~$m$ and $(i,j)=(0,0)$, where we have
		\[
		\frac{a_{0,0}}{a_{0,1}}=\frac{(q^2-1)(q^{(2n-1)(n-d+1)}-1)}{(q^{2n}-1)(q^{(2n-1)(n-d)}-1)}
		\geq\frac{q^2-1}{q}
		>1.
		\]
		This completes the proof for $\bil$ and $\alt$.
		
		Now, consider $\her$ with odd $d\geq 3$. For $n=d$, we immediately have $A_n\geq 0$. Henceforth, assume that $3\leq d<n$. Write
		\[
		a_{i,j}=(-1)^j b^{\binom{j}{2}} \bbin{n}{i}\bbin{n-i}{j} ((cb^n)^{m-j-i}-1)
		\]
		with $m\in\{n-d+1,d-1\}$ for all $i=0,1,\dots,m-1$ and $j=0,1,\dots,m-i-1$. The sign of $a_{i,j}$ is $(-1)^{\binom{j}{2}+ij+j}$ since $(-1)^{(n+1)m}=1$. Therefore, for all $i$, we have $a_{i,2j}\geq 0$ if $j\geq 0$ is even, and $a_{i,2j+1}\geq 0$ if $j+i\equiv 1\pmod 2$. We will show that
		\begin{equation}\label{eq:props_aij}
			\begin{aligned}
				a_{i,0}&\geq |a_{i,2}|\;\text{for all odd $i\geq 1$}\\
				a_{i,0}&\geq |a_{i,1}|+|a_{i,2}|\;\text{for all even $i\geq 0$},\\
				a_{i,2j}&\geq |a_{i,2j+2}|\;\text{for all $i\geq 0$ and all even $j\geq 2$}\\
				a_{i,2j+1}&\geq |a_{i,2j+3}|\;\text{for all $i,j\geq 0$ with $i+j\equiv 1\pmod 2$,}
			\end{aligned}
		\end{equation}
		which proves the nonnegativity of the inner and dual distribution. Take $i\in\{0,1,\dots,m-1\}$. For all $j=0,1,\dots,m-i-3$, by using~\eqref{eq:frac_qbin}, we have
		\begin{align*}
			\frac{|a_{i,j}|}{|a_{i,j+2}|}
			=\left\lvert \frac{b^{\binom{j}{2}}\bbin{n-i}{j}((cb^n)^{m-j-i}-1)}{b^{\binom{j+2}{2}}\bbin{n-i}{j+2}((cb^n)^{m-j-i-2}-1)}\right\rvert
			\geq q^{2j+2i+1}\,\frac{1-q^{-n(m-j-i)}}{1+q^{-n(m-j-i-2)}}.
		\end{align*}
		Since $n\geq 4$ and $3\leq m-j-i\leq m$, we obtain
		\[
		\frac{|a_{i,j}|}{|a_{i,j+2}|}< \frac{47}{50} q^{2j+2i+1}>1
		\]
		for all $i,j$. Let $i$ be even with $0\leq i\leq m-2$. It remains to show that
		\[
		\frac{|a_{i,0}|}{|a_{i,1}|+|a_{i,2}|}\geq 1.
		\]
		By using~\eqref{eq:ineq_bbin1} and \eqref{eq:ineq_bbin2}, we have
		\begin{align*}
			\frac{|a_{i,0}|}{|a_{i,1}|+|a_{i,2}|}
			&=\frac{|(cb^n)^{m-i}-1|}{\big\lvert\bbin{n-i}{1} ((cb^n)^{m-i-1}-1)\big\rvert+q\big\lvert\bbin{n-i}{2} ((cb^n)^{m-i-2}-1)\big\rvert}\\
			&\geq \frac{q^{n(m-i)}-1}{q^{n-i-1}(q^{n(m-i-1)}+1)+\frac 13 q^{2n-2i-1}(q^{n(m-i-2)}+1)}.
		\end{align*}
		This becomes
		\begin{align*}
			\frac{|a_{i,0}|}{|a_{i,1}|+|a_{i,2}|}
			\geq \frac{1-q^{-n(m-i)}}{q^{-i-1}(1+q^{-n(m-i-1)})+\frac 13 q^{-2i-1}(1+q^{-n(m-i-2)})}.
		\end{align*}
		Since $n\geq 4$, we obtain
		\[
		\frac{|a_{i,0}|}{|a_{i,1}|+|a_{i,2}|}>1
		\]
		if $i\leq m-3$, and if $i=m-2$ (where $a_{i,2}$ cannot occur), we have $|a_{i,0}|/|a_{i,1}|>1$. This completes the proof.
	\end{proof}
	
	\subsection{Ordinary \texorpdfstring{$q$}{q}-analogs}\label{subsec:feassol_primalLP_ordinaryq}
	The goal of this subsection is to prove the following proposition.
	\begin{prop}\label{prop:fsol_primal_qJohn_herm_hyp}
		Let $X$ be the set of $n$-spaces in $\qjohnson$ or generators in $\hermitianpol{2n-1}$ or $\hyperbolichalf$, where $n=\lfloor m/2\rfloor$ in the case of $\hyperbolichalf$. Then there exists a feasible solution of the primal LP~\eqref{eq:LPprimalCodes} for $d$-codes in $\qjohnson$, $\hermitianpol{2n-1}$, and $\hyperbolichalf$ with objective function value~\eqref{eq:LPbound_qJohn_herm_hyp} for $1\leq d\leq n$, where $d$ is required to be odd in the case of~$\hermitianpol{2n-1}$.
	\end{prop}
	The first part of Theorem~\ref{thm:LPopt_matrix_qJohn_herm_hyp} \ref{LPoptordinaryq} follows by combining Proposition~\ref{prop:fsol_primal_qJohn_herm_hyp} and~\ref{prop:fsol_dual_qJohn_herm_hyp} together with Theorem~\ref{thm:strongduality}. To prove Proposition~\ref{prop:fsol_primal_qJohn_herm_hyp}, we first derive the inner distribution of a code whose size equals the stated LP optimum~\eqref{eq:LPbound_qJohn_herm_hyp}.
	\begin{prop}\label{prop:innerdistr_qJohn_herm_hyp}
		Let $X$ be the set of $n$-spaces in $\qjohnson$ or generators in $\hermitianpol{2n-1}$ or $\hyperbolichalf$, where $n=\lfloor m/2\rfloor$ in the case of $\hyperbolichalf$. Assume that $Y$ is a $d$-code with $1\leq d\leq n$ in $\qjohnson$, $\hermitianpol{2n-1}$, or~$\hyperbolichalf$ of size~\eqref{eq:LPbound_qJohn_herm_hyp}, where $d$ is required to be odd in the case of $\hermitianpol{2n-1}$. Then $Y$ is an $(n-d+1)$-design, where the ordering of the dual distribution is imposed by~\eqref{eq:Qnum_qJohn_Herm_hyp}---in particular, the second ordering is taken for~$\hermitianpol{2n-1}$. Moreover for all $i=0,1,\dots,n-1$, the inner distribution $(A_i)$ of $Y$ satisfies
		\begin{align}\label{eq:innerdistr_qJohn_Herm_hyp}
			A_{n-i}=\sum_{j=i}^{n-d} (-1)^{j-i} b^{\binom{j-i}{2}} \bbin{j}{i}\bbin{n}{j}\left(\frac{\pochb{qcb^n}{n-j}\pochb{q}{d-1}}{\pochb{qcb^n}{d-1}\pochb{q}{n-j}}-1\right).
		\end{align}
	\end{prop}

	\begin{rem}
		The design property obtained in Proposition~\ref{prop:innerdistr_qJohn_herm_hyp} in the case of $\hermitianpol{2n-1}$ can be rephrased with respect to the standard ordering. Namely, if $Y$ is a $d$-code in $\hermitianpol{2n-1}$ for odd $d$ whose size equals~\eqref{eq:LPbound_qJohn_herm_hyp}, then its dual distribution with respect to the standard ordering satisfies $A'_i=A'_{n-i+1}=0$ for all $i=1,2,\dots,\tfrac{n-d+1}{2}$ if $n-d+1$ is even, and $A'_i=A'_{n-i+1}=0, A'_{n-(n-d)/2}=0$ for all $i=1,2,\dots,\tfrac{n-d}{2}$ if $n-d+1$ is odd.
	\end{rem}
	
	\begin{proof}[Proof of Proposition~\ref{prop:innerdistr_qJohn_herm_hyp}]
		Let $X$ be the set of $n$-spaces in $\qjohnson$ or generators in $\hermitianpol{2n-1}$ or $\hyperbolichalf$. Then we have
		\begin{align}\label{eq:sizeY_qjohn_herm_hyp}
			|Y|=\frac{|X|\pochb{q}{d-1}}{\pochb{qcb^n}{d-1}}.
		\end{align}
		Let $(A_i)$ and $(A'_k)$ denote the inner and dual distribution of~$Y$, respectively, in terms of the orderings imposed by~\eqref{eq:Pnum_qJohn_Herm_hyp} and~\eqref{eq:Qnum_qJohn_Herm_hyp}. From~\eqref{eq:dualdistr} and~\eqref{eq:Qnumidentity_qJohn_Herm_hyp}, we obtain for all $j=0,1,\dots,n$ that
		\begin{align}
			\sum_{k=0}^j b^{k(n-j)}\bbin{n-k}{n-j}\frac{\pochb{qcb^{n-k}}{n-j}}{\pochb{q}{n-j}} A'_k\notag
			&=\sum_{i=0}^n A_i\sum_{k=0}^j b^{k(n-j)}\bbin{n-k}{n-j}\frac{\pochb{qcb^{n-k}}{n-j}}{\pochb{q}{n-j}} Q'_k(i)\notag\\
			&=|X|\sum_{i=0}^n A_i\bbin{n-i}{j}.\label{eq:innerdualidentity_qJohn_Herm_hyp}
		\end{align}
		Since $A_1=\dots=A_{d-1}=0$ and $\bbin{n-i}{n-d+1}=0$ for $i\geq d$, we find from~\eqref{eq:innerdualidentity_qJohn_Herm_hyp} with $j=n-d+1$ that
		\[
		\sum_{k=0}^{n-d+1}b^{k(d-1)}\bbin{n-k}{d-1}\frac{\pochb{qcb^{n-k}}{d-1}}{\pochb{q}{d-1}}A'_k=|X|\bbin{n}{d-1}A_0.
		\]
		Since $A_0=1$ and $A'_0=|Y|$, we obtain
		\begin{align}\label{eq:innerdualidentity_qJohn_Herm_hyp_nd1}
			\sum_{k=1}^{n-d+1} b^{k(d-1)}\bbin{n-k}{d-1}\frac{\pochb{qcb^{n-k}}{d-1}}{\pochb{q}{d-1}}A'_k
			=\bbin{n}{d-1}\left(|X|-\frac{\pochb{qcb^n}{d-1}}{\pochb{q}{d-1}}|Y|\right),
		\end{align}
		where all coefficients of $A'_k$ on the left-hand side of~\eqref{eq:innerdualidentity_qJohn_Herm_hyp_nd1} are positive. Observe that the bracket on the right-hand side of \eqref{eq:innerdualidentity_qJohn_Herm_hyp_nd1} equals zero because of~\eqref{eq:sizeY_qjohn_herm_hyp}. Therefore, we have $A'_1=A'_2=\dots=A'_{n-d+1}=0$, which means that $Y$ is an $(n-d+1)$-design. By using $A'_0=|Y|=|X|\pochb{q}{d-1}/\pochb{qcb^n}{d-1}$ and~\eqref{eq:sizeY_qjohn_herm_hyp}, we then find from~\eqref{eq:innerdualidentity_qJohn_Herm_hyp} that the inner distribution is determined by the equations
		\[
		\bbin{n}{j}\frac{\pochb{qcb^n}{n-j}\pochb{q}{d-1}}{\pochb{qcb^n}{d-1}\pochb{q}{n-j}}=\bbin{n}{j}+\sum_{i=d}^n\bbin{n-i}{j}A_i
		\]
		for all $j=0,1,\dots,n-d$. This can be rewritten as
		\begin{align}\label{eq:systeminnerdistr_qJohn_Herm_hyp}
			\sum_{k=0}^{n-d}\bbin{k}{j}A_{n-k}=\bbin{n}{j}\left(\frac{\pochb{qcb^n}{n-j}\pochb{q}{d-1}}{\pochb{qcb^n}{d-1}\pochb{q}{n-j}}-1\right).
		\end{align}
		Applying the $q$-binomial inversion formula~\eqref{eq:qbininversion} gives the required expression.
	\end{proof}
	
	The derivation of the dual distribution needs the following lemma. To simplify the notation, we set $a=q^{-1}c^{-1}b^{-2n}$.
	\begin{lem}\label{lem:QCinverse}
		Let $C$ and $Q$ be the $(n+1)\times (n+1)$ matrices defined by $Q=(Q'_k(i))_{k,i}$ and $C=(c_{ji})_{j,i}$ with $Q'_k(i)$ as in \eqref{eq:Qnum_qJohn_Herm_hyp} and
		$c_{j,i}=\bbin{n-i}{j}$ for all $i,j,k=0,1,\dots,n$. Then $C$ is invertible and the product $QC^{-1}$ is given by
		\begin{align}\label{eq:productQCinv}
			(QC^{-1})_{k,j}
			=\mu'_k a^k b^{k^2+\binom{j}{2}} (-bc)^j \frac{\pochb{b^{-k}}{j} \pochb{ab^k}{j}\pochb{qb^{n-k}}{k}}{\pochb{b^{-n}}{j}\pochb{q^{-1}b^{1-n}}{j}\pochb{c^{-1}b^{-n}}{k}}
		\end{align}
		for all $k,j=0,1,\dots,n$ with $\mu'_k$ given as $\mu_k$ in Table~\ref{table:ClassicalASvalmul} in the case of $\qjohnson$ and otherwise in Table~\ref{table:Pnum_parameters}.
	\end{lem}
	To prove Lemma~\ref{lem:QCinverse}, we need an extension of the $q$-Pochhammer symbol for negative subscripts defined by
	\[
	(a;q)_{-k}=\prod\limits_{i=1}^k (1-aq^{-i})^{-1}
	\]
	for all positive integers $k$ if $a\neq q,q^2,\dots,q^k$.
	\begin{proof}[Proof of Lemma~\ref{lem:QCinverse}]
		The inverse of $C$ is given by
		\begin{align}\label{eq:Cinverse}
			(C^{-1})_{i,j}=(-1)^{i+j-n}b^{\binom{i+j-n}{2}}\bbin{j}{n-i}
		\end{align}
		for all $i,j=0,1,\dots,n$ since the $q$-binomial inversion formula~\eqref{eq:qbininversion} implies
		\[
		\sum\limits_{i=0}^n \bbin{n-i}{k} (-1)^{i+j-n} b^{\binom{i+j-n}{2}} \bbin{j}{n-i}
		=\delta_{k,j}.
		\]
		Let $k,j\in \{0,1,\dots,n\}$. By substituting~\eqref{eq:Cinverse} and~\eqref{eq:Qnum_qJohn_Herm_hyp}, we have
		\begin{align*}
			(QC^{-1})_{k,j}
			&=\sum_{i=0}^n Q'_k(i)(C^{-1})_{i,j}\\
			&=\sum_{i=0}^n \mu'_k \hyp{b^{-i},b^{-k},a b^k}{b^{-n},c^{-1}b^{-n}}{b;b}
			(-1)^{i+j-n}b^{\binom{i+j-n}{2}}\bbin{j}{n-i}.
		\end{align*}
		Use the definition of the $q$-hypergeometric function and~\eqref{eq:qbin_poch} to obtain
		\begin{align}\label{eq:proofQC_1}
			(QC^{-1})_{k,j}
			&=\mu'_k\sum_{i,\ell\geq 0} \frac{\pochb{b^{-i}}{\ell}\pochb{b^{-k}}{\ell}\pochb{a b^k}{\ell}}{\pochb{b^{-n}}{\ell}\pochb{c^{-1}b^{-n}}{\ell}\pochb{b}{\ell}}
			(-1)^{i+j-n}b^{\binom{i+j-n}{2}+\ell}\bbin{j}{n-i}\notag\\
			&=\mu'_k\sum_{\ell\geq 0} \frac{\pochb{b^{-k}}{\ell}\pochb{a b^k}{\ell}}{\pochb{b^{-n}}{\ell}\pochb{c^{-1}b^{-n}}{\ell}}
			(-1)^\ell b^{\binom{\ell}{2}+\ell} S_\ell,
		\end{align}
		where
		\[
		S_\ell=\sum_{i=0}^n(-1)^{i+j-n}b^{\binom{i+j-n}{2}-i\ell}\bbin{i}{\ell}\bbin{j}{n-i}.
		\]
		Interchanging the order of summation gives
		\begin{align*}
			S_\ell=b^{-n\ell}\sum_{i=0}^n (-1)^{j-i} b^{\binom{j-i}{2}+i\ell} \bbin{n-i}{\ell}\bbin{j}{i}.
		\end{align*}
		To compute this sum, we use the $q$-Chu--Vandermonde identity
		\begin{align}\label{eq:qChu-Vandermonde2}
			\bbin{x+y}{z}=\sum_{i=0}^x b^{i(y-z+i)}\bbin{x}{i}\bbin{y}{z-i}
		\end{align}
		(see~\cite[\S~2, 2.6.3(c)]{GouldenJackson}, for example). Applying the $q$-inversion formula~\eqref{eq:qbininversion} to \eqref{eq:qChu-Vandermonde2} reveals that
		\begin{align*}
			\sum_{k=0}^x (-1)^k b^{\binom{k}{2}}\bbin{x}{k}\bbin{x-k+y}{z}=b^{x(y-z+x)}\bbin{y}{z-x}.
		\end{align*}
		Put $x=i$, $y=n-i$, and $z=n-\ell$ to obtain
		\begin{align*}
			\sum_{k=0}^i (-1)^k b^{\binom{k}{2}} \bbin{i}{k}\bbin{n-k}{n-\ell}=b^{i\ell}\bbin{n-i}{\ell},
		\end{align*}
		which, by again applying the $q$-binomial inversion formula~\eqref{eq:qbininversion}, gives 
		\[
		\sum_{i=0}^j (-1)^{j-i} b^{\binom{j-i}{2}+i\ell} \bbin{n-i}{\ell}\bbin{j}{i}=(-1)^j b^{\binom{j}{2}}\bbin{n-j}{n-\ell}.
		\]
		Therefore, we have
		\[
		S_\ell=(-1)^j b^{\binom{j}{2}-n\ell}\bbin{n-j}{n-\ell}.
		\]
		Substitute into~\eqref{eq:proofQC_1} and interchange the order of summation to obtain
		\[
		(QC^{-1})_{k,j}=\mu'_k b^{j^2-nj}\sum_{\ell=0}^n \frac{\pochb{b^{-k}}{\ell+j}\pochb{a b^k}{\ell+j}}{\pochb{b^{-n}}{\ell+j}\pochb{c^{-1}b^{-n}}{\ell+j}}
		(-1)^{\ell} b^{\binom{\ell}{2}-(n-j-1)\ell} \bbin{n-j}{\ell}.
		\]
		By~\eqref{eq:qbin_poch}, this becomes
		\[
		(QC^{-1})_{k,j}=\mu'_kb^{j^2-nj}\sum_{\ell=0}^n \frac{\pochb{b^{-k}}{\ell+j}\pochb{a b^k}{\ell+j}\pochb{b^{-(n-j)}}{\ell}}{\pochb{b^{-n}}{\ell+j}\pochb{c^{-1}b^{-n}}{\ell+j}\pochb{b}{\ell}} b^\ell.
		\]
		Applying~\eqref{eq:pochidentity_indexsum} gives
		\begin{align*}
			(QC^{-1})_{k,j}
			&=\mu'_k b^{j^2-nj}\frac{\pochb{b^{-k}}{j}\pochb{a b^k}{j}}{\pochb{b^{-n}}{j}\pochb{c^{-1}b^{-n}}{j}}\sum_{\ell=0}^n \frac{\pochb{b^{-k+j}}{\ell}\pochb{a b^{k+j}}{\ell}}{\pochb{c^{-1}b^{-n+j}}{\ell}\pochb{b}{\ell}}b^\ell\\
			&=\mu'_k b^{j^2-nj}\frac{\pochb{b^{-k}}{j}\pochb{a b^k}{j}}{\pochb{b^{-n}}{j}\pochb{c^{-1}b^{-n}}{j}} \hyper{2}{1}{b^{-(k-j)},a b^{k+j}}{c^{-1}b^{-n+j}}{b;b}. 
		\end{align*}
		Using the $q$-Chu--Vandermonde identity of the form
		\[
		\hyper{2}{1}{b^{-i},x}{y}{b;b}=\frac{\pochb{x^{-1}y}{i}}{\pochb{y}{i}}\, x^i
		\]
		(see \cite[Eq.\ (1.11.5)]{Koekoek}, for example) implies
		\[
		(QC^{-1})_{k,j}
		=\mu'_k a^{k-j} b^{k^2-nj} \frac{\pochb{b^{-k}}{j}\pochb{a b^k}{j}\pochb{a^{-1}c^{-1} b^{-n-k}}{k-j}}{\pochb{b^{-n}}{j}\pochb{c^{-1}b^{-n}}{j}\pochb{c^{-1}b^{-n+j}}{k-j}}.
		\]
		From~\eqref{eq:pochidentity_indexsum} we have $\pochb{c^{-1}b^{-n+j}}{k-j}=\pochb{c^{-1}b^{-n}}{k}/\pochb{c^{-1}b^{-n}}{j}$ and from~\eqref{eq:pochidentity_indexdiff} we find
		\[
		(a^{-1}c^{-1}b^{-n-k})_{k-j}=\frac{(a^{-1}c^{-1}b^{-n-k})_k}{(acb^{n+1})_{j}}\,(-a^{-1}c^{-1}b^{-n-k})^{-j} b^{\binom{j}{2}-kj+j}.
		\]
		Substituting these and using $a=q^{-1}c^{-1}b^{-2n}$ give the required value of $(QC^{-1})_{k,j}$.
	\end{proof}
	
	We are now in position to derive the dual distribution.
	\begin{prop}
		Let $X$ be the set of $n$-spaces in $\qjohnson$ or generators in $\hermitianpol{2n-1}$ or $\hyperbolichalf$, where $n=\lfloor m/2\rfloor$ in the case of $\hyperbolichalf$. Assume that $Y$ is a $d$-code with $1\leq d\leq n$ in $\qjohnson$, $\hermitianpol{2n-1}$, or $\hyperbolichalf$ of size~\eqref{eq:LPbound_qJohn_herm_hyp}, where $d$ is required to be odd in the case of $\hermitianpol{2n-1}$. Let $(A'_k)$ be the dual distribution of~$Y$ in terms of the ordering imposed by~\eqref{eq:Qnum_qJohn_Herm_hyp}---in particular, the second ordering is taken for $\hermitianpol{2n-1}$. Then we have 
		\begin{align}\label{eq:dualdistr_qJohn_Herm_hyp}
			A'_{n-k}=c_k\sum_{j=0}^{d-2-k} 
			(-1)^j b^{\binom{n-k-j}{2}}\frac{\pochb{qb^k}{j}}{\pochb{qcb^{2k+1}}{j}}\bbin{n-k}{j}
			\left(1-\frac{\pochb{qb^{k+j}}{d-k-j-1}}{\pochb{qcb^{n+k+j}}{d-k-j-1}}\right)
		\end{align}
		for all $k=0,1,\dots,n-1$, where
		\[
		c_k=\mu'_{n-k} (-q^{-1}b^{-n+1})^{n-k}\frac{\pochb{qb^k}{n-k}\pochb{q^{-1}c^{-1}b^{-n-k}}{n-k}}{\pochb{c^{-1}b^{-n}}{n-k}\pochb{q^{-1}b^{1-n}}{n-k}}
		\]
		with $\mu'_{n-k}$ given as $\mu_{n-k}$ in Table~\ref{table:ClassicalASvalmul} in the case of $\qjohnson$ and otherwise in Table~\ref{table:Pnum_parameters}.
	\end{prop}
	\begin{proof}
		Similarly to the derivation of the inner distribution, we solve a system of linear equations. Because of~\eqref{eq:innerdualidentity_qJohn_Herm_hyp} and $A_0=1$, $A_1=\cdots=A_{d-1}=0$, we have
		\begin{align*}
			\sum_{k=0}^j b^{k(n-j)}\bbin{n-k}{n-j}\frac{\pochb{qcb^{n-k}}{n-j}}{\pochb{q}{n-j}} A'_k=|X| \left(\bbin{n}{j}+\sum_{i=d}^n A_i\bbin{n-i}{j}\right)
		\end{align*}
		for all $j=0,1,\dots,n$. Since $\bbin{n-i}{j}=0$ if $i\geq d$ and $j\geq n-d+1$, we see that
		\begin{align}\label{eq:dualdistr_Hermpol_LS}
			\sum_{k=0}^j b^{k(n-j)}\bbin{n-k}{n-j}\frac{\pochb{qcb^{n-k}}{n-j}}{\pochb{q}{n-j}} A'_k=|X|\bbin{n}{j}
		\end{align}
		for all $j=n-d+1,\dots,n$. Recall from Proposition~\ref{prop:innerdistr_qJohn_herm_hyp} that $Y$ is an $(n-d+1)$-design. Therefore, we have $A'_1=A'_2=\cdots=A'_{n-d+1}=0$. Because of $A'_0=|Y|=|X|\pochb{q}{d-1}/\pochb{qcb^n}{d-1}$, we obtain
		\begin{align}\label{eq:systemdualdistr_qJohn_Herm_hyp}
			\sum_{k=n-d+2}^jb^{k(n-j)}\bbin{n-k}{n-j}\frac{\pochb{qcb^{n-k}}{n-j}}{\pochb{q}{n-j}}A'_k=|X|\bbin{n}{j}\left(1-\frac{\pochb{qcb^n}{n-j}\pochb{q}{d-1}}{\pochb{qcb^n}{d-1}\pochb{q}{n-j}}\right)
		\end{align}
		for all $j=n-d+2,\dots,n$. In the case of the inner distribution, we simply applied the $q$-binomial inversion formula to solve the system~\eqref{eq:systeminnerdistr_qJohn_Herm_hyp}. Here, we need a slightly different approach. Recall that the $Q$-numbers are determined by~\eqref{eq:Qnumidentity_qJohn_Herm_hyp}. Define the $(n+1)\times (n+1)$ matrices $Q=(Q'_k(i))_{k,i}$, $B=(b_{jk})_{j,k}$, and $C=(c_{ji})_{j,i}$ by
		\[
		b_{jk}=b^{k(n-j)}\bbin{n-k}{n-j}\frac{\pochb{qcb^{n-k}}{n-j}}{\pochb{q}{n-j}}
		\quad\text{and}\quad
		c_{ji}=\bbin{n-i}{j}
		\]
		for all $i,j,k=0,1,\dots,n$. Then we can write~\eqref{eq:Qnumidentity_qJohn_Herm_hyp} as $BQ=|X|C$, which implies that the inverse of $B$ is determined by $B^{-1}=\frac{1}{|X|}QC^{-1}$. Multiplication of~\eqref{eq:systemdualdistr_qJohn_Herm_hyp} with $B^{-1}$ gives
		\begin{align*}
			A'_k=\sum_{j=n-d+2}^n (QC^{-1})_{k,j}\bbin{n}{j}\left(1-\frac{\pochb{qcb^n}{n-j}\pochb{q}{d-1}}{\pochb{qcb^n}{d-1}\pochb{q}{n-j}}\right)
		\end{align*}
		for all $k=n-d+2,\dots,n$. Substituting~\eqref{eq:productQCinv} and using~\eqref{eq:qbin_poch} imply
		\begin{align*}
			A'_k=\mu'_k a^k b^{k^2}\frac{\pochb{qb^{n-k}}{k}}{\pochb{c^{-1}b^{-n}}{k}}\sum_{j=n-d+2}^n 
			(b^{n+1}c)^j \frac{\pochb{b^{-k}}{j} \pochb{ab^k}{j}}{\pochb{b}{j}\pochb{q^{-1}b^{1-n}}{j}}
			\left(1-\frac{\pochb{qcb^n}{n-j}\pochb{q}{d-1}}{\pochb{qcb^n}{d-1}\pochb{q}{n-j}}\right).
		\end{align*}
		Using~\eqref{eq:qbin_poch} and interchanging the order of summation give
		\begin{align*}
			A'_k=\mu'_k a^k \frac{\pochb{qb^{n-k}}{k}}{\pochb{c^{-1}b^{-n}}{k}}\sum_{j=0}^{k-n+d-2} 
			(-b^{n+1}c)^{k-j} b^{\binom{k-j}{2}+kj}\frac{\pochb{ab^k}{k-j}}{\pochb{q^{-1}b^{1-n}}{k-j}}\bbin{k}{j}
			\left(1-\frac{\pochb{qcb^n}{n-k+j}\pochb{q}{d-1}}{\pochb{qcb^n}{d-1}\pochb{q}{n-k+j}}\right).
		\end{align*}
		Apply~\eqref{eq:pochidentity_indexdiff} and $a=q^{-1}c^{-1}b^{-2n}$ to obtain
		\[
		A'_{n-k}=c_k\sum_{j=0}^{d-2-k} 
		(-1)^j b^{\binom{n-k-j}{2}}\frac{\pochb{qb^k}{j}}{\pochb{qcb^{2k+1}}{j}}\bbin{n-k}{j}
		\left(1-\frac{\pochb{qcb^n}{k+j}\pochb{q}{d-1}}{\pochb{qcb^n}{d-1}\pochb{q}{k+j}}\right)
		\]
		with $c_k$ as stated in the proposition. The desired expression of $A'_{n-k}$ now follows by using~\eqref{eq:pochidentity_indexsum}.
	\end{proof}
	
	It remains to show that both distributions $(A_i)$ and $(A'_k)$ are nonnegative.
	\begin{prop}\label{prop:innerdistrnonneg_qJohn_hyp_her}
		For $1<d\leq n$, where $d$ is required to be odd in the case of $\hermitianpol{2n-1}$, all entries of the inner distribution $(A_i)$ given in~\eqref{eq:innerdistr_qJohn_Herm_hyp} are nonnegative.
	\end{prop}
	\begin{proof}
		Let $(A_i)$ be given in~\eqref{eq:innerdistr_qJohn_Herm_hyp}. By using~\eqref{eq:pochidentity_indexsum} and~\eqref{eq:qbinproduct}, and interchanging the order of summation, we have
		\[
		A_{n-i}=\sum_{j=0}^{n-d-i} (-1)^j b^{\binom{j}{2}} \bbin{n}{i}\bbin{n-i}{j}\left(\frac{\pochb{qcb^{n+d-1}}{n-i-j-d+1}}{\pochb{qb^{d-1}}{n-i-j-d+1}}-1\right)
		\]
		for all $i=0,1,\dots,n-d$. Set 
		\[
		a_{i,j}=(-1)^jb^{\binom{j}{2}}\bbin{n}{i}\bbin{n-i}{j}\left(\frac{\pochb{qcb^{n+d-1}}{n-i-j-d+1}}{\pochb{qb^{d-1}}{n-i-j-d+1}}-1\right)
		\]
		for all $i=0,1,\dots,n-d$ and $j=0,1,\dots,n-d-i$. 
		
		We begin with $\qjohnson$ and $\hyperbolichalf$.  Observe that $A_n\geq 0$ if $n=d$. Assume now that $n>d$. We will show that the sequence $(|a_{i,j}|)_j$ is decreasing for all $i=0,1,\dots,n-d$, which implies $A_{n-i}\geq 0$. Take $i\in\{0,1,\dots,n-d\}$. For all $j=0,1,\dots,n-d-i-1$, we have
		\begin{align*}
			\frac{|a_{i,j}|}{|a_{i,j+1}|}
			&=\frac{b^{\binom{j}{2}}\bbin{n-i}{j}\left(\frac{\pochb{qcb^{n+d-1}}{n-i-j-d+1}}{\pochb{qb^{d-1}}{n-i-j-d+1}}-1\right)}{b^{\binom{j+1}{2}}\bbin{n-i}{j+1}\left(\frac{\pochb{qcb^{n+d-1}}{n-i-j-d}}{\pochb{qb^{d-1}}{n-i-j-d}}-1\right)}\\
			&=b^{-j}\, \frac{(b^{j+1}-1)}{(b^{n-i-j}-1)}\,\frac{\left(\frac{(qcb^{2n-i-j-1}-1)}{(qb^{n-i-j-1}-1)}\,\frac{\pochb{qcb^{n+d-1}}{n-i-j-d}}{\pochb{qb^{d-1}}{n-i-j-d}}-1\right)}{\left(\frac{\pochb{qcb^{n+d-1}}{n-i-j-d}}{\pochb{qb^{d-1}}{n-i-j-d}}-1\right)}.
		\end{align*}
		Since $(qcb^{2n-i-j-1}-1)/(qb^{n-i-j-1}-1)>1$, we can apply \eqref{eq:ineqfraction1} to obtain
		\begin{align*}
			\frac{|a_{i,j}|}{|a_{i,j+1}|}>b^{-j}\, \frac{(b^{j+1}-1)(qcb^{2n-i-j-1}-1)}{(b^{n-i-j}-1)(qb^{n-i-j-1}-1)}.
		\end{align*}
		Then again from \eqref{eq:ineqfraction1} we find
		\[
		\frac{|a_{i,j}|}{|a_{i,j+1}|}>cb^{i+j+1}\left(1-\frac{1}{b^{j+1}}\right)\geq cb\left(1-\frac 1b\right)\geq 1
		\]
		for all $i,j\geq 0$, as required.
		
		Now, consider $\hermitianpol{2n-1}$ with odd $d\geq 3$. Write
		\[
		a'_{i,j}=(-1)^jb^{\binom{j}{2}}\bbin{n}{i}\bbin{n-i}{j}\qandq \epsilon_{i,j}=\frac{\pochb{b^{n+d}}{n-i-j-d+1}}{\pochb{-b^d}{n-i-j-d+1}},
		\]
		so that $a_{i,j}=a'_{i,j}(\epsilon_{i,j}-1)$. Observe that the signs of $a'_{i,j}$ and $\epsilon_{i,j}$ are $(-1)^{\binom{j}{2}+(i+j)(n-i)}$ and $(-1)^{(n+1)(i+j)}$, respectively, which implies $\sign(a_{i,j})=(-1)^{\binom{j}{2}+ij+j}$. For all $i\geq 0$, we have $a_{i,2j}\geq 0$ for all even $j\geq 0$ and $a_{i,2j+1}\geq 0$ for all~$j$ with $j+i\equiv 1\pmod 2$. If $n=d$, then we immediately obtain $A_n\geq 0$. Assume now that $n>d$. We will show that $a_{i,j}$ satisfies~\eqref{eq:props_aij}. Observe that this will prove the nonnegativity of the inner distribution. 
		
		Take $i\in\{0,1,\dots,n-d\}$. For all $j=0,1,\dots,n-d-i-2$, use~\eqref{eq:frac_qbin} to obtain
		\begin{align}\label{eq:quotajprime}
			\frac{|a'_{i,j}|}{|a'_{i,j+2}|}
			=\frac{q^{\binom{j}{2}}\left\lvert\bbin{n-i}{j}\right\rvert}{q^{\binom{j+2}{2}}\left\lvert\bbin{n-i}{j+2}\right\rvert}
			\geq q^{-2n+2j+2i+1}.
		\end{align}
		For all $j=0,1,\dots,n-d-i$, because of $i+j\leq n-d$ and $3\leq d<n$, we have
		\begin{align*}
			|\epsilon_{i,j}|
			=\prod_{\ell=0}^{n-j-i-d}\frac{q^{n+d+\ell}-(-1)^{n+d+\ell}}{q^{d+\ell}-(-1)^\ell}
			\geq \prod_{\ell=0}^{n-j-i-d}\frac{q^{n+d+\ell}-1}{q^{d+\ell}+1}
			\geq q^{(n-2)(n-j-i-d+1)}\geq 4.
		\end{align*}
		Since $\epsilon_{i,j}$ and $\epsilon_{i,j+2}$ have the same sign, we thus either have
		\[
		\frac{|\epsilon_{i,j}-1|}{|\epsilon_{i,j+2}-1|}=\frac{|\epsilon_{i,j}|+1}{|\epsilon_{i,j+2}|+1}\quad\text{or}\quad\frac{|\epsilon_{i,j}-1|}{|\epsilon_{i,j+2}-1|}=\frac{\epsilon_{i,j}-1}{\epsilon_{i,j+2}-1}.
		\]
		Because of $|\epsilon_{i,j}|\geq |\epsilon_{i,j+2}|\geq 1$, we find in both cases that
		\[
		\frac{|\epsilon_{i,j}-1|}{|\epsilon_{i,j+2}-1|}\geq \frac{|\epsilon_{i,j}|}{2|\epsilon_{i,j+2}|}=\frac12 \frac{(q^{2n-j-i}-(-1)^{j-i})(q^{2n-j-i-1}-(-1)^{j-i-1})}{(q^{n-j-i}-(-1)^{n-j-i-1})(q^{n-j-i-1}-(-1)^{n-j-i})},
		\]
		from which we obtain
		\begin{align*}
			\frac{|\epsilon_{i,j}-1|}{|\epsilon_{i,j+2}-1|}
			\geq \frac12 \frac{(q^{2n-j-i}+1)(q^{2n-j-i-1}-1)}{(q^{n-j-i}-1)(q^{n-j-i-1}+1)}
			\geq \frac12 q^{2n-1}
		\end{align*}
		by using~\eqref{eq:ineqfraction1}. Combining this with~\eqref{eq:quotajprime} gives
		\[
		\frac{|a_{i,j}|}{|a_{i,j+2}|} \geq \frac12 q^{2j+2i}\geq 2
		\]
		for all $(i,j)\neq (0,0)$.
		
		Let $i$ be even from $\{0,1,\dots,n-d-1\}$. It remains to prove that
		\[
		\frac{|a_{i,0}|}{|a_{i,1}|+|a_{i,2}|}\geq 1.
		\]
		Using~\eqref{eq:ineq_bbin1} and~\eqref{eq:ineq_bbin2} gives
		\begin{align*}
			\frac{|a_{i,0}|}{|a_{i,1}|+|a_{i,2}|}
			=\frac{(\epsilon_{i,0}-1)}{\big\lvert\bbin{n-i}{1}\big\rvert \, |\epsilon_{i,1}-1| + q \,\big\lvert\bbin{n-i}{2}\big\rvert\, (\epsilon_{i,2}-1)}
			\geq \frac{\left(1-\frac{1}{\epsilon_{i,0}}\right)}{D}
		\end{align*}
		with
		\[
		D=q^{n-i-1} \, \left(\left\lvert\frac{\epsilon_{i,1}}{\epsilon_{i,0}}\right\rvert+\frac{1}{\epsilon_{i,0}}\right) + \frac 13 q^{2n-2i-1}\, \left(\frac{\epsilon_{i,2}}{\epsilon_{i,0}}-\frac{1}{\epsilon_{i,0}}\right).
		\]
		We have
		\begin{align*}
			\bigg|\frac{\epsilon_{i,1}}{\epsilon_{i,0}}\bigg|
			=\bigg|\frac{-b^{n-i}-1}{b^{2n-i}-1} \bigg|
			\leq\frac{q^{n-i}+1}{q^{2n-i}-1}
			=q^{-n}\, \frac{\left(1+\frac{1}{q^{n-i}}\right)}{\left(1-\frac{1}{q^{2n-i}}\right)}
		\end{align*}
		and
		\begin{align*}
			\frac{\epsilon_{i,2}}{\epsilon_{i,0}}
			=\frac{(-b^{n-i}-1)(-b^{n-i-1}-1)}{(b^{2n-i}-1)(b^{2n-i-1}-1)}
			\leq q^{-n+1}\,\frac{(q^{n-i-1}+1)}{(q^{2n-i}-1)}
			=q^{-2n}\, \frac{\left(1+\frac{1}{q^{n-i-1}}\right)}{\left(1-\frac{1}{q^{2n-i}}\right)}.
		\end{align*}
		We thus obtain
		\begin{align*}
			D\leq q^{-i-1}\frac{\left(1+\dfrac{1}{q^{n-i}}\right)}{\left(1-\dfrac{1}{q^{2n-i}}\right)}+\frac 13 q^{-2i-1} \frac{\left(1+\dfrac{1}{q^{n-i-1}}\right)}{\left(1-\dfrac{1}{q^{2n-i}}\right)}
			-\frac{1}{\epsilon_{i,0}} \left(\frac 13 q^{2n-2i-1}-q^{n-i-1}\right).
		\end{align*}
		From $2n-i\geq 8$, $n-i\geq 4$, and $i\geq 0$, we see that
		\begin{align*}
			q^{-i-1}\frac{\left(1+\dfrac{1}{q^{n-i}}\right)}{\left(1-\dfrac{1}{q^{2n-i}}\right)}+\frac 13 q^{-2i-1} \frac{\left(1+\dfrac{1}{q^{n-i-1}}\right)}{\left(1-\dfrac{1}{q^{2n-i}}\right)}
			<1
		\end{align*}
		and
		\[
		\frac 13 q^{2n-2i-1}-q^{n-i-1}>1.
		\]
		This gives
		\[
		D<1-\frac{1}{\epsilon_{i,0}}
		\]
		and thus, 
		\[
		\frac{|a_{i,0}|}{|a_{i,1}|+|a_{i,2}|}>1
		\]
		for all $i\leq n-d-2$, and in particular, $|a_{i,0}|/|a_{i,1}|>1$ for $i=n-d-1$ since $a_{i,2}$ does not exist in this case. This completes the proof.
	\end{proof}
	
	Using a similar approach as in the preceding proof, we will now show that the dual distribution is also nonnegative.
	\begin{prop}\label{prop:dualdistrnonneg_qJohn_hyp_her}
		For $1<d\leq n$, where $d$ is required to be odd in the case of $\hermitianpol{2n-1}$, all entries of the dual distribution $(A'_k)$ given in~\eqref{eq:dualdistr_qJohn_Herm_hyp} are nonnegative.
	\end{prop}
	\begin{proof}
		Let $(A'_k)$ be given in~\eqref{eq:dualdistr_qJohn_Herm_hyp}. Write
		\[
		b_{k,j}=c_k(-1)^j b^{\binom{n-k-j}{2}}\frac{\pochb{qb^k}{j}}{\pochb{qcb^{2k+1}}{j}}\bbin{n-k}{j}
		\left(1-\frac{\pochb{qb^{k+j}}{d-k-j-1}}{\pochb{qcb^{n+k+j}}{d-k-j-1}}\right)
		\]
		for all $k=0,1,\dots,d-2$ and $j=0,1,\dots,d-2-k$. 
		
		We start with $\qjohnson$ and $\hyperbolichalf$. Observe that for all $k$, the factor $c_k$ in $A'_k$ is nonnegative. We will show that the sequence $(|b_{k,j}|)_j$ is decreasing for all $k=0,1,\dots,d-2$ implying $A'_{n-k}\geq 0$. Take $k\in\{0,1,\dots,d-2\}$. For all $j=0,1,\dots,d-3-k$, we have
		\begin{align*}
			\frac{|b_{k,j}|}{|b_{k,j+1}|}
			&=\frac{b^{\binom{n-k-j}{2}}\frac{\pochb{qb^k}{j}}{\pochb{qcb^{2k+1}}{j}}\bbin{n-k}{j}
				\left(1-\frac{\pochb{qb^{k+j}}{d-k-j-1}}{\pochb{qcb^{n+k+j}}{d-k-j-1}}\right)}{b^{\binom{n-k-j-1}{2}}\frac{\pochb{qb^k}{j+1}}{\pochb{qcb^{2k+1}}{j+1}}\bbin{n-k}{j+1}
				\left(1-\frac{\pochb{qb^{k+j+1}}{d-k-j-2}}{\pochb{qcb^{n+k+j+1}}{d-k-j-2}}\right)}\\
			&=b^{n-k-j-1}\frac{(qcb^{2k+j+1}-1)(b^{j+1}-1)}{(qb^{k+j}-1)(b^{n-k-j}-1)}\,\frac{\left(1-\frac{\pochb{qb^{k+j}}{d-k-j-1}}{\pochb{qcb^{n+k+j}}{d-k-j-1}}\right)}{\left(1-\frac{(qcb^{n+k+j}-1)\pochb{qb^{k+j}}{d-k-j-1}}{(qb^{k+j}-1)\pochb{qcb^{n+k+j}}{d-k-j-1}}\right)}.
		\end{align*}
		Since $(qcb^{n+k+j}-1)/(qb^{k+j}-1)>1$, we obtain
		\begin{align}\label{eq:bound_bkjbkj2_qjohn_herm_hyp}
			\frac{|b_{k,j}|}{|b_{k,j+1}|}\geq b^{n-k-j-1}\frac{(qcb^{2k+j+1}-1)(b^{j+1}-1)}{(qb^{k+j}-1)(b^{n-k-j}-1)}.
		\end{align}
		By using~\eqref{eq:ineqfraction1}, we have
		\[
		\frac{|b_{k,j}|}{|b_{k,j+1}|}\geq cb^{k+j}\geq 1
		\]
		for all $k,j\geq 0$ except for $\hyperbolichalf$ with even $m$ and $k=j=0$, in which case \eqref{eq:bound_bkjbkj2_qjohn_herm_hyp} gives
		\[
		\frac{|b_{0,0}|}{|b_{0,1}|}\geq \frac{q^{2n-2}(q+1)(q^2-1)}{(q^{2n}-1)}
		\geq\frac{q^2-1}{q}>1,
		\]
		as required.
		
		Now, consider $\hermitianpol{2n-1}$ with odd $d\geq 3$. Set
		\begin{align*}
			b'_{k,j}=b^{\binom{n-k-j}{2}}\frac{\pochb{qb^k}{j}}{\pochb{b^{2k+2}}{j}}\bbin{n-k}{j},
			\qquad
			\delta_{k,j}=\frac{\pochb{qb^{k+j}}{d-k-j-1}}{\pochb{b^{n+k+j+1}}{d-k-j-1}},
		\end{align*}
		so that $b_{k,j}=c_k(-1)^j \,b'_{k,j}(1-\delta_{k,j})$.
		For $0\leq k+j\leq d-2$ and $3\leq d\leq n$, we obtain
		\begin{align}\label{eq:bound_deltakj}
			|\delta_{k,j}|\leq \prod_{\ell=0}^{d-k-j-2}\frac{q^{k+j+\ell+2}}{q^{n+k+j+\ell}}=q^{(n-2)(k+j+1-d)}
		\end{align}
		and therefore,
		\begin{align}\label{eq:bound_deltakj2}
			|\delta_{k,j}|\leq \frac{1}{q^{n-2}}<\frac 12.
		\end{align}
		The sign of $c_k$ is $(-1)^{\binom{n-k}{2}}$ and hence, the entry $b_{k,j}$ has the sign~$(-1)^{\binom{j}{2}+kj+j}$. Thus, we have $b_{k,2j}=|b_{k,2j}|$ if $j$~is even and $b_{k,2j+1}=|b_{k,2j+1}|$ if $k+j\equiv 1\pmod 2$. Similarly to the inner distribution, we will show that $(b_{k,j})$ satisfies~\eqref{eq:props_aij}. Observe that this will prove the nonnegativity of $(A'_k)$.
		
		Take $k\in\{0,1,\dots,d-2\}$. For all \mbox{$j=0,1,\dots,d-4-k$}, we have
		\begin{align*}
			\frac{|b'_{k,j}|}{|b'_{k,j+2}|}
			=\frac{q^{\binom{n-k-j}{2}}\left\lvert \frac{\pochb{qb^k}{j}}{\pochb{b^{2k+2}}{j}}\bbin{n-k}{j}\right\rvert}{ q^{\binom{n-k-j-2}{2}}\left\lvert \frac{\pochb{qb^k}{j+2}}{\pochb{b^{2k+2}}{j+2}}\bbin{n-k}{j+2}\right\rvert}
			=q^{2(n-k-j)-3}\frac{\left\lvert (b^{2k+j+2}-1)(b^{2k+j+3}-1)\bbin{n-k}{j}\right\rvert}{\left\lvert (qb^{k+j}-1)(qb^{k+j+1}-1)\bbin{n-k}{j+2}\right\rvert}.
		\end{align*}
		Using \eqref{eq:frac_qbin} and \eqref{eq:ineqfraction1} gives
		\[
		\frac{|b'_{k,j}|}{|b'_{k,j+2}|}
		\geq q^{2(n-k-j)-3} \frac{(q^{2k+j+2}-1)(q^{2k+j+3}+1)}{(q^{k+j+1}+1)(q^{k+j+2}-1)} q^{-2n+4j+2k+2}
		\geq q^{2k+2j}.
		\]
		From \eqref{eq:bound_deltakj2} we see that $\frac 12<1-\delta_{k,j}<\frac 32$, which gives
		\[
		\frac{|1-\delta_{k,j}|}{|1-\delta_{k,\ell}|}>\frac 13
		\]
		for all $j,\ell=0,1,\dots, d-2-k$.
		Hence, we find 
		\[
		\frac{|b_{k,j}|}{|b_{k,j+2}|}\geq \frac 13 q^{2k+2j}>1
		\]
		for all $(k,j)\neq(0,0)$, as required.
		
		Let $k$ be even from $\{0,1,\dots,d-3\}$. It remains to show that
		\[
		\frac{|b_{k,0}|}{|b_{k,1}|+|b_{k,2}|}\geq 1.
		\]
		We have
		\begin{align*}
			\frac{|b_{k,0}|}{|b_{k,1}|+|b_{k,2}|}
			=\frac{\left\lvert b^{\binom{n-k}{2}} (1-\delta_{k,0}) \right\rvert}
			{\left\lvert b^{\binom{n-k-1}{2}} \frac{\pochb{qb^k}{1}}{\pochb{b^{2k+2}}{1}}\bbin{n-k}{1}(1-\delta_{k,1}) \right\rvert
				+
				\left\lvert b^{\binom{n-k-2}{2}} \frac{\pochb{qb^k}{2}}{\pochb{b^{2k+2}}{2}}\bbin{n-k}{2}(1-\delta_{k,2}) \right\rvert},
		\end{align*}
		which becomes
		\begin{align*}
			\frac{|b_{k,0}|}{|b_{k,1}|+|b_{k,2}|}
			=\frac{\lvert 1-\delta_{k,0} \rvert}
			{\frac{q^{-n+k+1}}{(q^{k+1}+1)}\lvert\bbin{n-k}{1}\rvert\,\lvert 1-\delta_{k,1}\rvert
				+
				\frac{q^{-2n+2k+3}(q^{k+2}+1)}{(q^{k+1}+1)(q^{2k+3}+1)}\lvert\bbin{n-k}{2}\rvert\,\lvert 1-\delta_{k,2}\rvert}.
		\end{align*}
		Applying \eqref{eq:ineq_bbin1} and \eqref{eq:ineq_bbin2} gives
		\begin{align*}
			\frac{|b_{k,0}|}{|b_{k,1}|+|b_{k,2}|}
			\geq 
			\frac{(q^{k+1}+1)\lvert 1-\delta_{k,0} \rvert}
				{\lvert1-\delta_{k,1}\rvert
				+
				\frac 13 q^{-k}\left(1+\frac{1}{q^{k+2}}\right)\lvert 1-\delta_{k,2} \rvert}
			\geq 
			\frac{3\lvert 1-\delta_{k,0} \rvert}
				{\lvert1-\delta_{k,1}\rvert
				+
				\frac 13 \left(1+\frac14\right)\lvert 1-\delta_{k,2} \rvert}.
		\end{align*}
		From~\eqref{eq:bound_deltakj} and $k\leq d-3$, we find
		\begin{align*}
			|\delta_{k,0}|&\leq q^{(n-2)(k+1-d)}\leq \frac{1}{q^{2n-4}}\leq \frac{1}{4}\\
			|\delta_{k,1}|&\leq q^{(n-2)(k+2-d)}\leq \frac{1}{q^{n-2}}\leq \frac{1}{2}\\
			|\delta_{k,2}|&\leq q^{(n-2)(k+3-d)}\leq \frac{1}{q^{2n-4}}\leq \frac{1}{4},
		\end{align*}
		where the latter follows from $k\leq d-5$ since otherwise $j=2$ cannot occur. This implies
		\[
		\frac{|b_{k,0}|}{|b_{k,1}|+|b_{k,2}|}
		\geq \frac{3\left(1-\frac 14 \right)}{1+\frac12+\frac 13 \left(1+\frac14\right)^2}
		>1
		\]
		for $k\leq d-5$ and $|b_{k,0}|/|b_{k,1}|>1$ for $k=d-3$ since $b_{k,2}$ does not exist in this case. This completes the proof.
	\end{proof}
	
	We can now prove the existence of a feasible solution of the primal LP with the required objective function value.
	\begin{proof}[Proof of Proposition~\ref{prop:fsol_primal_qJohn_herm_hyp}]
		For $d=1$, take $Y=X$ and thus, there exists a feasible solution of the primal LP~\eqref{eq:LPprimalCodes} with the required objective function value. For $d>1$, combine Proposition~\ref{prop:innerdistrnonneg_qJohn_hyp_her} and~\ref{prop:dualdistrnonneg_qJohn_hyp_her}.
	\end{proof}
	\subsection{Hermitian forms scheme and even \texorpdfstring{$d$}{d}}\label{subsec:feassol_primalLP_Hermat}
	In this subsection, we will prove the following proposition.
	\begin{prop}\label{prop:fsol_primal_hermat_deven}
		There exists a feasible solution of the primal LP~\eqref{eq:LPprimalCodes} for \mbox{$d$-codes} in $\her$ with objective function value~\eqref{eq:LPbound_Hermat_deven} for all even $d$ with $2\leq d\leq n$.
	\end{prop}
	Proposition~\ref{prop:fsol_primal_hermat_deven} and~\ref{prop:fsol_dual_hermat_deven} together with Theorem~\ref{thm:strongduality} imply the second part of Theorem~\ref{thm:LPopt_matrix_qJohn_herm_hyp}~\ref{LPoptaffineq}.
	
	As in the previous sections, we will compute the inner and dual distribution of a $d$-code with even $d$ whose size equals~\eqref{eq:LPbound_Hermat_deven}. Afterwards, we will show that these distributions are nonnegative, implying the existence of a feasible solution of the primal LP. We start with the computation of the inner distribution.
	\begin{prop}\label{prop:innerdistrHermat}
		Let $n$ and $d$ be integers with $2\leq d\leq n$ and even $d$. Assume that $Y$ is a $d$-code in $\her$ of size~\eqref{eq:LPbound_Hermat_deven}. Let $(A_i)$ and $(A'_k)$ be the inner and dual distribution of $Y$, respectively. Then we have
		\begin{multline*}
			A_{n-i}=\sum_{j=i}^{n-d} (-1)^{j-i} b^{\binom{j-i}{2}} \bbin{j}{i}\bbin{n}{j}\\
			\times\Bigg((-1)^j\frac{|Y|}{b^{nj}}-1
			-\frac{b^j-1}{b^{n-d+1}-1}
			\left((-1)^j\frac{|Y|}{b^{nj}}+(-1)^{n+j}b^{n(n-d+1-j)}\right)\Bigg)
		\end{multline*}
		for all $i=0,1,\dots,n-1$. Moreover, we have $A'_2=\cdots=A'_{n-d+2}=0$ and
		\begin{align}\label{eq:dualdistrHermevenFirst}
			A'_1=\frac{b^n-1}{b^{n-d+1}-1}\left((-1)^{n+1}q^{n(n-d+1)}-|Y|\right).
		\end{align}
	\end{prop}
	\begin{proof}
		By using \eqref{eq:dualdistr}, \eqref{eq:Pnum_affineschemes}, and \eqref{eq:Pnumidentity_affineschemes}, we obtain for all $j=0,1,\dots,n$ that
		\begin{align}\label{eq:innerdualidentityHerm}
			\sum_{k=0}^j \bbin{n-k}{n-j} A'_k
			=(-1)^j b^{nj}\sum_{i=0}^n A_i\bbin{n-i}{j}.
		\end{align}
		Put
		\[
		x_k=\bbin{n-k}{d-1}\bbin{n-1}{d-2},
		\qquad
		y_k=\bbin{n-k}{d-2}\bbin{n-1}{d-1}.
		\]
		Using~\eqref{eq:innerdualidentityHerm} with $j=n-d+1$ and $j=n-d+2$ gives
		\begin{align*}
			(-1)^{n+1}\sum_{k=0}^{n-d+2}\left(x_k-y_k\right) A'_k
			=q^{n(n-d+1)}\left(\bbin{n}{d-1}\bbin{n-1}{d-2}+b^n\bbin{n}{d-2}\bbin{n-1}{d-1}\right).
		\end{align*}
		It was shown in~\cite[Proof of Theorem~2]{SchmidtHermitian} that the coefficients of $A'_k$ on the left-hand side are nonnegative. Since $Y$ is of size~\eqref{eq:LPbound_Hermat_deven}, we have
		\[
		\sum_{k=1}^{n-d+2} (x_k-y_k) A'_k=0.
		\]
		Because of $x_1=y_1$ and $x_k\neq y_k$ for all $k=2,3,\dots,n-d+2$, we obtain $A'_2=A'_3=\cdots=A'_{n-d+2}=0$. From~\eqref{eq:innerdualidentityHerm} we thus find
		\begin{align*}
			\bbin{n}{j}|Y|+\bbin{n-1}{j-1}A'_1=(-1)^j b^{nj} \left(\bbin{n}{j}+\sum_{k=0}^{n-d}A_{n-k}\bbin{k}{j}\right)
		\end{align*}
		for all $j=0,1,\dots,n-d$. Applying the $q$-binomial inversion formula~\eqref{eq:qbininversion} gives 
		\begin{align}\label{eq:hermat_innerdistr_comp}
			A_{n-i}=\sum_{j=i}^{n-d} (-1)^{j-i} b^{\binom{j-i}{2}} \bbin{j}{i}
			\left(\bbin{n}{j}\frac{(-1)^j|Y|}{b^{nj}}+\bbin{n-1}{j-1}\frac{(-1)^jA'_1}{b^{nj}}-\bbin{n}{j}\right)
		\end{align}
		for all $i=0,1,\dots,n-d$. It remains to compute $A'_1$. Set $j=n-d+1$ in~\eqref{eq:innerdualidentityHerm} to obtain
		\[
		\bbin{n}{d-1}|Y|+\bbin{n-1}{d-1}A'_1=(-1)^{n+1} q^{n(n-d+1)}\bbin{n}{d-1}.
		\]
		This gives
		\begin{align*}
			A'_1=\frac{b^n-1}{b^{n-d+1}-1}\left((-1)^{n+1} q^{n(n-d+1)}-|Y|\right).
		\end{align*}
		By substituting this into~\eqref{eq:hermat_innerdistr_comp}, we have
		\begin{multline*}
			A_{n-i}=\sum_{j=i}^{n-d} (-1)^{j-i} b^{\binom{j-i}{2}} \bbin{j}{i}
			\Bigg(\bbin{n}{j}\left((-1)^j\frac{|Y|}{b^{nj}}-1\right)\\
			-\bbin{n-1}{j-1}\,\frac{b^n-1}{b^{n-d+1}-1}
			\left((-1)^j\frac{|Y|}{b^{nj}}+(-1)^{n+j}b^{n(n-d+1-j)}\right)\Bigg)
		\end{multline*}
		for all $i=0,1,\dots,n-d$. After some elementary manipulations, the stated expression of the inner distribution follows.
	\end{proof}
	
	We now derive the dual distribution.
	\begin{prop}\label{prop:dualdistrHermat}
		Let $n$ and $d$ be integers with $2\leq d\leq n$ and even $d$. Assume that~$Y$ is a $d$-code in $\her$ of size~\eqref{eq:LPbound_Hermat_deven}. Then the dual distribution $(A'_k)$ of $Y$ satisfies
		\begin{align*}
			A'_{n-k}=\sum_{j=0}^{d-k-3}(-1)^{n-k}b^{\binom{j}{2}+n(n-j-k)}\bbin{n}{k}\bbin{n-k}{j} (1-\delta_{k,j})
		\end{align*}
		for all $k=0,1,\dots,n-2$, where
		\[
		\delta_{k,j}=(-1)^{n-j-k}\frac{|Y|}{b^{n(n-j-k)}}-(-1)^{n-j-k}\frac{(b^{n-j-k}-1)}{(b^{n-d+1}-1)}\,\frac{(|Y|+(-1)^n q^{n(n-d+1)})}{b^{n(n-j-k)}}
		\]
		for all $j=0,1,\dots,d-k-3$, and $A'_1$ is given in \eqref{eq:dualdistrHermevenFirst}.
	\end{prop}
	\begin{proof}
		Let $(A_i)$ and $(A'_k)$ denote the inner and dual distribution of $Y$, respectively. By Proposition~\ref{prop:innerdistrHermat}, we have $A'_2=\cdots=A'_{n-d+2}=0$. Since $Y$ is a $d$-code, we also have $A_1=\cdots=A_{d-1}=0$. Combine this with~\eqref{eq:innerdualidentityHerm} to obtain
		\[
		\sum_{k=n-d+3}^j \bbin{n-k}{n-j}A'_k=\bbin{n}{j}((-1)^j b^{nj}-|Y|)-\bbin{n-1}{j-1}A'_1
		\]
		for all $j=n-d+3,\dots,n$. Changing the index $j$ and interchanging the order of summation give
		\[
		\sum_{k=j}^{d-3}\bbin{k}{j}A'_{n-k}=\bbin{n}{j}\left((-1)^{n-j}b^{n(n-j)}-|Y|\right)-\bbin{n-1}{j}A'_1
		\]
		for all $j=0,1,\dots,d-3$. By applying the $q$-binomial inversion formula~\eqref{eq:qbininversion}, we have
		\[
		A'_{n-k}=\sum_{j=k}^{d-3}(-1)^{j-k} b^{\binom{j-k}{2}} \bbin{j}{k}\Bigg(
		\bbin{n}{j}\left((-1)^{n-j}b^{n(n-j)}-|Y|\right)-\bbin{n-1}{j} A'_1
		\Bigg)
		\]
		for all $k=0,1,\dots,d-3$. Using \eqref{eq:dualdistrHermevenFirst} gives
		\begin{multline*}
			A'_{n-k}=\sum_{j=k}^{d-3} (-1)^{j-k} b^{\binom{j-k}{2}} \bbin{j}{k}\bbin{n}{j} \\
			\times\left((-1)^{n-j} b^{n(n-j)}-|Y|-\frac{b^{n-j}-1}{b^{n-d+1}-1}\left((-1)^{n+1}q^{n(n-d+1)}-|Y|\right)\right).
		\end{multline*}
		Interchanging the order of summation implies the stated expression of the dual distribution.
	\end{proof}
	
	To show that the inner and dual distribution given in Proposition~\ref{prop:innerdistrHermat} and~\ref{prop:dualdistrHermat}, respectively, are nonnegative, we first bound the size of the corresponding code.
	\begin{lem}\label{lem:bounds_gamma}
		Let $n$, $d$, and $q$ be integers, where $q\geq 2$, $2\leq d\leq n$, and $d$ is even. Let $Y$ be a $d$-code in $\her$ of size~\eqref{eq:LPbound_Hermat_deven}. Then we have
		\begin{align}\label{eq:bounds_gamma}
			\frac 13 q^{n(n-d+2)-1}\leq |Y|\leq \frac 12 q^{n(n-d+2)}.
		\end{align}
	\end{lem}
	\begin{proof}
		For odd $n$, we have
		\begin{align}\label{eq:gamma_nodd}
			|Y|=q^{n(n-d+1)}\;\frac{q^{n-d+2}+1+q^n(q^{n-d+1}-1)}{q^{n-d+2}+q^{n-d+1}},
		\end{align}
		and for even $n$,
		\begin{align}\label{eq:gamma_neven}
			|Y|=q^{n(n-d+1)}\;\frac{-q^{n-d+2}+1+q^n(q^{n-d+1}+1)}{q^{n-d+2}+q^{n-d+1}}.
		\end{align}
		Hence, we obtain
		\begin{align*}
			|Y|
			\leq q^{n(n-d+1)}\;\frac{q^{2n-d+1}+q^n-q^{n-d+2}+1}{q^{n-d+1}(q+1)}
			\leq \frac 13 q^{n(n-d+1)}\;\frac{q^{2n-d+1}+q^n}{q^{n-d+1}}.
		\end{align*}
		This gives
		\[
		|Y|
		\leq \frac 13 q^{n(n-d+2)}\left(1+\frac{1}{q^{n-d+1}}\right)
		\leq \frac 12 q^{n(n-d+2)}.
		\]	
		From~\eqref{eq:gamma_nodd} and~\eqref{eq:gamma_neven}, we see that
		\[
		|Y|
		\geq q^{n(n-d+1)}\;\frac{q^{2n-d+1}-q^n+q^{n-d+2}+1}{q^{n-d+1}(q+1)}.
		\]
		We have
		\[
		\frac{q^{2n-d+1}-q^n+q^{n-d+2}+1}{q^{n-d+1}(q+1)}
		\geq \frac{q^{2n-d+1}-q^n}{q^{n-d+1}(q+1)}
		=q^{n-1}\frac{\left(1-\frac{1}{q^{n-d+1}}\right)}{\left(1+\frac1q\right)}
		\geq \frac 13 q^{n-1}.
		\]
		This implies the required lower bound on $|Y|$.
	\end{proof}
	
	To prove the nonnegativity of the inner distribution $(A_i)$ given in Proposition~\ref{prop:innerdistrHermat}, we rewrite it by doing some elementary manipulations and interchanging the order of summation to obtain
	\begin{align}\label{eq:innerdistr_hermat_deven_2}
		A_{n-i}=\sum_{j=0}^{n-d-i} a_{i,j}
	\end{align}
	with $a_{i,j}=a'_{i,j}(1-\epsilon_{i,j})$, where
	\begin{align}\label{eq:aij'_Hermat}
		a'_{i,j}=(-1)^i b^{\binom{j}{2}-n(j+i)}\bbin{n}{i}\bbin{n-i}{j}|Y|
	\end{align}
	and
	\[
	\epsilon_{i,j}=(-1)^{j+i}\frac{b^{n(j+i)}}{|Y|}+\frac{b^{j+i}-1}{b^{n-d+1}-1}\left(1+(-1)^n\frac{q^{n(n-d+1)}}{|Y|}\right)
	\]
	for all $i=0,1,\dots,n-d$ and $j=0,1,\dots,n-d-i$. We start by deriving bounds on $1-\epsilon_{i,j}$.
	\begin{lem}\label{lem:bounds_errorterm_Hermat}
		Let $q$, $n$, and $d$ be integers, where $q\geq 2$, $2\leq d<n$, and $d$ is even. For $0\leq i+j\leq n-d$, we have
		\begin{align*}
			1-\epsilon_{i,j}\geq\begin{cases}
				\frac{29}{32}&\text{for $n$ and $i+j$ odd}\\
				\frac{125}{128} &\text{for $n$ and $i+j$ even}\\
				\frac{186}{256}&\text{for odd $n$ and even $i+j$}\\
				\frac{31}{64}&\text{for even $n$ and odd $i+j$}
			\end{cases}
		\end{align*}
		and
		\begin{align*}
			1-\epsilon_{i,j}\leq\begin{cases}
				2&\text{for $n$ and $i+j$ odd}\\
				\frac{27}{16}&\text{for $n$ and $i+j$ even}\\
				1&\text{for odd $n$ and even $i+j$}\\
				\frac{2051}{2048}&\text{for even $n$ and odd $i+j$}.
			\end{cases}
		\end{align*}
		Moreover, for odd $n\geq 3$, we have
		\[
		1-\epsilon_{0,0}\geq \frac{253}{256}.
		\]
	\end{lem}
	\begin{proof}
		We have to look at four different cases depending on the parity of $n$ and~$i+j$. 
		
		Assume that $n\geq 3$ and $i+j$ are odd, hence $1\leq i+j\leq n-d$. We have
		\begin{align}\label{eq:Herm_eps_oo}
			1-\epsilon_{i,j}
			= 1-\frac{q^{n(j+i)}}{|Y|}+\frac{q^{j+i}+1}{q^{n-d+1}-1}\left(1-\frac{q^{n(n-d+1)}}{|Y|}\right).
		\end{align}
		By using the lower bound from~\eqref{eq:bounds_gamma}, we obtain
		\begin{align*}
			1-\epsilon_{i,j}
			\geq 1-\frac{3}{q^{n(n-d-j-i+2)-1}}
			\geq \frac{29}{32}.
		\end{align*}
		Moreover, applying the lower bound from~\eqref{eq:bounds_gamma} to~\eqref{eq:Herm_eps_oo} gives
		\[
		1-\epsilon_{i,j}
		\leq 1+\frac{q^{j+i}+1}{q^{n-d+1}-1}
		\leq 1+\frac{q^{n-d}+1}{q^{n-d+1}-1}
		\leq 2,
		\]
		where the last inequality holds since $n-d\geq 1$.
		
		Assume that $n\geq 4$ and $i+j$ are even, thus $0\leq i+j\leq n-d$. We have
		\begin{align}\label{eq:Herm_eps_ee}
			1-\epsilon_{i,j}
			= 1-\frac{q^{n(j+i)}}{|Y|}+\frac{q^{j+i}-1}{q^{n-d+1}+1}\left(1+\frac{q^{n(n-d+1)}}{|Y|}\right).
		\end{align}
		Because of the lower bound from~\eqref{eq:bounds_gamma}, we obtain
		\[
		1-\epsilon_{i,j}
		\geq 1-\frac{3}{q^{n(n-d-j-i+2)-1}}
		\geq \frac{125}{128}.
		\]
		Applying the lower bound from~\eqref{eq:bounds_gamma} to~\eqref{eq:Herm_eps_ee} gives
		\[
		1-\epsilon_{i,j}
		\leq  1+\frac{q^{j+i}-1}{q^{n-d+1}+1}\left(1+\frac{3}{q^{n-1}}\right)
		\leq \frac{27}{16}.
		\]
		
		Assume that $n\geq 3$ is odd and $i+j$ is even, thus $0\leq i+j\leq n-d-1$. Then we have
		\begin{align*}
			1-\epsilon_{i,j}
			= 1-\frac{q^{n(j+i)}}{|Y|}-\frac{q^{j+i}-1}{q^{n-d+1}-1}\left(1-\frac{q^{n(n-d+1)}}{|Y|}\right).
		\end{align*}
		The bounds from~\eqref{eq:bounds_gamma} imply
		\[
		1-\epsilon_{i,j}
		\geq 1-\frac{3}{q^{n(n-d-j-i+2)-1}}-\frac{q^{j+i}-1}{q^{n-d+1}-1}.
		\]
		Applying~\eqref{eq:ineqfraction} and $i+j\leq n-d-1$ gives
		\[
		1-\epsilon_{i,j}
		\geq 1-\frac{3}{q^{3n-1}}-\frac{1}{q^2}
		\geq \frac{186}{256}.
		\]
		Moreover, we obtain
		\[
		1-\epsilon_{0,0}
		=1-\frac{1}{|Y|}
		\geq 1-\frac{3}{q^{n(n-d+2)-1}}
		\geq \frac{253}{256}.
		\]
		Since $q^{n(n-d+1)}/|Y|<1$, we also have
		\[
		1-\epsilon_{i,j}
		\leq 1.
		\]
		
		Assume that $n\geq 4$ is even and $i+j$ is odd, hence $1\leq i+j\leq n-d-1$. Then we have
		\begin{align*}
			1-\epsilon_{i,j}
			= 1+\frac{q^{n(j+i)}}{|Y|}-\frac{q^{j+i}+1}{q^{n-d+1}+1}\left(1+\frac{q^{n(n-d+1)}}{|Y|}\right).
		\end{align*}
		Using the lower bound from~\eqref{eq:bounds_gamma} gives
		\begin{align*}
			1-\epsilon_{i,j}
			&\geq 1-\frac{q^{j+i}+1}{q^{n-d+1}+1}\left(1+\frac{3}{q^{n-1}}\right)\\
			&\geq 1-\frac{1}{q^2}\left(1+\frac{1}{q^{n-d-1}}\right)\left(1+\frac{3}{q^{n-1}}\right)
			\geq \frac{31}{64}.
		\end{align*}
		Due to \eqref{eq:bounds_gamma}, we also obtain
		\[
		1-\epsilon_{i,j}
		\leq 1+\frac{3}{q^{n(n-d-j-i+2)-1}}
		\leq 1+\frac{3}{q^{3n-1}}
		< \frac{2051}{2048}.
		\]
		This finishes the proof.
	\end{proof}
	
	We can now prove the nonnegativity of the inner distribution.
	\begin{prop}\label{prop:innerdistrHermat_deven_nonneg}
		For all even $d$ with $2\leq d\leq n$, all entries of the inner distribution~$(A_i)$ given in Proposition~\ref{prop:innerdistrHermat} are nonnegative.
	\end{prop}
	\begin{proof}
		Let $(A_i)$ be given in the form of~\eqref{eq:innerdistr_hermat_deven_2}. First, assume that $n=d$. Then $n$ is even and we only have to show that $A_n\geq 0$. We have
		\[
		A_n=|Y|-1\geq \frac 13 q^{2n-1}\geq 0.
		\]
		So, we henceforth assume that $2\leq d<n$. Since $1-\epsilon_{i,j}>0$ for all $i,j$ by Lemma~\ref{lem:bounds_errorterm_Hermat}, the sign of $a_{i,j}$ is $(-1)^{\binom{j}{2}+ij+j}$. Hence, for all $i=0,1,\dots,n-d$, we have $\sign(a_{i,2j})=1$ if $j\geq 0$ is even and $\sign(a_{i,2j+1})=1$ if $j+i\equiv 1\pmod 2$. In all other cases, we have $\sign(a_{i,j})=-1$. We will show that $a_{i,j}$ satisfies~\eqref{eq:props_aij}, which proves the nonnegativity of the inner distribution.
		
		Take $i\in\{0,1,\dots,n-d\}$. For all $j=0,1,\dots,n-d-i-2$, by using~\eqref{eq:frac_qbin}, we obtain 
		\begin{align}\label{eq:frac_aij'_Hermat}
			\frac{|a'_{i,j}|}{|a'_{i,j+2}|}
			=\frac{q^{\binom{j}{2}-n(j+i)}\big\lvert \bbin{n-i}{j} \big\rvert}{q^{\binom{j+2}{2}-n(j+i+2)}\big\lvert \bbin{n-i}{j+2} \big\rvert}
			\geq q^{2j+2i+1}.
		\end{align}
		From Lemma~\ref{lem:bounds_errorterm_Hermat}, we deduce
		\[
		\frac{1-\epsilon_{i,j}}{1-\epsilon_{i,j+2}}\geq\begin{cases}
			\frac{29}{64}&\text{for $n$ and $i+j$ odd}\\
			\frac{125}{216}&\text{for $n$ and $i+j$ even}\\
			\frac{186}{256}&\text{for odd $n$ and even $i+j$}\\
			\frac{992}{2051}&\text{for even $n$ and odd $i+j$}
		\end{cases}
		\]
		for all $i,j$ with $0\leq i+j\leq n-d$. We thus obtain
		\[
		\frac{|a_{i,j}|}{|a_{i,j+2}|}
		=\frac{|a'_{i,j}| (1-\epsilon_{i,j})}{|a'_{i,j+2}|(1-\epsilon_{i,j+2})}
		\geq q^{2j+2i+1}\frac{1-\epsilon_{i,j}}{1-\epsilon_{i,j+2}}
		\geq 1
		\]
		for all $i,j$ with $i+j\geq 0$. 
		
		Let $i\geq 0$ be even. It remains to show that
		\[
		\frac{|a_{i,0}|}{|a_{i,1}|+|a_{i,2}|}\geq 1.
		\]
		By using~\eqref{eq:ineq_bbin1} and~\eqref{eq:ineq_bbin2}, we have
		\begin{align*}
			\frac{|a_{i,0}|}{|a_{i,1}|+|a_{i,2}|}
			&=\frac{1-\epsilon_{i,0}}{q^{-n}\lvert\bbin{n-i}{1}\rvert (1-\epsilon_{i,1})+q^{1-2n}\lvert\bbin{n-i}{2}\rvert (1-\epsilon_{i,2})}\\
			&\geq \frac{1-\epsilon_{i,0}}{q^{-i-1}(1-\epsilon_{i,1})+\frac13 q^{-2i-1}(1-\epsilon_{i,2})}.
		\end{align*}
		Together with Lemma~\ref{lem:bounds_errorterm_Hermat}, this gives 
		\[
		\frac{|a_{i,0}|}{|a_{i,1}|+|a_{i,2}|}\geq 1
		\]
		for all $n\geq 3$ and all even $i\geq 0$ except for odd $n\geq 3$ and $i=0$. In the latter case, by using Lemma~\ref{lem:bounds_errorterm_Hermat}, we obtain
		\begin{align*}
			\frac{|a_{0,0}|}{|a_{0,1}|+|a_{0,2}|}
			&=\frac{1-\epsilon_{0,0}}{q^{-n}\lvert\bbin{n}{1}\rvert (1-\epsilon_{0,1})+q^{1-2n}\lvert\bbin{n}{2}\rvert (1-\epsilon_{0,2})}\\
			&\geq \frac{1-\epsilon_{0,0}}{\frac{(1+q^{-n})}{(q+1)} (1-\epsilon_{0,1})+\frac{(1+q^{-n})}{(q+1)(q^2-1)} (1-\epsilon_{0,2})}
			\geq 1.
		\end{align*}
		This completes the proof.
	\end{proof}
	
	It remains to show the nonnegativity of the dual distribution.
	\begin{prop}\label{prop:dualdistrHermat_deven_nonneg}
		For all even $d$ with $2\leq d\leq n$, all entries of the dual distribution $(A'_k)$ given in Proposition~\ref{prop:dualdistrHermat} are nonnegative.
	\end{prop}
	\begin{proof}
		Let $(A'_k)$ be given in Proposition~\ref{prop:dualdistrHermat}. We have
		\begin{align*}
			A'_1
			= \frac{q^n-(-1)^n}{q^{n-d+1}+(-1)^n}\left(|Y|+(-1)^n q^{n(n-d+1)}\right).
		\end{align*}
		From Lemma~\ref{lem:bounds_gamma}, we find $A'_1\geq 0$. In particular, this proves the proposition for $d=2$. Henceforth, we assume $4\leq d\leq n$. We write
		\begin{align*}
			A'_{n-k}=\sum_{j=0}^{d-k-3} b_{k,j}
		\end{align*}
		with $b_{k,j}=b'_{k,j}(1-\delta_{k,j})$, where $\delta_{k,j}$ is as in Proposition~\ref{prop:dualdistrHermat} and
		\[
		b'_{k,j}=(-1)^{n-k}b^{\binom{j}{2}+n(n-j-k)}\bbin{n}{k}\bbin{n-k}{j}
		\]
		for all $k=0,1,\dots,d-3$ and $j=0,1,\dots,{d-k-3}$. First, we look at $\delta_{k,j}$. For $0\leq k+j\leq d-3$, applying the upper bound from~\eqref{eq:bounds_gamma} gives
		\begin{align*}
			|\delta_{k,j}|
			&\leq \frac{|Y|}{q^{n(n-j-k)}}+\frac{(q^{n-j-k}+1)(|Y|+q^{n(n-d+1)})}{q^{n(n-j-k)}(q^{n-d+1}-1)}\\
			&\leq \frac{q^{n(n-d+2)}}{2q^{n(n-j-k)}}+\frac{(q^{n-j-k}+1)(\frac12 q^{n(n-d+2)}+q^{n(n-d+1)})}{q^{n(n-j-k)}(q^{n-d+1}-1)}.
		\end{align*}
		Using $k+j\leq d-3$ and doing some elementary manipulations imply
		\begin{align*}
			|\delta_{k,j}|
			\leq \frac{1}{2q^n}
			+
			\frac{1}{q^{n-3}}\left(1+\frac{1}{q^{n-d+3}}\right)\left(\frac12+\frac{1}{q^n}\right)
			\leq \frac{89}{256}
		\end{align*}
		for all $n\geq 4$. This gives
		\begin{align}\label{eq:hermat_deltakj}
			\frac{1-|\delta_{k,j}|}{1+|\delta_{k,\ell}|}\geq \frac{167}{345}
		\end{align}
		for all $k,j,\ell\geq 0$ with $k+j\leq d-3$ and $k+\ell\leq d-3$. In particular, we have $1-\delta_{k,j}>0$ and hence, the sign of $b_{k,j}$ is $(-1)^{\binom{j}{2}+kj+j}$. Thus, for all $k=0,1,\dots,d-3$, we have $\sign(b_{k,2j})=1$ if $j$ is even and $\sign(b_{k,2j+1})=1$ if $j+k\equiv 1\pmod 2$. In the other cases, we have $\sign(b_{k,j})=-1$. We will now show that $b_{k,j}$ has the same properties as $a_{k,j}$ in~\eqref{eq:props_aij}. Observe that this will prove the nonnegativity of the dual distribution.
		
		We have $b'_{k,j}=\frac{q^{n^2}}{|Y|}a'_{k,j}$ for $a'_{k,j}$ as in~\eqref{eq:aij'_Hermat}. Because of \eqref{eq:frac_aij'_Hermat}, we conclude
		\[
		\frac{|b_{k,j}|}{|b_{k,j+2}|}
		\geq \frac{|a'_{k,j}|}{|a'_{k,j+2}|}\,\frac{(1-|\delta_{k,j}|)}{(1+|\delta_{k,j+2}|)}
		\geq \frac{167}{345}q^{2j+2k+1}
		>1
		\]
		for all $k,j\geq 0$ with $(k,j)\neq (0,0)$.	Let $k\geq 0$ be even. It remains to show that
		\[
		\frac{|b_{k,0}|}{|b_{k,1}|+|b_{k,2}|}\geq 1.
		\]
		From~\eqref{eq:hermat_deltakj},~\eqref{eq:ineq_bbin1}, and~\eqref{eq:ineq_bbin2}, we have
		\begin{align*}
			\frac{|b_{k,0}|}{|b_{k,1}|+|b_{k,2}|}
			\geq \frac{167}{345}\,\frac{1}{\left(q^{-n}\lvert\bbin{n-k}{1}\rvert+q^{1-2n}\lvert\bbin{n-k}{2}\rvert\right)}
			\geq \frac{167}{345}\,\frac{1}{(q^{-k-1}+\frac 13 q^{-2k-1})}.
		\end{align*}
		This gives
		\[
		\frac{|b_{k,0}|}{|b_{k,1}|+|b_{k,2}|}\geq 1
		\]
		for all $k\geq 2$. Similarly to the estimation of $|\delta_{k,j}|$, we obtain
		\[
		|\delta_{0,0}|
		\leq \frac{1}{2q^{n(d-2)}}+\left(\frac{1}{q^{n(d-3)}}+\frac{1}{q^{n(d-2)}}\right)\left(\frac 12+\frac{1}{q^n}\right)
		\leq \frac{161}{4096}
		\]
		for all $n,d$ with $4\leq d\leq n$. Hence, we have
		\[
		\frac{|b_{0,0}|}{|b_{0,1}|+|b_{0,2}|}
		\geq \frac{1}{(q^{-1}+\frac 13 q^{-1})}\,\frac{(1-\frac{161}{4096})}{(1+\frac{89}{256})}
		\geq 1.
		\]
		This completes the proof.
	\end{proof}
	
	Combining Proposition~\ref{prop:innerdistrHermat_deven_nonneg} and~\ref{prop:dualdistrHermat_deven_nonneg} proves Proposition~\ref{prop:fsol_primal_hermat_deven}.
	
	\subsection{Hermitian polar space and even \texorpdfstring{$d$}{d}}\label{subsec:feassol_primalLP_Hermpol}
	In this subsection, we will prove the following proposition.
	\begin{prop}\label{prop:fsol_primal_hermpol_deven}
		Let $X$ be the set of generators in $\hermitianpol{2n-1}$ and let $d$ be even with $2\leq d\leq n$. Then there exists a feasible solution of the primal LP~\eqref{eq:LPprimalCodes} for $d$-codes in $\hermitianpol{2n-1}$ with objective function value~\eqref{eq:LPbound_hermpol_deven}.
	\end{prop}
	
	Observe that Proposition~\ref{prop:fsol_primal_hermpol_deven} and~\ref{prop:fsol_dual_hermpol_deven} together with Theorem~\ref{thm:strongduality} imply the second part of Theorem~\ref{thm:LPopt_matrix_qJohn_herm_hyp} \ref{LPoptordinaryq}.
	
	We will use the same approach as for odd~$d$ and compute the inner and dual distribution of a $d$-code in $\hermitianpol{2n-1}$ of size~\eqref{eq:LPbound_hermpol_deven}. However, the expressions for the inner and dual distribution are more complicated if $d$ is even. This comes from the fact that for even $d$, similarly to $\her$, the first entry of the dual distribution is nonzero whereas for odd $d$, we obtain an $(n-d+1)$-design.
	
	We start by computing the inner distribution. Recall that a $d$-code $Y$ of size~\eqref{eq:LPbound_hermpol_deven} satisfies
	\begin{align}\label{eq:sizeY_hermPol}
		|Y|=|X|\frac{\pochb{q}{d-1}}{\pochb{b^{n+1}}{d-1}}\epsilon(n,d).
	\end{align} 
	
	\begin{prop}\label{prop:innerdistrHermpol_deven}
		Let $X$ be the set of generators in $\hermitianpol{2n-1}$ and let $d$ be even with $2\leq d\leq n$. Assume that $Y$ is a $d$-code in $\hermitianpol{2n-1}$ of size~\eqref{eq:LPbound_hermpol_deven}. Let $(A_i)$ and $(A'_k)$ be the inner and dual distribution of $Y$, respectively, in terms of the orderings imposed by~\eqref{eq:Pnum_qJohn_Herm_hyp} and~\eqref{eq:Qnum_qJohn_Herm_hyp}. Then we have
		\begin{multline*}
			A_{n-i}=\frac{|Y|}{|X|}\sum_{j=i}^{n-d} (-1)^{j-i} b^{\binom{j-i}{2} }\bbin{j}{i}\bbin{n}{j}\frac{\pochb{b^{n+1}}{n-j}}{\pochb{q}{n-j}}\\
			\times\Bigg(1-(1-\epsilon(n,d)^{-1})b^{n-j-d+1}\frac{(b^{n+d-1}-1)(b^j-1)}{(b^{n-d+1}-1)(b^{2n-j}-1)}
			-\frac{\pochb{qb^{d-1}}{n-j-d+1}}{\pochb{b^{n+d}}{n-j-d+1}}\epsilon(n,d)^{-1}
			\Bigg)
		\end{multline*}
		for all $i=0,1,\dots,n-1$. Moreover, we have $A'_2=\cdots=A'_{n-d+2}=0$ and
		\begin{align}\label{eq:dualdistr1_Hermpol}
			A'_1=|X|b^{-d+1}\frac{(q)_{d-1}}{(b^n)_{d-1}}\,\frac{(b^n-1)}{(b^{n-d+1}-1)}\, (1-\epsilon(n,d)).
		\end{align}
	\end{prop}
	\begin{proof}
		Put
		\begin{align*}
				x_k&=b^{k(d-1)+d-2}\frac{\pochb{b^{n-k+1}}{d-1}\pochb{b^n}{d-2}}{\pochb{q}{d-1}\pochb{q}{d-2}}\bbin{n-k}{d-1}\bbin{n-1}{d-2}\\
				y_k&=b^{k(d-2)+d-1}\frac{\pochb{b^{n-k+1}}{d-2}\pochb{b^n}{d-1}}{\pochb{q}{d-2}\pochb{q}{d-1}}\bbin{n-k}{d-2}\bbin{n-1}{d-1}.
		\end{align*}
		Observe that \eqref{eq:innerdualidentity_qJohn_Herm_hyp} also holds for even $d$. Then use~\eqref{eq:innerdualidentity_qJohn_Herm_hyp} with $j=n-d+1$ and $j=n-d+2$ to obtain
		\begin{align}\label{eq:hermpoldeven_1}
			\sum_{k=0}^{n-d+2}(x_k-y_k)A'_k
			=|X|b^{d-2}\frac{\pochb{b^n}{d-2}}{\pochb{q}{d-2}}\left(\bbin{n-1}{d-2}\bbin{n}{d-1}+q\bbin{n-1}{d-1}\bbin{n}{d-2}\frac{b^{n+d-2}-1}{qb^{d-2}-1}\right).
		\end{align}
		It was shown in \cite[Proof of Thm.~3.1]{SchmidtWeissSteiner} that all coefficients of $A'_k$ on the left-hand side of \eqref{eq:hermpoldeven_1} are nonnegative. Using $A'_0=|Y|$ and doing some elementary manipulations give
		\[
		\sum_{k=1}^{n-d+2}(x_k-y_k)A'_k=0.
		\]
		Since $x_1=y_1$ and $x_k\neq y_k$ for all $k\geq 2$, we obtain $A'_2=A'_3=\dots=A'_{n-d+2}=0$. We therefore find from~\eqref{eq:innerdualidentity_qJohn_Herm_hyp} that
		\[
		\bbin{n}{j} \frac{\pochb{b^{n+1}}{n-j}}{\pochb{q}{n-j}}|Y|+b^{n-j}\bbin{n-1}{j-1}\frac{\pochb{b^n}{n-j}}{\pochb{q}{n-j}}A'_1=|X|\left(\bbin{n}{j}+\sum_{k=0}^{n-d}A_{n-k}\bbin{k}{j}\right)
		\]
		for all $j=0,1,\dots,n-d$. Applying the $q$-binomial inversion formula~\eqref{eq:qbininversion} gives
		\begin{align}\label{eq:innerdistr_HermPol_1}
			A_{n-i}=\frac{1}{|X|}\sum_{j=i}^{n-d}(-1)^{j-i}b^{\binom{j-i}{2}} \bbin{j}{i}
			\left(\bbin{n}{j} \frac{\pochb{b^{n+1}}{n-j}}{\pochb{q}{n-j}}|Y|+b^{n-j}\bbin{n-1}{j-1}\frac{\pochb{b^n}{n-j}}{\pochb{q}{n-j}}A'_1-|X|\bbin{n}{j}\right)
		\end{align}
		for all $i=0,1,\dots,n-d$. It remains to compute $A'_1$. Set $j=n-d+1$ in~\eqref{eq:innerdualidentity_qJohn_Herm_hyp} to obtain
		\[
		\bbin{n}{d-1}\frac{\pochb{b^{n+1}}{d-1}}{\pochb{q}{d-1}}|Y|+b^{d-1}\bbin{n-1}{d-1}\frac{\pochb{b^n}{d-1}}{\pochb{q}{d-1}}A'_1=|X|\bbin{n}{d-1}.
		\]
		This implies
		\begin{align*}
			A'_1 =b^{-d+1}\frac{\pochb{q}{d-1}}{\pochb{b^n}{d-1}}\,\frac{\bbin{n}{d-1}}{\bbin{n-1}{d-1}}\left(|X|-\frac{\pochb{b^{n+1}}{d-1}}{\pochb{q}{d-1}}|Y|\right).
		\end{align*}
		Using~\eqref{eq:sizeY_hermPol} gives
		\[
		A'_1=|X| b^{-d+1}\frac{\pochb{q}{d-1}}{\pochb{b^n}{d-1}} \,\frac{(b^n-1)}{(b^{n-d+1}-1)} (1-\epsilon(n,d)).
		\]
		By substituting this into \eqref{eq:innerdistr_HermPol_1} and doing some elementary manipulations, we obtain
		\begin{multline*}
			A_{n-i}=\frac{|Y|}{|X|}\sum_{j=i}^{n-d} (-1)^{j-i} b^{\binom{j-i}{2} }\bbin{j}{i}\bbin{n}{j}\frac{\pochb{b^{n+1}}{n-j}}{\pochb{q}{n-j}}\\
			\times \bigg(1-\frac{|X|}{|Y|}\frac{\pochb{q}{n-j}}{\pochb{b^{n+1}}{n-j}}
			+\frac{|X|}{|Y|} (1-\epsilon(n,d)) b^{n-j-d+1} \,\frac{(b^j-1)(b^n-1)}{(b^{n-d+1}-1)(b^{2n-j}-1)}\,\frac{\pochb{q}{d-1}}{\pochb{b^n}{d-1}}
			\bigg)
		\end{multline*}
		Using \eqref{eq:sizeY_hermPol} and applying \eqref{eq:pochidentity_indexsum} give the desired expression of $A_{n-i}$.
	\end{proof}
	
	We can now compute the dual distribution.
	\begin{prop}\label{prop:dualdistrHermpol_deven}
		Let $X$ be the set of generators in $\hermitianpol{2n-1}$ and let $d$ be even with $2\leq d\leq n$. Assume that $Y$ is a $d$-code in $\hermitianpol{2n-1}$ of size~\eqref{eq:LPbound_hermpol_deven}. Let $(A'_k)$ be the dual distribution of $Y$ in terms of the second ordering imposed by~\eqref{eq:Qnum_qJohn_Herm_hyp}. Then $A'_1$ is given in~\eqref{eq:dualdistr1_Hermpol} and for all $k=0,1,\dots,n-2$, we have
		\begin{align}\label{eq:dualdistr_Hermpol_deven}
			A'_{n-k}=c_k
			\sum_{j=0}^{d-k-3} b^{\binom{j}{2}-nj}\frac{\pochb{b^{-n-1-k}}{n-k-j}}{\pochb{-b^{-n}}{n-k-j}}\bbin{n-k}{j}(1-\delta_{k,j}),
		\end{align}
		where $c_k=\mu'_{n-k} b^{\binom{n-k}{2}-n(n-k)}\frac{\pochb{qb^k}{n-k}}{\pochb{-b^{-n}}{n-k}}$ with $\mu'_{n-k}$ given in Table~\ref{table:Pnum_parameters} and
		\begin{align}\label{eq:delta_Hermpol_deven}
			\delta_{k,j}=\epsilon(n,d)\frac{\pochb{qb^{k+j}}{d-k-j-1}}{\pochb{b^{n+k+j+1}}{d-k-j-1}}
			+(1-\epsilon(n,d))b^{k+j-d+1}\frac{(b^{n-k-j}-1)}{(b^{n-d+1}-1)}\,\frac{\pochb{qb^{k+j}}{d-k-j-1}}{\pochb{b^{n+k+j}}{d-k-j-1}}
		\end{align}
		for all $j=0,1,\dots,d-k-3$. In particular, we have $A'_2=A'_3=\cdots=A'_{n-d+2}=0$.
	\end{prop}
	\begin{proof}
		First, observe that $A'_1, A'_2,\dots,A'_{n-d+2}$ were already determined in Proposition~\ref{prop:innerdistrHermpol_deven} with $A'_2=A'_3=\cdots=A'_{n-d+2}=0$. To obtain $A'_{n-k}$ for $k=0,1,\dots d-3$, we proceed similarly as in the case of odd $d$ and solve a system of linear equations by using the inverse matrix from Lemma~\ref{lem:QCinverse}. As in~\eqref{eq:dualdistr_Hermpol_LS}, we have
		\[
		\sum_{k=0}^j b^{k(n-j)}\bbin{n-k}{n-j}\frac{\pochb{b^{n-k+1}}{n-j}}{\pochb{q}{n-j}} A'_k=|X|\bbin{n}{j}
		\]
		for all  $j=n-d+1,\dots,n$. Use $A'_2=A'_3=\cdots=A'_{n-d+2}=0$ to obtain
		\begin{multline}\label{eq:dualdistrHermpol_LS}
			\sum_{k=n-d+3}^j b^{k(n-j)}\bbin{n-k}{n-j}\frac{\pochb{b^{n-k+1}}{n-j}}{\pochb{q}{n-j}} A'_k\\
			=|X|\bbin{n}{j}-\bbin{n}{j}\frac{\pochb{b^{n+1}}{n-j}}{\pochb{q}{n-j}}|Y|-b^{n-j}\bbin{n-1}{n-j}\frac{\pochb{b^n}{n-j}}{\pochb{q}{n-j}} A'_1
		\end{multline}
		for all $j=n-d+3,\dots,n$. By using~\eqref{eq:sizeY_hermPol} and~\eqref{eq:dualdistr1_Hermpol}, the right-hand side of~\eqref{eq:dualdistrHermpol_LS} can be written as
		\begin{align*}
			|X|\bbin{n}{j}\Bigg(1&-\epsilon(n,d)\frac{\pochb{b^{n+1}}{n-j}\pochb{q}{d-1}}{\pochb{b^{n+1}}{d-1}\pochb{q}{n-j}}
			-(1-\epsilon(n,d))b^{n-j-d+1}\,\frac{(b^j-1)}{(b^{n-d+1}-1)}\,\frac{\pochb{b^n}{n-j}\pochb{q}{d-1}}{\pochb{b^n}{d-1}\pochb{q}{n-j}}\Bigg).
		\end{align*}
		By applying~\eqref{eq:pochidentity_indexsum}, we thus have
		\begin{multline*}
			\sum_{k=n-d+3}^j b^{k(n-j)}\bbin{n-k}{n-j}\frac{\pochb{b^{n-k+1}}{n-j}}{\pochb{q}{n-j}} A'_k
			=|X|\bbin{n}{j} \Bigg(1-\epsilon(n,d)\frac{\pochb{qb^{n-j}}{d-n+j-1}}{\pochb{b^{2n-j+1}}{d-n+j-1}}
			\\
			-(1-\epsilon(n,d))b^{n-j-d+1}\frac{(b^j-1)}{(b^{n-d+1}-1)}\,\frac{\pochb{qb^{n-j}}{d-n+j-1}}{\pochb{b^{2n-j}}{d-n+j-1}}\Bigg).
		\end{multline*}
		Similarly as in the case of odd $d$, we multiply the matrix $QC^{-1}$ from Lemma~\ref{lem:QCinverse}, which gives
		\begin{multline*}
			A'_k=\mu'_k b^{-(2n+1)k+k^2}\frac{\pochb{qb^{n-k}}{k}}{\pochb{-b^{-n}}{k}}
			\sum_{j=n-d+3}^n b^{\binom{j}{2}+j}\frac{\pochb{b^{-k}}{j}\pochb{b^{-2n-1+k}}{j}}{\pochb{b^{-n}}{j}\pochb{-b^{-n}}{j}}\bbin{n}{j}\\
			\times\Bigg(1-\epsilon(n,d)\frac{\pochb{qb^{n-j}}{d-n+j-1}}{\pochb{b^{2n-j+1}}{d-n+j-1}}
			-(1-\epsilon(n,d))b^{n-j-d+1}\frac{(b^j-1)}{(b^{n-d+1}-1)}\,\frac{\pochb{qb^{n-j}}{d-n+j-1}}{\pochb{b^{2n-j}}{d-n+j-1}}\Bigg).
		\end{multline*}
		Apply \eqref{eq:qbin_poch} and shift the index from $k$ to $n-k$ to obtain
		\begin{multline*}
			A'_{n-k}=\mu'_{n-k} b^{-(n-k)(n+k+1)}\frac{\pochb{qb^k}{n-k}}{\pochb{-b^{-n}}{n-k}}
			\sum_{j=n-d+3}^{n-k} b^{\binom{j}{2}+j+kj}\frac{\pochb{b^{-n-1-k}}{j}}{\pochb{-b^{-n}}{j}}\bbin{n-k}{j}\\
			\times\Bigg(1-\epsilon(n,d)\frac{\pochb{qb^{n-j}}{d-n+j-1}}{\pochb{b^{2n-j+1}}{d-n+j-1}}
			-(1-\epsilon(n,d))b^{n-j-d+1}\frac{(b^j-1)}{(b^{n-d+1}-1)}\,\frac{\pochb{qb^{n-j}}{d-n+j-1}}{\pochb{b^{2n-j}}{d-n+j-1}}\Bigg).
		\end{multline*}
		Changing the order of summation and doing some elementary manipulations give the stated expression of $A'_{n-k}$.
	\end{proof}
	
	Proving the nonnegativity of both distributions requires the following bounds on~$\epsilon(n,d)$.
	\begin{lem}\label{lem:eps_hermpol}
		Let $q$, $n$, and $d$ be integers, where $q\geq 2$, $2\leq d\leq n$, and $d$ is even. Then we have
		\begin{align}\label{eq:boundlower_epsHermpol}
			\epsilon(n,d)>\begin{cases}
				-\dfrac{q^{n+d-1}+1}{q^{d-1}-1}&\text{for even $n$}\\[2ex]
				\frac 12 q^{d-2}(q^{n-d+1}-1) &\text{for odd $n$}
			\end{cases}
		\end{align}
		and
		\begin{align}\label{eq:boundupper_epsHermpol}
			\epsilon(n,d)<\begin{cases}
				-\dfrac{q^n}{q+1}&\text{for even $n$}\\[2.5ex]
				\dfrac{q^{n+d-1}-1}{q^{d-1}-1}&\text{for odd $n$.}
			\end{cases}
		\end{align}
	\end{lem}
	\begin{proof}
		Since $d$ is even, we find from~\eqref{eq:epsnd_Hermpol_deven} that
		\[
		(-1)^{n+1}\epsilon(n,d)=\frac{-q^{n-d+2}+(-1)^n+(-1)^n q\,\frac{q^{n+d-2}-(-1)^n}{q^{d-1}-1}(q^{n-d+1}+(-1)^n)}{(-1)^nq^{n-d+2}-1+(-1)^n q\,\frac{q^{n+d-2}-(-1)^n}{q^{n+d-1}+(-1)^n}(q^{n-d+1}+(-1)^n)}.
		\]
		In \cite[Eq.~(3.19)]{SchmidtWeissSteiner}, it was shown that
		\[
		(-1)^{n+1}\epsilon(n,d)<\frac{q^{n+d-1}+(-1)^n}{q^{d-1}-1},
		\]
		which implies the stated lower and upper bound for even $n$ and odd $n$, respectively.
		Assume that $n$ is even. Then we have
		\[
		(-1)^{n+1}\epsilon(n,d)=\frac{-q^{n-d+2}+1+q\frac{q^{n+d-2}-1}{q^{d-1}-1}(q^{n-d+1}+1)}{q^{n-d+2}-1+q\frac{q^{n+d-2}-1}{q^{n+d-1}+1}(q^{n-d+1}+1)}.
		\]
		Using \eqref{eq:ineqfraction1} gives
		\begin{align*}
			\frac{q^{n+d-2}-1}{q^{d-1}-1}\geq q^{n-1}
		\end{align*}
		and we also have
		\[
		\frac{q^{n+d-2}-1}{q^{n+d-1}+1}\leq \frac1q.
		\]
		We thus obtain
		\begin{align*}
			(-1)^{n+1}\epsilon(n,d)&\geq\frac{-q^{n-d+2}+1+q^{2n-d+1}+q^n}{q^{n-d+2}+q^{n-d+1}}
			>\frac{q^{2n-d+1}}{q^{n-d+2}+q^{n-d+1}}
			=\frac{q^n}{q+1},
		\end{align*}
		as stated. Assume that $n$ is odd. We have
		\[
		(-1)^{n+1}\epsilon(n,d)=\frac{q^{n-d+2}+1+q\frac{q^{n+d-2}+1}{q^{d-1}-1}(q^{n-d+1}-1)}{q^{n-d+2}+1+q\frac{q^{n+d-2}+1}{q^{n+d-1}-1}(q^{n-d+1}-1)}.
		\]
		Because of
		\begin{align*}
			\frac{q^{n+d-2}+1}{q^{d-1}-1}\geq q^{n-1}\qandq\frac{q^{n+d-2}+1}{q^{n+d-1}-1}\leq 1,
		\end{align*}
		we obtain
		\begin{align*}
			(-1)^{n+1}\epsilon(n,d)&\geq \frac{q^{n-d+2}+1+q^n(q^{n-d+1}-1)}{q^{n-d+2}+1+q(q^{n-d+1}-1)}.
		\end{align*}
		Observe that it holds
		\[
		q^{n-d+2}+1>\frac 12 (-q^n+q^{n-1}+q^{d-1}-q^{d-2})
		\]
		since $-q^n+q^{n-1}+q^{d-1}\leq 0$. This can be used to show, by elementary manipulations, that
		\begin{align*}
			q^{n-d+2}+1+q^n(q^{n-d+1}-1)
			>\frac 12 q^{d-2}(q^{n-d+1}-1)\left(q^{n-d+2}+1+q(q^{n-d+1}-1)\right),
		\end{align*}
		which implies the required bound.
	\end{proof}
	
	To prove the nonnegativity of the inner distribution given in Proposition~\ref{prop:innerdistrHermpol_deven}, we rewrite it as
	\begin{align}\label{eq:innerdistr_Hermpoleven_rewritten}
		A_{n-i}=\frac{|Y|}{|X|}\sum_{j=0}^{n-d-i} a_{i,j}
	\end{align}
	with $a_{i,j}=a'_{i,j} (1-\epsilon_{i,j})$, where
	\[
	a'_{i,j}=(-1)^j b^{\binom{j}{2}} \bbin{n}{i}\bbin{n-i}{j}\frac{\pochb{b^{n+1}}{n-j-i}}{\pochb{q}{n-j-i}}
	\]
	and
	\begin{align*}
		\epsilon_{i,j}=(1-\epsilon(n,d)^{-1})b^{n-j-i-d+1}\frac{(b^{n+d-1}-1)(b^{j+i}-1)}{(b^{n-d+1}-1)(b^{2n-j-i}-1)}
		+\frac{\pochb{qb^{d-1}}{n-j-i-d+1}}{\pochb{b^{n+d}}{n-j-i-d+1}}\,\epsilon(n,d)^{-1}
	\end{align*}
	for all $i=0,1,\dots,n-d$ and $j=0,1,\dots,n-d-i$. We first give lower and upper bounds on $\epsilon_{i,j}$.
	
	\begin{lem}\label{lem:bounds_errorterm_Hermpol}
		Let $q\geq 2$, $n\geq 4$, and $d$ be integers, where $d$ is even with $2\leq d<n$. For all $i,j$ with $0\leq i+j\leq n-d$, we have
		\begin{align*}
			1-\epsilon_{i,j}>\begin{cases}
				\frac{31}{32}&\text{for $n$ and $i+j$ odd}\\
				\frac{61}{64} &\text{for $n$ and $i+j$ even}\\
				\frac{191}{256} &\text{for odd $n$ and even $i+j$}\\
				\frac{109}{128} &\text{for even $n$ and odd $i+j$}
			\end{cases}
		\end{align*}
		and
		\begin{align*}
			1-\epsilon_{i,j}<\begin{cases}
				2&\text{for $n$ and $i+j$ odd}\\
				\frac{51}{32} &\text{for $n$ and $i+j$ even}\\
				1 &\text{for odd $n$ and even $i+j$}\\
				\frac{259}{256} &\text{for even $n$ and odd $i+j$}.
			\end{cases}
		\end{align*}
		Moreover, for odd $n$, we have
		\[
		1-\epsilon_{0,0}>\frac{255}{256}.
		\]
	\end{lem}
	\begin{proof}
		From Lemma~\ref{lem:eps_hermpol}, we find
		\begin{align}\label{eq:bounds_eps_Hermpol_neven}
			-\frac{q+1}{q^n}<\epsilon(n,d)^{-1}<-\frac{q^{d-1}-1}{q^{n+d-1}+1}\quad\text{for even $n$}
		\end{align}
		and
		\begin{align}\label{eq:bounds_eps_Hermpol_nodd}
			\frac{q^{d-1}-1}{q^{n+d-1}-1}<\epsilon(n,d)^{-1}<2\,\frac{q^{-d+2}}{q^{n-d+1}-1}\quad\text{for odd $n$}.
		\end{align}
		Note that the sign of $(qb^{d-1})_{n-j-i-d+1}/(b^{n+d})_{n-j-i-d+1}$ is $(-1)^{(n+1)(n-j-i+1)}$. Moreover, we have
		\begin{align}\label{eq:epsij_bounds_prod}
			\left\lvert 
			\frac{\pochb{qb^{d-1}}{n-j-i-d+1}}{\pochb{b^{n+d}}{n-j-i-d+1}}
			\right\rvert
			=\left\lvert
			\prod_{\ell=0}^{n-j-i-d}\frac{qb^{d-1+\ell}-1}{b^{n+d+\ell}-1}
			\right\rvert
			\leq \prod_{\ell=0}^{n-j-i-d} \frac{q^{d+\ell}+1}{q^{n+d+\ell}-1}
			\leq q^{-(n-2)(n-j-i-d+1)}
		\end{align}
		for all $i,j$ with $0\leq i+j\leq n-d$. 
		
		Assume that $n\geq 5$ and $i+j$ are odd. We have
		\begin{multline}\label{eq:errorterm_Hermpol_nodd_jiodd}
			1-\epsilon_{i,j}=1+(1-\epsilon(n,d)^{-1})q^{n-j-i-d+1}\frac{(q^{n+d-1}-1)(q^{j+i}+1)}{(q^{n-d+1}-1)(q^{2n-j-i}+1)}\\
			-\frac{(qb^{d-1})_{n-j-i-d+1}}{(b^{n+d})_{n-j-i-d+1}}\epsilon(n,d)^{-1}.
		\end{multline}
		Use that the second summand is nonnegative and
		\[
		\frac{\pochb{qb^{d-1}}{n-j-i-d+1}}{\pochb{b^{n+d}}{n-j-i-d+1}}\geq 0
		\]
		together with the upper bound from~\eqref{eq:bounds_eps_Hermpol_nodd} to obtain
		\begin{align*}
			1-\epsilon_{i,j}
			>1-\frac{(qb^{d-1})_{n-j-i-d+1}}{(b^{n+d})_{n-j-i-d+1}} \frac{2q^{-d+2}}{(q^{n-d+1}-1)}
			>1-\frac{(qb^{d-1})_{n-j-i-d+1}}{(b^{n+d})_{n-j-i-d+1}} 2q^{-n+2}.
		\end{align*}
		Because of~\eqref{eq:epsij_bounds_prod}, we have
		\begin{align*}
			1-\epsilon_{i,j}>1-\frac{2}{q^{(n-2)(n-j-i-d+2)}}.
		\end{align*}
		Since $n\geq 5$ and $i+j\leq n-d$, we deduce
		\[
		1-\epsilon_{i,j}>\frac{31}{32}.
		\]
		From~\eqref{eq:errorterm_Hermpol_nodd_jiodd} we also obtain
		\[
		1-\epsilon_{i,j}<1+(1-\epsilon(n,d)^{-1})q^{n-j-i-d+1}\frac{(q^{n+d-1}-1)(q^{j+i}+1)}{(q^{n-d+1}-1)(q^{2n-j-i}+1)}.
		\]
		Using the lower bound from~\eqref{eq:bounds_eps_Hermpol_nodd} gives
		\begin{align*}
			1-\epsilon_{i,j}
			<1+q^{n-j-i}\frac{(q^n-1)(q^{j+i}+1)}{(q^{n-d+1}-1)(q^{2n-j-i}+1)}
			<1+\frac{q^{j+i}+1}{q^{n-d+1}-1}.
		\end{align*}
		Because of $i+j\leq n-d$ and $n-d\geq 1$, we have
		\[
		1-\epsilon_{i,j}<1+\frac{q^{n-d}+1}{q^{n-d+1}-1}\leq 2,
		\]
		as required.
		
		Assume that $n$ and $j+i$ are even. We have
		\begin{multline}\label{eq:errorterm_Hermpol_neven_jieven}
			1-\epsilon_{i,j}
			=1+(1-\epsilon(n,d)^{-1})q^{n-j-i-d+1}\frac{(q^{n+d-1}+1)(q^{j+i}-1)}{(q^{n-d+1}+1)(q^{2n-j-i}-1)}\\
			-\frac{(qb^{d-1})_{n-j-i-d+1}}{(b^{n+d})_{n-j-i-d+1}}\epsilon(n,d)^{-1}.
		\end{multline}
		The second summand is nonnegative, whereas $\epsilon(n,d)^{-1}$ and
		\[
		\frac{(qb^{d-1})_{n-j-i-d+1}}{(b^{n+d})_{n-j-i-d+1}}
		\]
		are negative. Therefore, we obtain
		\begin{align*}
			1-\epsilon_{i,j}
			&\geq 
			1-\left\lvert\frac{(qb^{d-1})_{n-j-i-d+1}}{(b^{n+d})_{n-j-i-d+1}}\right\rvert\,|\epsilon(n,d)^{-1}|.
		\end{align*}
		Using~\eqref{eq:epsij_bounds_prod} and the lower bound from~\eqref{eq:bounds_eps_Hermpol_neven} gives
		\begin{align*}
			1-\epsilon_{i,j}
			>1-\frac{q+1}{q^{n+(n-2)(n-j-i-d+1)}}
			\geq \frac{61}{64},
		\end{align*}
		where the last inequality follows from $j+i\leq n-d$ and $n\geq 4$. From~\eqref{eq:errorterm_Hermpol_neven_jieven} we also obtain
		\[
		1-\epsilon_{i,j}<1+(1-\epsilon(n,d)^{-1})q^{n-j-i-d+1}\frac{(q^{n+d-1}+1)(q^{j+i}-1)}{(q^{n-d+1}+1)(q^{2n-j-i}-1)}.
		\]
		Using the lower bound from~\eqref{eq:bounds_eps_Hermpol_neven} gives
		\[
		1-\epsilon_{i,j}<1+q^{-j-i-d+1}\frac{(q^n+q+1)(q^{n+d-1}+1)(q^{j+i}-1)}{(q^{n-d+1}+1)(q^{2n-j-i}-1)}.
		\]
		By~\eqref{eq:ineqfraction1} and~\eqref{eq:ineqfraction}, we have
		\[
		1-\epsilon_{i,j}<1+\frac{q^n+q+1}{q^{2n-d-j-i+1}}.
		\]
		Since $i+j\leq n-d$ and $n\geq 4$, we obtain
		\[
		1-\epsilon_{i,j}<1+\frac 1q+\frac{1}{q^n}+\frac{1}{q^{n+1}}\leq \frac{51}{32}.
		\]
		
		Assume that $n\geq 5$ is odd and $i+j$ is even, which implies $i+j\leq n-d-1$. We have
		\begin{multline}\label{eq:errorterm_Hermpol_nodd_jieven}
			1-\epsilon_{i,j}=1-(1-\epsilon(n,d)^{-1})q^{n-j-i-d+1}\frac{(q^{n+d-1}-1)(q^{j+i}-1)}{(q^{n-d+1}-1)(q^{2n-j-i}-1)}\\
			-\frac{(qb^{d-1})_{n-j-i-d+1}}{(b^{n+d})_{n-j-i-d+1}}\,\epsilon(n,d)^{-1}.
		\end{multline}
		Since $(qb^{d-1})_{n-j-i-d+1}/(b^{n+d})_{n-j-i-d+1}\geq 0$, we use the bounds from~\eqref{eq:bounds_eps_Hermpol_nodd} to obtain
		\begin{align}\label{eq:bound_epsij_nodd_ijeven}
			1-\epsilon_{i,j}>1-q^{n-j-i}\frac{(q^n-1)(q^{j+i}-1)}{(q^{n-d+1}-1)(q^{2n-j-i}-1)}
			-2\,\frac{q^{-d+2}}{(q^{n-d+1}-1)}\,\frac{(qb^{d-1})_{n-j-i-d+1}}{(b^{n+d})_{n-j-i-d+1}}.
		\end{align}
		From~\eqref{eq:ineqfraction} we find
		\begin{align*}
			\frac{q^n-1}{q^{2n-j-i}-1}\leq q^{-n+j+i},
			\qquad
			\frac{q^{j+i}-1}{q^{n-d+1}-1}\leq q^{-n+j+i+d-1},
		\end{align*}
		which implies
		\begin{align*}
			1-\epsilon_{i,j}>1-\frac{1}{q^{n-j-i-d+1}}-2\,\frac{q^{-d+2}}{(q^{n-d+1}-1)}\,\frac{(qb^{d-1})_{n-j-i-d+1}}{(b^{n+d})_{n-j-i-d+1}}.
		\end{align*}
		Using~\eqref{eq:epsij_bounds_prod} gives
		\begin{align*}
			1-\epsilon_{i,j}>1-\frac{1}{q^{n-j-i-d+1}}-\frac{2}{q^{(n-2)(n-j-i-d+2)}}.
		\end{align*}
		Because of $n\geq 5$ and $j+i\leq n-d-1$, we have
		\begin{align*}
			1-\epsilon_{i,j}>1-\frac 14-\frac{1}{2^8}=\frac{191}{256}.
		\end{align*}
		In particular, from~\eqref{eq:bound_epsij_nodd_ijeven} we find similarly
		\[
		1-\epsilon_{0,0}>1-\frac{1}{2^8}=\frac{255}{256}.
		\]
		By applying $1-\epsilon(n,d)^{-1}>0$ and 
		\[
		\frac{(qb^{d-1})_{n-j-i-d+1}}{(b^{n+d})_{n-j-i-d+1}}\,\epsilon(n,d)^{-1}>0
		\]
		to~\eqref{eq:errorterm_Hermpol_nodd_jieven}, we obtain $1-\epsilon_{i,j}<1$.
		
		Assume that $n$ is even and $i+j$ is odd, which implies $n-d\geq 2$ and $1\leq i+j\leq n-d-1$. We have
		\begin{multline}\label{eq:errorterm_Hermpol_neven_jiodd}
			1-\epsilon_{i,j}=1-(1-\epsilon(n,d)^{-1})q^{n-j-i-d+1}\frac{(q^{n+d-1}+1)(q^{j+i}+1)}{(q^{n-d+1}+1)(q^{2n+j+i}+1)}\\
			-\frac{(qb^{d-1})_{n-j-i-d+1}}{(b^{n+d})_{n-j-i-d+1}}\,\epsilon(n,d)^{-1}.
		\end{multline}
		Since $\epsilon(n,d)^{-1}<0$ and $(qb^{d-1})_{n-j-i-d+1}/(b^{n+d})_{n-j-i-d+1}\geq 0$, the third summand is nonnegative. Together with the lower bound from~\eqref{eq:bounds_eps_Hermpol_neven} we obtain
		\begin{align*}
			1-\epsilon_{i,j}>1-q^{-j-i-d+1}\frac{(q^n+q+1)(q^{n+d-1}+1)(q^{j+i}+1)}{(q^{n-d+1}+1)(q^{2n+j+i}+1)}.
		\end{align*}
		By \eqref{eq:ineqfraction1} we have
		\[
		\frac{q^{n+d-1}+1}{q^{n-d+1}+1}<q^{2d-2},
		\]
		which gives us
		\begin{align*}
			1-\epsilon_{i,j}
			>1-\frac{q^n+q+1}{q^{2n+j+i-d}}
			=1-\frac{1}{q^{n+j+i-d}}-\frac{1}{q^{2n+j+i-d-1}}-\frac{1}{q^{2n+j+i-d}}.
		\end{align*}
		Because of $n-d\geq 2$, $n\geq 4$, and $i+j\geq 1$, we deduce
		\[
		1-\epsilon_{i,j}>\frac{109}{128}.
		\]
		The second summand in~\eqref{eq:errorterm_Hermpol_neven_jiodd} is nonpositive since $\epsilon(n,d)^{-1}<0$. Together with the lower bound from~\eqref{eq:bounds_eps_Hermpol_neven}, we obtain
		\[
		1-\epsilon_{i,j}<1+\frac{(q+1)}{q^n}\, \frac{(qb^{d-1})_{n-j-i-d+1}}{(b^{n+d})_{n-j-i-d+1}}.
		\]
		Using~\eqref{eq:epsij_bounds_prod} and $i+j\leq n-d-1$ as well as $n\geq 4$ gives
		\[
		1-\epsilon_{i,j}<1+\frac{(q+1)}{q^n}q^{-(n-2)(n-j-i-d+1)}\leq 1+\frac{q+1}{q^{3n-4}}\leq \frac{259}{256},
		\]
		which completes the proof.
	\end{proof}
	
	We now show that the inner distribution is nonnegative.
	\begin{prop}\label{prop:innerdistrHermpol_deven_nonneg}
		For all even $d$ with $2\leq d\leq n$, all entries of the inner distribution $(A_i)$ given in Proposition~\ref{prop:innerdistrHermpol_deven} are nonnegative.
	\end{prop}
	\begin{proof}
		Let $(A_i)$ be given in the form of~\eqref{eq:innerdistr_Hermpoleven_rewritten}. First, assume that $d=n$. Then $n$ is even and we only have to show that $A_n\geq 0$. By using~\eqref{eq:sizeY_hermPol}, we have
		\begin{align*}
			A_n
			&=\frac{|Y|}{|X|}\,\frac{\pochb{b^{n+1}}{n}}{\pochb{q}{n}}\,\left( 1-\frac{\pochb{qb^{n-1}}{1}}{\pochb{b^{2n}}{1}}\epsilon(n,n)^{-1} \right)\\
			&=\frac{\pochb{q}{n-1}\pochb{b^{n+1}}{n}}{\pochb{b^{n+1}}{n-1}\pochb{q}{n}}\,\epsilon(n,n)
			-\frac{\pochb{q}{n-1} \pochb{b^{n+1}}{n} \pochb{qb^{n-1}}{1} }{\pochb{b^{n+1}}{n-1}\pochb{q}{n}\pochb{b^{2n}}{1}}\\
			&=-(q^n-1)\epsilon(n,n)-1.
		\end{align*}
		From Lemma~\ref{lem:eps_hermpol} we obtain
		\[
		A_n>\frac{q^n(q^n-1)}{q+1}-1>0.
		\]
		
		Assume now that $2\leq d<n$. Since $1-\epsilon_{i,j}>0$ for all $i,j$ by Lemma~\ref{lem:bounds_errorterm_Hermpol} and the sign of $\pochb{b^{n+1}}{n-j-i}/\pochb{q}{n-j-i}$ is $(-1)^{(n+1)(n-j-i)}$, the sign of $a_{i,j}$ is $(-1)^{\binom{j}{2}+ij+j}$.  Hence for all $i=0,1,\dots,n-d$, we have $\sign(a_{i,2j})=1$ if $j\geq 0$ is even and $\sign(a_{i,2j+1})=1$ if $j+i\equiv 1\pmod 2$. In all other cases, we have $\sign(a_{i,j})=-1$. Similarly as in the case of odd $d$, we will show that $(a_{i,j})$ satisfies~\eqref{eq:props_aij}. Observe that this will prove the nonnegativity of the inner distribution.
		
		Take $i\in\{0,1,\dots,n-d\}$. For all $j=0,1,\dots,n-d-i-2$, by using~\eqref{eq:frac_qbin}, we obtain
		\begin{align*}
			\frac{|a'_{i,j}|}{|a'_{i,j+2}|}
			=\frac{q^{\binom{j}{2}} \left\lvert\bbin{n-i}{j}\frac{\pochb{b^{n+1}}{n-j-i}}{\pochb{q}{n-j-i}}\right\rvert}
			{q^{\binom{j+2}{2}}\left\lvert\bbin{n-i}{j+2}\frac{\pochb{b^{n+1}}{n-j-i-2}}{\pochb{q}{n-j-i-2}}\right\rvert}
			\geq q^{2j+2i-2n+1}\,\frac{|(b^{2n-j-i-1}-1)(b^{2n-j-i}-1)|}{|(qb^{n-j-i-2}-1)(qb^{n-j-i-1}-1)|}.
		\end{align*}
		Applying~\eqref{eq:ineqfraction1} gives
		\begin{align*}
			\frac{|a'_{i,j}|}{|a'_{i,j+2}|}\geq q^{2j+2i-2n+1}\,\frac{(q^{2n-j-i-1}-1)(q^{2n-j-i}+1)}{(q^{n-j-i-1}+1)(q^{n-j-i}-1)}>q^{2j+2i}.
		\end{align*}
		
		Assume that $n\geq 4$. By using Lemma~\ref{lem:bounds_errorterm_Hermpol}, we deduce
		\[
		\frac{1-\epsilon_{i,j}}{1-\epsilon_{i,j+2}}>\begin{cases}
			\frac{31}{64}&\text{for $n$ and $i+j$ odd}\\
			\frac{61}{102} &\text{for $n$ and $i+j$ even}\\
			\frac{191}{256} &\text{for odd $n$ and even $i+j$}\\
			\frac{218}{259} &\text{for even $n$ and odd $i+j$}
		\end{cases}
		\]
		for all $i,j$ with $0\leq i+j\leq n-d-2$. For all $i,j$ with $1\leq i+j\leq n-d-2$, we thus obtain
		\[
		\frac{|a_{i,j}|}{|a_{i,j+2}|}
		=\frac{|a'_{i,j}|(1-\epsilon_{i,j})}{|a'_{i,j+2}|(1-\epsilon_{i,j+2})}
		>q^{2i+2j}\frac{1-\epsilon_{i,j}}{1-\epsilon_{i,j+2}}
		>1,
		\]
		as wanted. Let $i$ be even with $0\leq i\leq n-d$. In the case of $n\geq 4$, it remains to show that
		\[
		\frac{|a_{i,0}|}{|a_{i,1}|+|a_{i,2}|}\geq 1.
		\]
		We have
		\[
		\frac{|a_{i,0}|}{|a_{i,1}|+|a_{i,2}|}
		=\frac{\left\lvert \frac{(b^{n+1})_{n-i}}{(q)_{n-i}}\right\rvert(1-\epsilon_{i,0})}{\left\lvert \bbin{n-i}{1}\frac{(b^{n+1})_{n-i-1}}{(q)_{n-i-1}}\right\rvert(1-\epsilon_{i,1})+q\left\lvert \bbin{n-i}{2}\frac{(b^{n+1})_{n-i-2}}{(q)_{n-i-2}}\right\rvert(1-\epsilon_{i,2})}.
		\]
		It holds that
		\begin{align*}
			\left\lvert \frac{(b^{n+1})_{n-i-1}}{(q)_{n-i-1}}\right\rvert
			=\frac{(q^{n-i}+(-1)^n)}{(q^{2n-i}-1)}\left\lvert \frac{(b^{n+1})_{n-i}}{(q)_{n-i}}\right\rvert
		\end{align*}
		and
		\begin{align*}
			\left\lvert \frac{(b^{n+1})_{n-i-2}}{(q)_{n-i-2}}\right\rvert=\frac{(q^{n-i}+(-1)^n)(q^{n-i-1}-(-1)^n)}{(q^{2n-i}-1)(q^{2n-i-1}+1)}\left\lvert \frac{(b^{n+1})_{n-i}}{(q)_{n-i}}\right\rvert.
		\end{align*}
		Combining this with
		\begin{align*}
			\left\lvert\bbin{n-i}{1}\right\rvert&=\frac{q^{n-i}-(-1)^n}{q+1}\\
			\left\lvert\bbin{n-i}{2}\right\rvert&=\frac{(q^{n-i}-(-1)^n)(q^{n-i-1}+(-1)^n)}{(q+1)(q^2-1)}
		\end{align*}
		gives
		\begin{align*}
			\frac{|a_{i,0}|}{|a_{i,1}|+|a_{i,2}|}
			&=\frac{1-\epsilon_{i,0}}{\frac{(q^{2n-2i}-1)}{(q+1)(q^{2n-i}-1)}(1-\epsilon_{i,1})+\frac{q(q^{2n-2i}-1)(q^{2n-2i-2}-1)}{(q+1)(q^2-1)(q^{2n-i}-1)(q^{2n-i-1}+1)}(1-\epsilon_{i,2})}.
		\end{align*}
		Because of~\eqref{eq:ineqfraction}, we obtain
		\begin{align*}
			\frac{|a_{i,0}|}{|a_{i,1}|+|a_{i,2}|}
			\geq \frac{1-\epsilon_{i,0}}{\frac{1}{q^i(q+1)}\,(1-\epsilon_{i,1})+\frac{1}{q^{2i}(q+1)(q^2-1)}\,(1-\epsilon_{i,2})}.
		\end{align*}
		For $i\geq 2$, we find from Lemma~\ref{lem:bounds_errorterm_Hermpol} that
		\begin{align*}
			\frac{|a_{i,0}|}{|a_{i,1}|+|a_{i,2}|}
			\geq \frac{1-\epsilon_{i,0}}{\frac{1}{12}\,(1-\epsilon_{i,1})+\frac{1}{144}\,(1-\epsilon_{i,2})}
			\geq 1.
		\end{align*}
		For $i=0$, we find from Lemma~\ref{lem:bounds_errorterm_Hermpol} that
		\begin{align*}
			\frac{|a_{0,0}|}{|a_{0,1}|+|a_{0,2}|}
			\geq \frac{1-\epsilon_{0,0}}{\frac13\,(1-\epsilon_{0,1})+\frac19\,(1-\epsilon_{0,2})}.
		\end{align*}
		This shows that
		\begin{align*}
			\frac{|a_{i,0}|}{|a_{i,1}|+|a_{i,2}|}
			\geq 1
		\end{align*}
		for all $n\geq 4$ and even $i\geq 0$, as required.

		It remains to look at $n=3$ and $d=2$. Because of the sign of $a_{1,0}$, we immediately have $A_2=\frac{|Y|}{|X|}a_{1,0}\geq 0$. For the entry $A_3$, we only need to show $|a_{0,0}|/|a_{0,1}|\geq 1$. We have
		\[
		\frac{|a'_{0,0}|}{|a'_{0,1}|}=\frac{\left\lvert \frac{(q^4)_3}{(q)_3}\right\rvert}{\left\lvert \bbin{3}{1}\frac{(q^4)_2}{(q)_2}\right\rvert}=q+1\geq 3.
		\]
		We also require bounds on $1-\epsilon_{0,0}$ and $1-\epsilon_{0,1}$. Use the bound from~\eqref{eq:boundlower_epsHermpol} to obtain
		\begin{align*}
			1-\epsilon_{0,0}
			&=1-\frac{(q^2+1)}{(q^5+1)(q^3+1)}\,\epsilon(3,2)^{-1}\\
			&\geq 1-\frac{2(q^2+1)}{(q^5+1)(q^3+1)(q^2-1)}\\
			&\geq 1-\frac{2}{q^6}\\
			&\geq \frac{31}{32}.
		\end{align*}
		Applying the bound from~\eqref{eq:boundupper_epsHermpol} gives
		\begin{align*}
			1-\epsilon_{0,1}
			&=1+(1-\epsilon(3,2)^{-1})\,q\frac{(q^4-1)(q+1)}{(q^2-1)(q^5+1)}-\frac{(q^2+1)}{(q^5+1)}\,\epsilon(3,2)^{-1}\\
			&\leq 1+\left(1-\frac{q-1}{q^4-1}\right)q\,\frac{(q^2+1)(q+1)}{(q^5+1)}\\
			&=1+\frac{q^2(q^3-1)}{(q-1)(q^5+1)}\\
			&\leq 2.
		\end{align*}
		In summary, we conclude
		\[
		\frac{|a_{0,0}|}{|a_{0,1}|}
		=\frac{|a'_{0,0}|}{|a'_{0,1}|}\,\frac{(1-\epsilon_{0,0})}{(1-\epsilon_{0,1})}
		>1.
		\]
		This completes the proof.
	\end{proof}
	
	It remains to show that the dual distribution is nonnegative. This requires the inequality
	\begin{align}\label{eq:ineq_product}
		\prod_{i=1}^n \left(1+\frac{1}{q^i}\right)<\frac52
	\end{align}
	for all integers $n\geq 1$ and $q\geq 2$, which was proved in \cite[Lem.~3.6]{SchmidtWeissSteiner}.
	
	\begin{prop}\label{prop:dualdistrHermpol_deven_nonneg}
		For all even $d$ with $2\leq d\leq n$, all entries of the dual distribution $(A'_k)$ given in Proposition~\ref{prop:dualdistrHermpol_deven} are nonnegative.
	\end{prop}
	\begin{proof}
		Let $(A'_k)$ be given in Proposition~\ref{prop:dualdistrHermpol_deven}. We first show that $A'_1\geq 0$ for $2\leq d\leq n$. The sign of $(q)_{d-1}/(b^n)_{d-1}$ is $(-1)^n$ and from~\eqref{eq:dualdistr1_Hermpol}, we thus find
		\[
		A'_1=|X| q^{-d+1} \left\lvert\frac{(q)_{d-1}}{(b^n)_{d-1}}\right\rvert \frac{(q^n-(-1)^n)}{(q^{n-d+1}+(-1)^n)}((-1)^n-(-1)^n\epsilon(n,d)).
		\]
		From Lemma \ref{lem:eps_hermpol}, we see that $\epsilon(n,d)\geq 1$ for odd $n$ and $\epsilon(n,d)<0$ for even~$n$. Thus, we have $(-1)^n-(-1)^n\epsilon(n,d)\geq 0$ implying $A'_1\geq 0$, as required. 
		
		Observe that we can now consider $d\geq 4$ since for $d=2$, we only need to show $A'_1\geq 0$. Set $b_{k,j}=b'_{k,j} (1-\delta_{k,j})$ with $\delta_{k,j}$ as in~\eqref{eq:delta_Hermpol_deven} and
		\[
		b'_{k,j}=b^{\binom{j}{2}-nj} \frac{(b^{-n-1-k})_{n-k-j}}{(-b^{-n})_{n-k-j}}\bbin{n-k}{j}
		\]
		for all $k=0,1,\dots,d-3$ and $j=0,1,\dots,d-k-3$, so that
		\[
		A'_{n-k}=c_k\sum\limits_{j=0}^{d-k-3} b_{k,j}
		\]
		as in~\eqref{eq:dualdistr_Hermpol_deven}. Since the sign of $(qb^k)_{n-k}/(-b^{-n})_{n-k}$ is $(-1)^{\binom{n-k}{2}+(n-k)(k+1)}$, we obtain $c_k\geq 0$. Hence, we need to show that $\sum_{j=0}^{d-k-3} b_{k,j}\geq 0$. This will follow by proving that $(b_{k,j})$ satisfies~\eqref{eq:props_aij}. The proof is split into some intermediate results, starting with bounds on $|\delta_{k,j}|$.
		
		\textbf{Claim 1}. \textit{For all integers $q$, $n$, and $d$ with $q\geq 2$, $n\geq 8$, and $4\leq d\leq n$, we have 
			\begin{align*}
				|\delta_{k,j}|\leq\begin{cases}
					0.0032 & \text{if\, $0\leq k+j\leq d-4$},\\
					0.05 & \text{if\, $k+j=d-3$}.
				\end{cases}
			\end{align*}
			For all $q\geq 4$ and $n\geq 4$, we have
			\begin{align*}
				|\delta_{k,j}|\leq\begin{cases}
					0.014 &\text{if\, $0\leq k+j\leq d-4$},\\
					0.095 &\text{if\, $k+j=d-3$.}
				\end{cases}
		\end{align*}}
		\textit{Proof of Claim 1.} For $0\leq k+j\leq d-3$, we have
		\begin{multline*}
			|\delta_{k,j}|
			\leq |\epsilon(n,d)|\left\lvert \frac{(qb^{k+j})_{d-k-j-1}}{(b^{n+k+j+1})_{d-k-j-1}}\right\rvert
			\\
			+(1+|\epsilon(n,d)|)q^{k+j-d+1}\frac{(q^{n-k-j}+1)}{(q^{n-d+1}-1)}\left\lvert \frac{(qb^{k+j})_{d-k-j-1}}{(b^{n+k+j})_{d-k-j-1}}\right\rvert.
		\end{multline*}
		Using $|\epsilon(n,d)|\leq (q^{n+d-1}+1)/(q^{d-1}-1)$ from Lemma~\ref{lem:eps_hermpol} gives
		\begin{align}\label{eq:deltabound_HermPoleven}
			|\delta_{k,j}|
			\leq \frac{(q^{n+d-1}+1)}{(q^{d-1}-1)}\left\lvert \frac{(qb^{k+j})_{d-k-j-1}}{(b^{n+k+j+1})_{d-k-j-1}}\right\rvert
			+q^{k+j}\frac{(q^n+1)(q^{n-k-j}+1)}{(q^{d-1}-1)(q^{n-d+1}-1)}\left\lvert \frac{(qb^{k+j})_{d-k-j-1}}{(b^{n+k+j})_{d-k-j-1}}\right\rvert.
		\end{align}
		For $m\in\{n,n+1\}$, by using \eqref{eq:ineq_product_minus} and \eqref{eq:ineq_product}, we find
		\[
		\left\lvert \frac{(qb^{k+j})_{d-k-j-1}}{(b^{m+k+j})_{d-k-j-1}}\right\rvert\leq \prod_{\ell=0}^{d-k-j-2}\frac{q^{k+j+1+\ell}+1}{q^{m+k+j+\ell}-1}\leq 10 q^{-(m-1)(d-k-j-1)}.
		\]
		This gives us
		\begin{align*}
			|\delta_{k,j}|
			\leq 10 q^{-n(d-k-j-1)}\,\frac{q^{n+d-1}+1}{q^{d-1}-1}
			+10q^{k+j-(n-1)(d-k-j-1)}\frac{(q^n+1)(q^{n-k-j}+1)}{(q^{d-1}-1)(q^{n-d+1}-1)},
		\end{align*}
		which becomes
		\begin{align*}
			|\delta_{k,j}|
			\leq 10q^{-n(d-k-j-2)}\,\frac{1+q^{-n-d+1}}{1-q^{-d+1}}
			+10q^{-(n-1)(d-k-j-2)+1}\frac{(1+q^{-n})(1+q^{-n+k+j})}{(1-q^{-d+1})(1-q^{-n+d-1})}.
		\end{align*}
		For all $n\geq 8$, $q\geq 2$, and $0\leq k+j\leq d-4$, we obtain $|\delta_{k,j}|\leq 0.0032$. For all $n\geq 4$, $q\geq 4$, and $0\leq k+j\leq d-4$, we obtain $|\delta_{k,j}|\leq 0.014$. In the case of $k+j=d-3$, we find from~\eqref{eq:deltabound_HermPoleven} that
		\begin{align*}
			|\delta_{k,j}|
			&\leq \frac{q^{d-2}+1}{q^{n+d-2}-1}+q^{d-3}\frac{(q^n+1)(q^{n-d+3}+1)(q^{d-2}+1)}{(q^{n-d+1}-1)(q^{n+d-3}-1)(q^{n+d-2}+1)}\\[2ex]
			&\leq q^{-n}\,\frac{1+q^{-d+2}}{1-q^{-n-d+2}}+q^{-n+2}\frac{(1+q^{-n})(1+q^{-n+d-3})(1+q^{-d+2})}{(1-q^{-n+d-1})(1-q^{-n-d+3})}
		\end{align*}
		For all $n\geq 8$ and $q\geq 2$, we obtain $|\delta_{k,j}|\leq 0.05$, and for all $n\geq 4$ and $q\geq 4$, we obtain $|\delta_{k,j}|\leq 0.095$. This completes the proof of Claim 1.
		
		\textbf{Claim 2.} \textit{For all integers $q\geq 2$, $n\geq 4$, $k=0,1,\dots,d-3$, and $j=0,1,\dots,d-k-5$, we have
			\[
			\frac{|b'_{k,j}|}{|b'_{k,j+2}|}\geq q^{2j+2k+1}\,\frac{1-q^{-2k-j-2}}{1+q^{-k-j-1}}.
			\]}
		\textit{Proof of Claim 2.} Using the definition of the $q$-Pochhammer symbol, we see that
		\[
		\frac{|b'_{k,j}|}{|b'_{k,j+2}|}=\left\lvert \frac{b^{\binom{j}{2}-nj}\bbin{n-k}{j}(1-b^{-2k-j-3})(1-b^{-2k-j-2})}{b^{\binom{j+2}{2}-n(j+2)}\bbin{n-k}{j+2}(1+b^{-k-j-2})(1+b^{-k-j-1})}\right\rvert.
		\]
		Applying~\eqref{eq:frac_qbin} gives
		\[
		\frac{|b'_{k,j}|}{|b'_{k,j+2}|}
		\geq q^{2j+2k+1}\frac{(1+q^{-2k-j-3})(1-q^{-2k-j-2})}{(1-q^{-k-j-2})(1+q^{-k-j-1})}
		\geq q^{2j+2k+1}\,\frac{1-q^{-2k-j-2}}{1+q^{-k-j-1}},
		\]
		as wanted.
		
		\textbf{Claim 3.} \textit{For all $q\geq 2$ and $n\geq 8$, we have $|b_{k,0}|/|b_{k,1}|\geq 2$ for all $k\geq 1$. For all $q\geq 4$ and $n\geq 4$, we have $|b_{k,0}|/|b_{k,1}|\geq 2$ for all $k\geq 0$.}\\
		\textit{Proof of Claim 3.} We have
		\[
		\frac{|b'_{k,0}|}{|b'_{k,1}|}
		=\left\lvert \frac{1-b^{-2k-2}}{b^{-n}(1+b^{-k-1})\bbin{n-k}{1}}\right\rvert
		\geq \frac{(1-q^{-2k-2})(q+1)}{q^{-n}(1+q^{-k-1})(q^{n-k}+1)}.
		\]
		For $q\geq 2$, $n\geq 8$, and $k\geq 1$, we have
		\[
		\frac{|b'_{k,0}|}{|b'_{k,1}|}\geq \frac{(1-q^{-4})(q+1)}{q^{-n}(1+q^{-2})(q^{n-1}+1)}\geq \frac{129}{43}.
		\]
		For $q\geq 4$, $n\geq4$, and $k\geq 0$, we obtain
		\[
		\frac{|b'_{k,0}|}{|b'_{k,1}|}\geq \frac{(1-q^{-2})(q+1)}{q^{-n}(1+q^{-2})(q^n+1)}\geq \frac{960}{257}.
		\]
		Together with Claim 1, this proves Claim 3.
		
		\textbf{Claim 4.} \textit{The sign of $b_{k,j}$ is $(-1)^{\binom{j}{2}+kj+j}$ for all $j,k\geq 0$, if $q\geq 2$ and $n\geq 8$, or if $q\geq 4$ and $n\geq 4$.}\\
		\textit{Proof of Claim 4.} This follows from Claim~1, which gives $1-\delta_{k,j}\geq 0$, and from
		\[
		\frac{(b^{-n-1-k})_{n-k-j}}{ (-b^{-n})_{n-k-j}}\geq 0.
		\]
		
		We now combine the results from the previous claims to show that $b_{k,j}$ satisfies~\eqref{eq:props_aij}. Assume that $q\geq 2$ and $n\geq 8$, or that $q \geq 4$ and $n\geq 4$. Because of Claim~4, we then have
		\[
		\frac{b_{k,2j}}{|b_{k,2j+2}|}=\frac{|b_{k,2j}|}{|b_{k,2j+2}|}\quad\text{for all $k\geq 0$ and even $j\geq 0$}
		\]
		and
		\[
		\frac{b_{k,2j+1}}{|b_{k,2j+3}|}=\frac{|b_{k,2j+1}|}{|b_{k,2j+3}|}\quad\text{for all $k,j\geq 0$ with $k+j\equiv 1\pmod 2$.}
		\]
		We can thus look at $|b_{k,j}|/|b_{k,j+2}|$. Claim~1 and~2 imply
		\begin{align*}
			\frac{|b_{k,j}|}{|b_{k,j+2}|}\geq \frac{|b'_{k,j}|}{|b'_{k,j+2}|}\,\frac{(1-|\delta_{k,j}|)}{(1+|\delta_{k,j+2}|)}\geq 2
		\end{align*}
		for all $k=0,1,\dots,d-3$, $j=0,1,\dots,d-k-5$ if $q\geq 4$, $n\geq 4$, or for all $k=0,1,\dots,d-3$, $j=0,1,\dots,d-k-5$ with $k+j\neq 0$ if $q\geq 2$, $n\geq 8$, which also gives
		\[
		\frac{|b_{k,0}|}{|b_{k,1}|+|b_{k,2}|}\geq 1
		\]
		in the respective cases by using Claim 3.
		
		It remains to look at $q\geq 2$, $n\geq 8$, and $k=j=0$, that is
		\[
		\frac{|b_{0,0}|}{|b_{0,1}|+|b_{0,2}|}.
		\]
		Using the definition of the $q$-Pochhammer symbol, we see that
		\begin{align*}
			\frac{|b_{0,0}|}{|b_{0,1}|+|b_{0,2}|}
			&=\frac{1-\delta_{0,0}}{q^{-n}\frac{(1-q^{-1})}{(1-q^{-2})}\left\lvert\bbin{n}{1}\right\rvert (1-\delta_{0,1})+q^{1-2n}\frac{(1-q^{-1})(1+q^{-2})}{(1-q^{-2})(1+q^{-3})}\left\lvert\bbin{n}{2}\right\rvert (1-\delta_{0,2})}\\
			&\geq \frac{1-\delta_{0,0}}{q^{-n}\frac{1}{(1-q^{-2})}\left\lvert\bbin{n}{1}\right\rvert (1-\delta_{0,1})+q^{1-2n}\frac{(1+q^{-2})}{(1-q^{-2})}\left\lvert\bbin{n}{2}\right\rvert (1-\delta_{0,2})}.
		\end{align*}
		By applying
		\[
		\left\lvert\bbin{n}{1}\right\rvert\leq \frac{q^n+1}{q+1}\qandq\left\lvert \bbin{n}{2}\right\rvert\leq \frac{(q^n-1)(q^{n-1}+1)}{(q+1)(q^2-1)},
		\]
		we find
		\[
		\frac{|b_{0,0}|}{|b_{0,1}|+|b_{0,2}|}\geq \frac{1-\delta_{0,0}}{\frac{(1+q^{-n})}{(1-q^{-2})(q+1)} (1-\delta_{0,1})+\frac{(1+q^{-2})(1+q^{-n+1})}{(1-q^{-2})(q+1)(q^2-1)}(1-\delta_{0,2})}.
		\]
		By using $q\geq 2$ and Claim~1, this becomes
		\[
		\frac{|b_{0,0}|}{|b_{0,1}|+|b_{0,2}|}\geq\frac{1-|\delta_{0,0}|}{\frac{257}{576}(1+|\delta_{0,1}|)+\frac{215}{1152}(1+|\delta_{2,2}|)}\geq1,
		\]
		as wanted.
		
		For the remaining cases with $q\in\{2,3\}$, $n\in\{4,5,6,7\}$, and $d\in\{4,\dots,n\}$, the dual distribution is nonnegative as well, which was checked with a computer. This concludes the proof.
	\end{proof}
	
	Combining Proposition~\ref{prop:innerdistrHermpol_deven_nonneg} and~\ref{prop:dualdistrHermpol_deven_nonneg} proves Proposition~\ref{prop:fsol_primal_hermpol_deven}.
	
	\renewcommand{\UrlFont}{\small\ttfamily}
	\printbibliography

\end{document}